\documentclass[authoryear]{elsarticle}

\usepackage[utf8]{inputenc}
\usepackage{a4wide, amsmath, amssymb}
\usepackage{algorithmic}
\usepackage[active]{srcltx}
\usepackage{tikz}
\usetikzlibrary{shapes,arrows,backgrounds,positioning,plotmarks,calc,patterns}
\usepackage{picins}
\usepackage{float}
\floatstyle{ruled}
\usepackage{enumerate}

\newfloat{Decompose}{tbp}{lop}[section]

\usepackage{maple2e}
\def\emptyline{\vspace{6pt}}

\usepackage[T1]{fontenc}

\newtheorem{theorem}{Theorem}[section]
\newtheorem{myalgorithm}[theorem]{Algorithm}{\bfseries}{\rmfamily}
\newtheorem{mydefinition}[theorem]{Definition}{\bfseries}{\rmfamily}
\newtheorem{myremark}[theorem]{Remark}{\bfseries}{\rmfamily}
\newtheorem{mycorollary}[theorem]{Corollary}{\bfseries}{\rmfamily}
\newtheorem{mylemma}[theorem]{Lemma}{\bfseries}{\rmfamily}
\newtheorem{defrem}[theorem]{Definition and Remark}{\bfseries}{\rmfamily}
\newtheorem{myexample}[theorem]{Example}{\bfseries}{\itshape}
\newproof{proof}{Proof}

\newcommand{\ld}{\operatorname{ld}}
\newcommand{\ini}{\operatorname{init}}
\newcommand{\mdeg}{\operatorname{mdeg}}
\newcommand{\prem}{\mathsf{prem}}
\newcommand{\pquo}{\mathsf{pquo}}
\newcommand{\sepa}{\operatorname{sep}}
\newcommand{\der}[1]{\ensuremath{_{\mbox{\textnormal{\tiny #1}}}}}

\newcommand{\res}{\operatorname{res}}

\newcommand{\lcm}{\ensuremath{\operatorname{lcm}}}

\newcommand{\C}{\ensuremath{\mathbb{C}}}

\newcommand{\R}{\ensuremath{\mathbb{R}}}
\newcommand{\Z}{\ensuremath{\mathbb{Z}}}

\newcommand{\sol}{\ensuremath{\mathfrak{Sol}}}

\newcommand{\picturesize}{0.8cm}

\newcommand{\prs}{\ensuremath{\operatorname{PRS}}}
\newcommand{\ma}{\ensuremath{\mathbf{a}}}
\newcommand{\pma}[1][]{\ensuremath{\phi_{\ifthenelse{\equal{#1}{}}{}{#1,}\mathbf{a}}}}
\def\osim{
  \mathrel{
    \mathchoice {
      \vcenter{
        \offinterlineskip
        \halign{
          \hfil $\displaystyle##$\hfil\cr\not\succ\cr\noalign{\vskip-3pt}\not\prec\cr
        }
      }
    }
    {
      \vcenter{
        \offinterlineskip
        \halign{
          \hfil$\textstyle##$\hfil\cr\not\succ\cr\noalign{\vskip-3pt}\not\prec\cr
        }
      }
    }
    {
      \vcenter{
        \offinterlineskip
        \halign{
          \hfil$\scriptstyle##$\hfil\cr\not\succ\cr\noalign{\vskip-2.4pt}\not\prec\cr
        }
      }
    }
    {
      \vcenter{
        \offinterlineskip
        \halign{
          \hfil$\scriptscriptstyle##$\hfil\cr\not\succ\cr\noalign{\vskip-0.9pt}\not\prec\cr
        }
      }
    }
  }
}

\usepackage[pdfauthor   = {B\"achler and Gerdt and Lange-Hegermann and Robertz},
            pdftitle    = {Algorithmic Thomas Decomposition of Algebraic and Differential Systems},
            pdfsubject  = {},
            pdfkeywords = {},
            bookmarks=true,
            bookmarksopen=true,
            pagebackref=true,
            hyperindex=true,
            colorlinks=true,
            linkcolor=blue,
            citecolor=green,
            filecolor=magenta,
            urlcolor=cyan,
            ]{hyperref}

\title{Algorithmic Thomas Decomposition of Algebraic and Differential Systems}

\begin{document}

\begin{frontmatter}

\author[lbfm]{Thomas B\"achler}
\ead{thomas@momo.math.rwth-aachen.de}
\author[lbfm]{Vladimir Gerdt\fnref{fn1}}
\ead{gerdt@jinr.ru}
\author[lbfm]{Markus Lange-Hegermann}
\ead{markus@momo.math.rwth-aachen.de}
\author[lbfm]{Daniel Robertz}
\ead{daniel@momo.math.rwth-aachen.de}

\fntext[fn1]{Permanent address: Joint Institute for Nuclear Research, Dubna, Russia.}
\address[lbfm]{Lehrstuhl B f\"ur Mathematik, RWTH-Aachen University, Templergraben 64, 52062 Aachen, Germany}

\begin{abstract}
 In this paper, we consider systems of algebraic and non-linear partial differential equations and inequations.
 We decompose these systems into so-called simple subsystems and thereby partition the set of solutions.
 For algebraic systems, simplicity means triangularity, square-freeness and non-vanishing initials.
 Differential simplicity extends algebraic simplicity with involutivity.
 We build upon the constructive ideas of J. M. Thomas and develop them into a new algorithm for disjoint decomposition.
 The given paper is a revised version of \cite{thomasalg_casc} and includes the proofs of correctness and termination of our decomposition algorithm.
 In addition, we illustrate the algorithm with further instructive examples and describe its Maple implementation together with an experimental comparison to some other triangular decomposition algorithms.
\end{abstract}


\end{frontmatter}

\section{Introduction}\label{introduction}

Nowadays, triangular decomposition algorithms, which go back to the characteristic set method by \cite{Ritt} and \cite{Wu}, have become powerful tools for investigating and solving systems of multivariate polynomial equations.
In many cases these methods are computationally more efficient than those based on construction of \textsc{Gr\"obner} bases.
For an overview over triangular decomposition methods for polynomial and differential-polynomial systems we refer to the tutorial papers by \cite{Hubert1,Hubert2} and to the bibliographical references therein.

Among numerous triangular decompositions the \textsc{Thomas} one stands by itself.
It was suggested by the American mathematician \cite{Tho1,Tho2} and decomposes a finite system of polynomial equations and/or inequations into finitely many triangular subsystems, which he called \emph{simple}.
The \textsc{Thomas} decomposition splits a given quasi-affine variety into a finite number of quasi-affine varieties defined by simple systems.
Unlike other decomposition algorithms, the \textsc{Thomas} decomposition always yields a \emph{disjoint} decomposition of the solution set.

Wang was the first to design and implement an algorithm that constructs the \textsc{Thomas} decomposition (cf.\ \cite{wang_simple,WangMethods,WangLi}).
For polynomial systems he implemented his algorithm in \textsc{Maple} (cf.\ \cite{WangPractice}) as part of the software package $\epsilon$\textsf{psilon} (cf.\ \cite{epsilon}), which also contains implementations of a number of other triangular decomposition algorithms.
\cite{delliere_wang_dyn} has shown that the ``dynamic constructible closure'' introduced in the thesis by \cite{gomez-diaz_constructible_closure} can be modeled using simple systems.
Nonetheless, according to the remark after \cite[Thm.~5.2]{delliere_wang_dyn}, simple systems are more general.

Every simple system is a regular system and its equations form a regular chain.
The \textsf{RegularChains} package (cf.~\cite{regularchains}) includes procedures to decompose the solution set of the input by means of regular chains (if the input only consists of equations) or regular systems.
However, the \textsc{Thomas} decomposition differs noticeably from this decomposition, since the \textsc{Thomas} decomposition is finer and demands disjointness of the solution set.
For a detailed description of algorithms related to regular chains, we refer the reader to \cite{MazaMEGA2000}.

The disjointness of the \textsc{Thomas} decomposition combined with the structural properties of simple systems provide a useful platform for counting solutions of polynomial systems.
In fact, the \textsc{Thomas} decomposition is the only known method to compute the {\em counting polynomial} introduced by \cite{PleskenCounting}.
We refer to \textsection\ref{sect_compare_simple_and_regular} for details on this structure, counting and their applications.

During his research on triangular decomposition, Thomas was motivated by the \textsc{Riquier}-\textsc{Janet} theory (cf.\ \cite{Riquier,Janet}), extending it to \emph{non-linear systems of partial differential equations}.
For this purpose he developed a theory of (\textsc{Thomas}) monomials, which generate an involutive monomial division nowadays called \textsc{Thomas} division (cf.\ \cite{GB1}).
He gave a recipe for decomposing a non-linear differential system into algebraically simple and passive subsystems (cf.\ \cite{Tho1}).
A modified version of the differential \textsc{Thomas} decomposition was considered by \cite{GerdtSimple} with its link to the theory of involutive bases (cf.\  \cite{GB1,GerI,Ger3,Seiler}).
In this decomposition, the output systems are \textsc{Janet}-involutive in accordance to the involutivity criterion from \cite{GerdtSimple} and hence they are coherent.
For a linear differential system it is a \textsc{Janet} basis of the corresponding differential ideal, as computed by the \textsc{Maple} package \textsf{Janet} (cf.\ \cite{JanetPackage}).

The differential \textsc{Thomas} decomposition differs from that computed by the \textsc{Rosenfeld}-\textsc{Gr\"ob\-ner} algorithm (cf.\ \cite{BLOP,BLOP2}).
The latter decomposition forms a basis of the \textsf{diffalg}, \textsf{DifferentialAlgebra} and \textsf{BLAD} packages (cf.\ \cite{diffalg,blad}).
Experimentally, we found that these three packages are optimized and well-suited for ordinary differential equations.
Furthermore, $\epsilon$\textsf{psilon} also allows to treat ordinary differential systems.
\cite{Bouziane} mentions another implementation not available to the authors.
However, all these methods give a zero decomposition, which, unlike the \textsc{Thomas} decomposition, is not necessarily disjoint.

In the given paper we present a new algorithmic version of the \textsc{Thomas} decomposition for polynomial and (ordinary and partial) differential systems.
In this unified algorithm, only two changes to the algebraic version are necessary to adapt it for the treatment of differential systems.
We briefly describe our implementation of this algorithm in \textsc{Maple}.

This paper is organized as follows.
In \textsection\ref{algebraic}, we present the algebraic part of our algorithm for the \textsc{Thomas} decomposition with its main objects defined in \textsection\ref{algebraic_definition_notation}.
In \textsection\ref{algebraic_algorithms}, we describe the main algorithm and its subalgorithms, followed by the correctness and termination proof.
Decomposition of differential systems is considered in \textsection\ref{section_differential}.
Here, we briefly introduce some basic notions and concepts from differential algebra (\textsection\ref{differential_preliminaries}) and from the theory of involutive bases specific to \textsc{Janet} division (\textsection\ref{Janet}).
In \textsection\ref{differential reduction}, we present our version of the differential pseudo reduction and, building upon it, the definition of differential simple systems.
Subsection \textsection\ref{differential algorithm} contains a description of the differential \textsc{Thomas} decomposition algorithm and the proof of its correctness and termination.
Some implementation issues are discussed in \textsection\ref{implementation}, followed by a comparison of our implementation to some other implementations of triangular decompositions with the help of benchmarks.

\section{Algebraic Systems}\label{algebraic}

This section introduces the concepts of simple systems and the \textsc{Thomas} decomposition for algebraic systems.
These concepts are based on properties of the set of solutions of a system.
We conclude the section with an algorithm for constructing a \textsc{Thomas} decomposition.

\begin{myexample}\label{ex_elliptic_curve} We give an easy example of a \textsc{Thomas} decomposition.
Consider the equation\\
\[
  p=\underline{x}^3 + (3y+1)x^2 + (3y^2+2y)x + y^3 = 0\enspace.
\]
\begin{tabular}{cc}  
  \begin{minipage}{0.55\textwidth}
  A {\sc Thomas} decomposition of $\{p=0\}$ is given by:
  \[\begin{array}{rl}
    S_1 := \{ & \underline{x}^3 + (3y+1)x^2 + (3y^2+2y)x + y^3 = 0, \\
              & 27\underline{y}^3 - 4y \neq 0 \} \\
    S_2 := \{ & 6\underline{x}^2 + (-27y^2 + 12y + 6)x - 3y^2 + 2y = 0, \\
              & 27\underline{y}^3 - 4y = 0 \}
  \end{array}\]
  \end{minipage}
  
  &
  
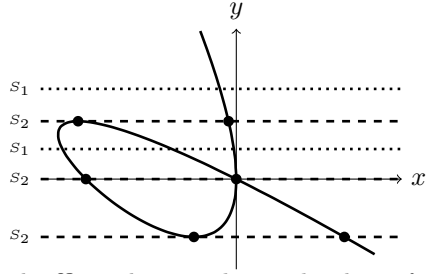
  \begin{minipage}{0.39\textwidth}
    \begin{tikzpicture}[domain=-1.4:1.4, smooth, x=2cm, y=2cm]
      \draw[line width=1pt, style=dashed] (-1.3,0.3849) node[left] {\tiny{$S_2$}} -- (1.1,0.3849);
      \draw[line width=1pt, style=dashed] (-1.3,0) node[left] {\tiny{$S_2$}} -- (1.1,0);
      \draw[line width=1pt, style=dashed] (-1.3,-0.3849) node[left] {\tiny{$S_2$}} -- (1.1,-0.3849);
      \draw[line width=1pt, style=dotted] (-1.3,0.2) node[left] {\tiny{$S_1$}} -- (1.1,0.2);
      \draw[line width=1pt, style=dotted] (-1.3,.6) node[left] {\tiny{$S_1$}} -- (1.1,.6);

      \draw[->] (-1.3,0) -- (1.1,0) node[right] {$x$};
      \draw[->] (0,-.6) -- (0,1) node[above] {$y$};
      \draw[line width=1pt] plot[] coordinates{ (-0.239629,0.984848) (-0.227273,0.952257) (-0.206107,0.893939) (-0.193948,0.860615) (-0.174083,0.803030) (-0.161924,0.767984) (-0.143747,0.712121) (-0.132560,0.678015) (-0.115310,0.621212) (-0.105958,0.590806) (-0.088996,0.530303) (-0.075758,0.481933) (-0.065026,0.439394) (-0.051379,0.384713) (-0.043553,0.348485) (-0.034607,0.307335) (-0.026093,0.257576) (-0.015791,0.197609) (-0.011931,0.166667) (-0.006399,0.127611) (-0.004241,0.075758) (0,0) (-0.005818,-0.106061) (-0.012236,-0.139279) (-0.017785,-0.166667) (-0.031959,-0.210465) (-0.045455,-0.242733) (-0.061588,-0.271745) (-0.075758,-0.292439) (-0.102481,-0.321761) (-0.136364,-0.347609) (-0.156728,-0.358423) (-0.196970,-0.373606) (-0.225728,-0.380333) (-0.257576,-0.383768) (-0.313458,-0.383512) (-0.348485,-0.379595) (-0.384766,-0.372810) (-0.439394,-0.358175) (-0.470348,-0.348485) (-0.530303,-0.324842) (-0.573198,-0.305590) (-0.621212,-0.281176) (-0.663814,-0.257576) (-0.712121,-0.228670) (-0.760233,-0.196970) (-0.803030,-0.167216) (-0.843861,-0.136364) (-0.893939,-0.096036) (-0.952348,-0.045455) (-0.984848,-0.014907) (-1.043257,0.045455) (-1.075758,0.083324) (-1.116588,0.136364) (-1.138334,0.168637) (-1.168885,0.227273) (-1.182775,0.273684) (-1.182023,0.318182) (-1.167040,0.348858) (-1.136364,0.369752) (-1.107041,0.379768) (-1.075758,0.383707) (-1.019946,0.383582) (-0.984848,0.380108) (-0.950633,0.374876) (-0.893939,0.363098) (-0.834698,0.348485) (-0.803030,0.339221) (-0.746003,0.321761) (-0.712121,0.310387) (-0.653616,0.289979) (-0.621212,0.277854) (-0.567659,0.257576) (-0.530303,0.242484) (-0.492270,0.227273) (-0.439394,0.204882) (-0.399889,0.187768) (-0.348485,0.165408) (-0.290794,0.139279) (-0.257576,0.123985) (-0.217940,0.106061) (-0.166667,0.081373) (-0.117482,0.056876) (0,0) (0.106061,-0.053826) (0.147274,-0.075758) (0.196970,-0.100842) (0.248823,-0.127611) (0.287879,-0.148773) (0.321402,-0.166667) (0.378788,-0.197292) (0.433739,-0.227273) (0.469697,-0.246790) (0.512231,-0.269807) (0.560606,-0.296808) (0.598881,-0.318182) (0.651515,-0.347372) (0.707157,-0.378788) (0.742424,-0.398705) (0.788002,-0.424365) (0.833333,-0.450482) (0.866887,-0.469697) (0.919443,-0.500000) };

      \draw plot[only marks,mark=*,mark size=1.5,mark options={line width=1pt}] coordinates{(0,0) (-1,0) (-1.051566,0.384900) (-0.051566,0.384900) (-0.281766,-0.384900) (0.718233,-0.384900)};
    \end{tikzpicture}
  \end{minipage}
\end{tabular}

The picture shows the solutions of $\{p=0\}$ in the real affine plane.
The cardinality of the fibers of the projection onto the $y$-component depends on $y$.
However, if we consider all solutions in the complex plane, this cardinality is constant within each system, i.e., $3$ and $2$ in $S_1$ and $S_2$, respectively.
This property is formalized in the definition of simple systems.

\end{myexample}

\subsection{Preliminaries}\label{algebraic_definition_notation}

Let $F$ be a computable field of characteristic 0 and $R:=F[x_1,\dots,x_n]$ be the polynomial ring in $n$ variables.
A total order $<$ on $\{1,x_1,\dots,x_n\}$ with $1<x_i$ for all $i$ is called a \textbf{ranking}.
The indeterminate $x$ is called \textbf{leader} of $p \in R$ if $x$ is the $<$-largest variable occurring in $p$.\footnote{In the context of triangular decompositions, the leader is usually called \textbf{main variable}. The term leader is used in \cite{Tho1} and has later been adopted in differential algebra.}
In this case we write $\ld(p)=x$.
If $p \in F$, we define $\ld(p)=1$.
The degree of $p$ in $\ld(p)$ is called \textbf{main degree} of $p$ ($\mdeg(p)$) and the leading coefficient $\ini(p) \in F[\ y\ |\ y<\ld(p)\ ]$ of $\ld(p)^{\mdeg(p)}$ in $p$ is called \textbf{initial} of $p$.

For $\ma \in \overline{F}^n$, where $\overline{F}$ denotes the algebraic closure of $F$, define the following evaluation homomorphisms:
\[\pma: F[x_1,\dots,x_n] \to \overline{F}: x_i \mapsto a_i\]\label{substitution_homomorphism}
\[\pma[<x_k]: F[x_1,\ldots,x_n] \to \overline{F}[x_k,\ldots,x_n]: \left\{\begin{array}{ll} x_i \mapsto a_i,&i<k \\x_i \mapsto x_i,&\mbox{otherwise}\end{array}\right.\]\label{part_substitution_homomorphism}

Given a polynomial $p \in R$, the symbols $p_{=}$ and $p_{\neq}$ denote the equation $p=0$ and inequation $p\neq0$, respectively.
A finite set of equations and inequations is called an \textbf{(algebraic) system} over $R$.
Abusing notation, we sometimes treat $p_{=}$ or $p_{\neq}$ as the underlying polynomial $p$.
A \textbf{solution} of $p_=$ or $p_{\not=}$ is a tuple $\ma \in \overline{F}^n$ with $\pma(p)=0$ or $\pma(p)\not=0$, respectively.
We call $\ma \in \overline{F}^n$ a solution of a system $S$, if it is a solution of each element in $S$.
The set of all solutions of $S$ is denoted by $\sol(S)$.

The subsets of all equations $p_{=} \in S$ and all inequations $p_{\neq}\in S$ are denoted by $S^{=}$ and $S^{\neq}$, respectively.
Define $S_x := \{ p \in S\ |\ \ld(p)=x \}$.
In a situation where it is clear that $|S_x|=1$, we also write $S_x$ to denote the unique element of $S_x$.
The subset $S_{<x} := \{ p \in S\ |\ \ld(p)<x \}$ is a system over $F[\ y\ |\ y<x\ ]$.

The {\sc Thomas} approach uses the homomorphisms $\pma[<x]$ to treat each polynomial $p\in S_x$ as the family of \emph{univariate} polynomials $\pma[<x](p) \in \overline{F}[x]$ for $\ma \in \sol(S_{<x})$.
This idea forms the basis of our central object, the \textbf{simple system}:

\begin{mydefinition}[Simple Systems]\label{simple_system} Let $S$ be a system.
 \begin{enumerate}
  \item $S$ is \textbf{triangular} if $|S_{x_i}|\leq 1\ \forall\ 1\leq i \leq n$ and $S \cap \{ c_=, c_{\neq}\mid c \in F\}=\emptyset$.
  \item \label{simple_nonzeroinitials} $S$ has \textbf{non-vanishing initials} if $\pma(\ini(p)) \not=0\ \forall\ \ma \in \sol(S_{<x_i})$ and $p \in S_{x_i}$ for $1 \leq i \leq n$.
  \item \label{simple_squarefree} $S$ is \textbf{square-free}\footnote{Square-freeness has an important side-effect in the differential case. A square-free polynomial and its separant have no common roots. Thus, the separants do not vanish on solutions of the lower-ranking subsystems.} if the univariate polynomial $\pma[<x_i](p) \in \overline{F}[x_i]$ is square-free $\forall\ma \in \sol(S_{<x_i})$ and $p \in S_{x_i}$ for  $1\leq i\leq n$.
  \item $S$ is called \textbf{simple} if it is \emph{triangular}, has \emph{non-vanishing initials} and is \emph{square-free}.
 \end{enumerate}
\end{mydefinition}

Properties (\ref{simple_nonzeroinitials}) and (\ref{simple_squarefree}) are characterized via solutions of lower-ranking equations and inequations.
However, the {\sc Thomas} decomposition algorithm does not calculate roots of polynomials.
Instead, it uses polynomial equations and inequations to \emph{partition} the set of solutions of the lower-ranking system to ensure the above properties.

\begin{myremark}\label{exist_sol}
Every simple system has a solution.
In particular, if $\mathbf{b} \in \sol(S_{<x})$ and $S_x$ is not empty, then $\phi_{<x,\mathbf{b}}(S_x)$ is a univariate polynomial with exactly $\mdeg(S_x)$ \emph{distinct} roots.
When $S_x$ is an equation, each solution $\mathbf{b} \in \sol(S_{<x})$ extends to a solution $(\mathbf{b},a) \in \sol(S_{\le x})$ with $\mdeg(S_x)$ possible choices $a \in \overline{F}$.
Otherwise, all but finitely many $a \in \overline{F}$ yield a solution $(\mathbf{b},a) \in \sol(S_{\le x})$, because an inequation $S_x$ excludes $\mdeg(S_x)$ different $a$ and $S_x=\emptyset$ imposes no restriction on $a$.

Conversely, if $(a_1,\ldots,a_n) \in \sol(S)$ where $S$ is a system over $F[x_1,\ldots,x_n]$ with $x_1 < \ldots < x_n$, then $(a_1,\ldots,a_i) \in \sol(S_{\le x_i})$.
\end{myremark}

To transform a system into a simple system, it is in general necessary to partition the set of solutions.
This leads to a so-called decomposition into simple systems.

\begin{mydefinition}
 A family $(S_i)_{i=1}^m$ is called \textbf{decomposition} of $S$ if $\sol(S)=\bigcup_{i=1}^m \sol(S_i)$.
 A decomposition is called \textbf{disjoint} if $\sol(S_i) \cap \sol(S_j) = \emptyset\ \forall\ i\neq j$.
 A \emph{disjoint} decomposition of a system into \emph{simple systems} is called \textbf{(algebraic) \textsc{Thomas} decomposition}.
\end{mydefinition}

For any algebraic system $S$, there exists a {\sc Thomas} decomposition (cf.\ \cite{Tho1,Tho2}, \cite{wang_simple}).
The algorithm presented in the following section provides another proof of this fact.

\begin{myexample}\label{ex_mitternacht}
  We compute a \textsc{Thomas} decomposition of
  $\left\{\left(p:=a\underline{x}^2+bx+c\right)_=\right\} \subseteq \mathbb{Q}[a,b,c,x]$ with respect to $a < b < c < x$.
  We highlight the highest power of the leader by underlining it.
  
  First, we ensure that the initial $\ini(p)$ of $p$ is not zero.
  Therefore, we insert $\left(\ini(p)\right)_{\neq}=\left(\underline{a}\right)_{\neq}$ into the system.
  Since we restricted the solution set of this system, we also have to consider the system $\left\{p_=,\left(\underline{a}\right)_=\right\}$, which simplifies to $\left\{\left(b\underline{x}+c\right)_=,\left(\underline{a}\right)_=\right\}$.
  Similarly, we add $\left(\underline{b}\right)_{\neq}$ to ensure $\ini(b\underline{x}+c)\neq0$ and get the special case system $\left\{\left(\underline{c}\right)_=,\left(\underline{b}\right)_=,\left(\underline{a}\right)_=\right\}$.
  Up to this point, we have three systems, where the second and third one are easily checked to be simple:
  
  \begin{minipage}{.30\textwidth}
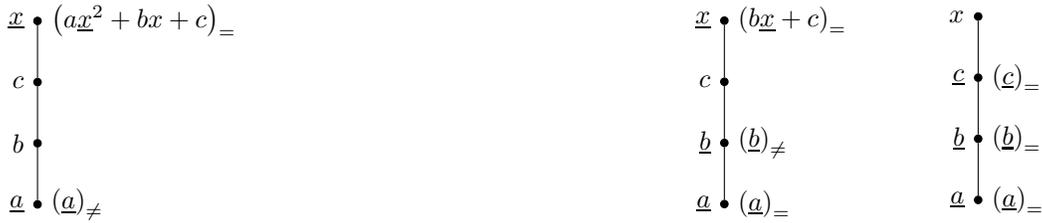
    \begin{tikzpicture}[dot/.style={circle,draw=black,fill=black,inner sep=1pt}]
    \node[dot,label=180:$\underline{x}$] (x) {};
    \node[dot,label=180:{$c$},below=7mm of x] (c) {};
    \node[dot,label=180:{$b$},below=7mm of c] (b) {};
    \node[dot,label=180:{$\underline{a}$},below=7mm of b] (a) {};
    \draw (x) -- (c);
    \draw (c) -- (b);
    \draw (b) -- (a);
    \node[right=0pt of x] (S2x) { $\left(a\underline{x}^2+bx+c\right)_=$ };
    \node[right=0pt of a] (S2a) { $\left(\underline{a}\right)_{\neq}$ };
    \end{tikzpicture}  \end{minipage}
    \begin{minipage}{.30\textwidth}\ \end{minipage}
    \begin{minipage}{.22\textwidth}
    \begin{tikzpicture}[dot/.style={circle,draw=black,fill=black,inner sep=1pt}]
    \node[dot,label=180:{$\underline{x}$}] (x) {};
    \node[dot,label=180:{$c$},below=7mm of x] (c) {};
    \node[dot,label=180:{$\underline{b}$},below=7mm of c] (b) {};
    \node[dot,label=180:{$\underline{a}$},below=7mm of b] (a) {};
    \draw (x) -- (c);
    \draw (c) -- (b);
    \draw (b) -- (a);
    \node[right=0pt of x] (S2x) { $\left(b\underline{x}+c\right)_=$ };
    \node[right=0pt of c] (S2c) { };
    \node[right=0pt of b] (S2b) { $\left(\underline{b}\right)_{\neq}$ };
    \node[right=0pt of a] (S2a) { $\left(\underline{a}\right)_=$ };
    \end{tikzpicture}  \end{minipage}
    \begin{minipage}{.18\textwidth}
    \begin{tikzpicture}[dot/.style={circle,draw=black,fill=black,inner sep=1pt}]
    \node[dot,label=180:$x$] (x) {};
    \node[dot,label=180:{$\underline{c}$},below=7mm of x] (c) {};
    \node[dot,label=180:{$\underline{b}$},below=7mm of c] (b) {};
    \node[dot,label=180:{$\underline{a}$},below=7mm of b] (a) {};
    \draw (x) -- (c);
    \draw (c) -- (b);
    \draw (b) -- (a);
    \node[right=0pt of c] (S2c) { $\left(\underline{c}\right)_=$ };
    \node[right=0pt of b] (S2b) { $\left(\underline{b}\right)_=$ };
    \node[right=0pt of a] (S2a) { $\left(\underline{a}\right)_=$ };
  \end{tikzpicture}  \end{minipage}
  
  Second, we ensure that $p$ is square-free by insertion of $\left(4a\underline{c}-b^2\right)_{\neq}$ into the first system.
  Again, we also need to consider the system $\{\left(p\right)_=,\left(4a\underline{c}-b^2\right)_=,\left(\underline{a}\right)_{\neq}\}$.
  As $p$ is a square in this system, we can replace it by its square-free part $2a\underline{x}-b$.
  Now, all systems are easily verified to be simple and we obtain the following \textsc{Thomas} decomposition:
  
  \begin{minipage}{.30\textwidth}
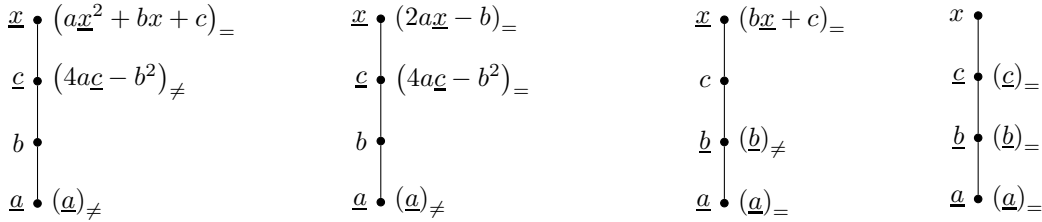
    \begin{tikzpicture}[dot/.style={circle,draw=black,fill=black,inner sep=1pt}]
    \node[dot,label=180:$\underline{x}$] (x) {};
    \node[dot,label=180:{$\underline{c}$},below=7mm of x] (c) {};
    \node[dot,label=180:{$b$},below=7mm of c] (b) {};
    \node[dot,label=180:{$\underline{a}$},below=7mm of b] (a) {};
    \draw (x) -- (c);
    \draw (c) -- (b);
    \draw (b) -- (a);
    \node[right=0pt of x] (S2x) { $\left(a\underline{x}^2+bx+c\right)_=$ };
    \node[right=0pt of c] (S2c) { $\left(4a\underline{c}-b^2\right)_{\neq}$};
    \node[right=0pt of a] (S2a) { $\left(\underline{a}\right)_{\neq}$ };
    \end{tikzpicture}  \end{minipage}
    \begin{minipage}{.30\textwidth}
    \begin{tikzpicture}[dot/.style={circle,draw=black,fill=black,inner sep=1pt}]
    \node[dot,label=180:$\underline{x}$] (x) {};
    \node[dot,label=180:{$\underline{c}$},below=7mm of x] (c) {};
    \node[dot,label=180:{$b$},below=7mm of c] (b) {};
    \node[dot,label=180:{$\underline{a}$},below=7mm of b] (a) {};
    \draw (x) -- (c);
    \draw (c) -- (b);
    \draw (b) -- (a);
    \node[right=0pt of x] (S2x) {$\left(2a\underline{x}-b\right)_=$};
    \node[right=0pt of c] (S2c) {$\left(4a\underline{c}-b^2\right)_=$ };
    \node[right=0pt of a] (S2a) {$\left(\underline{a}\right)_{\neq}$ };
    \end{tikzpicture}  \end{minipage}
    \begin{minipage}{.22\textwidth}
    \begin{tikzpicture}[dot/.style={circle,draw=black,fill=black,inner sep=1pt}]
    \node[dot,label=180:{$\underline{x}$}] (x) {};
    \node[dot,label=180:{$c$},below=7mm of x] (c) {};
    \node[dot,label=180:{$\underline{b}$},below=7mm of c] (b) {};
    \node[dot,label=180:{$\underline{a}$},below=7mm of b] (a) {};
    \draw (x) -- (c);
    \draw (c) -- (b);
    \draw (b) -- (a);
    \node[right=0pt of x] (S2x) { $\left(b\underline{x}+c\right)_=$ };
    \node[right=0pt of c] (S2c) { };
    \node[right=0pt of b] (S2b) { $\left(\underline{b}\right)_{\neq}$ };
    \node[right=0pt of a] (S2a) { $\left(\underline{a}\right)_=$ };
    \end{tikzpicture}  \end{minipage}
    \begin{minipage}{.18\textwidth}
    \begin{tikzpicture}[dot/.style={circle,draw=black,fill=black,inner sep=1pt}]
    \node[dot,label=180:$x$] (x) {};
    \node[dot,label=180:{$\underline{c}$},below=7mm of x] (c) {};
    \node[dot,label=180:{$\underline{b}$},below=7mm of c] (b) {};
    \node[dot,label=180:{$\underline{a}$},below=7mm of b] (a) {};
    \draw (x) -- (c);
    \draw (c) -- (b);
    \draw (b) -- (a);
    \node[right=0pt of c] (S2c) { $\left(\underline{c}\right)_=$ };
    \node[right=0pt of b] (S2b) { $\left(\underline{b}\right)_=$ };
    \node[right=0pt of a] (S2a) { $\left(\underline{a}\right)_=$ };
    \end{tikzpicture}  \end{minipage}
\end{myexample}

\subsection{Algebraic Thomas Decomposition}\label{algebraic_algorithms}

This section presents our main algorithm for algebraic systems and its subalgorithms.
The algorithm represents each system as a pair consisting of a candidate simple system and a queue of unprocessed equations and inequations.\footnote{This approach has been adapted from \cite{GB2}, where $T$ was an intermediate Janet basis and $Q$ a queue of new prolongations to be checked. A similar approach was later used for triangular decompositions in \cite{MazaMEGA2000}.}
During each step, the algorithm chooses a suitable polynomial from the queue, pseudo-reduces it and afterwards combines it with the polynomial from the candidate simple system having the same leader.
In this process, the algorithm may split the system, i.e., add a new polynomial into the queue as an inequation and at the same time create a new subsystem with the same polynomial added to the queue as an equation.
This way, we ensure that no solutions are lost and the solution sets are disjoint.
The algorithm considers a system inconsistent and discards it when an equation of the form $c_=$ with $c \in F\setminus\{0\}$ or the inequation $0_{\neq}$ is produced.

We consider a system $S$ as a pair of sets $(S_T, S_Q)$, where $S_T$ represents the candidate simple system and $S_Q$ is the queue.
We require $S_T$ to be triangular and thus $(S_T)_x$ denotes the unique equation or inequation of leader $x$ in $S_T$, if any.
Moreover, $S_T$ must fulfill a weaker form of the other two simplicity conditions, in particular, in conditions (\ref{simple_system})(\ref{simple_nonzeroinitials}) and (\ref{simple_squarefree}), the tuple $\ma$ can be a solution of $(S_T)_{<x} \cup (S_Q)_{<x}$ instead of just $(S_T)_{<x}$.
Obviously, $S_Q=\emptyset$ implies simplicity of $S$.

From now on, let $\prem$ be a \textbf{pseudo remainder algorithm}\footnote{In our context $\prem$ does not necessarily have to be the classical pseudo remainder, but any sparse pseudo remainder with property (\ref{prempquo}) will suffice.} in $R$ and $\pquo$ the corresponding \textbf{pseudo quotient algorithm}.
To be precise, if $p, q\in R$ with $\ld(p)=\ld(q)=x$, then \begin{equation}\label{prempquo}m\cdot p = \pquo(p,q,x) \cdot q+\prem(p,q,x)\end{equation} holds, where $\deg_x(q)>\deg_x(\prem(p,q,x))$, $\ld(m)<x$ and $m \mid \ini(q)^k$ for some $k \in \mathbb{Z}_{\geq0}$.
Note that $\pma(\ini(p))\not=0$ and $\pma(\ini(q))\not=0$ imply $\pma(\pquo(p,q,x))\not=0$ and $\pma(m)\not=0$.

The following algorithm employs \textsf{prem} to reduce a polynomial modulo $S_T$:

\begin{myalgorithm}[{\sf Reduce}]\label{algo_reduce}\ \\
 \textit{Input:} A system $S$, a polynomial $p \in R$ \\
 \textit{Output:} A polynomial $q$ with $\pma(p)=0$ if and only if $\pma(q)=0$ for each $\ma\in\sol(S)$.\\
 \textit{Algorithm:} \begin{algorithmic}[1]
                      \STATE $x \gets \ld(p)$; $q\gets p$ 
                      \WHILE{$x>1$ and $(S_T)_x$ is an equation and $\mdeg(q) \geq \mdeg((S_T)_x)$}
                       \STATE $q \gets \prem(q, (S_T)_x, x)$
                       \STATE $x \gets \ld(q)$
                      \ENDWHILE
                      \IF{$x>1$ and $\textrm{\sf Reduce}(S, \ini(q)) = 0$}
                       \RETURN $\textrm{\sf Reduce}(S, q-\ini(q)x^{\mdeg(q)})$
                      \ELSE
                       \RETURN $q$
                      \ENDIF
                     \end{algorithmic}
\end{myalgorithm}
 \begin{proof}[Correctness]
  There exist $m \in R\setminus\{0\}$ with $\ld(m)<\ld(p)$ and $\pma(m) \neq 0$ for all $\ma \in \sol(S_{\le\ld(p)})$ such that
  \[
    \mathsf{Reduce}(S, p)=mp-\sum_{y \le \ld(p)} c_y\cdot (S_T)_y
  \]
  with $c_y \in R$ and $\ld(c_y)\le \ld(p)$ if $(S_T)_y$ is an equation and $c_y=0$ otherwise. This implies
  \[
    \pma(\mathsf{Reduce}(S, p)) = \underbrace{\pma(m)}_{\neq0}\pma(p)-\sum_{y \le x} \pma(c_y)\underbrace{\pma((S_T)_y)}_{=0}
  \]
  and therefore $\pma(p) = 0$ if and only if $\pma(\mathsf{Reduce}(S, p)) = 0$.\qed
 \end{proof}

Note that this algorithm only uses the equation part of the triangular system in $S$, i.e. $S_T^=$.
For ease of notation in the following algorithms, we write $\mathsf{Reduce}(S,p)$ instead of $\mathsf{Reduce}(S_T^=,p)$.

A polynomial $p$ \textbf{reduces to $q$ modulo $S_T$} if $\mathsf{Reduce}(S, p)=q$.
A polynomial is \textbf{reduced modulo $S_T$} if it reduces to itself.

The \textsf{Reduce} algorithm differs slightly from the classical $\prem(p, S_T^=)$ as defined in \cite{AubryTriTheo}.
While $\prem(p, S_T^=)$ fully reduces $p$ modulo all variables, $\mathsf{Reduce}(S, p)$ only reduces modulo the leader and ensures that the initial of the reduced form doesn't vanish.
Performing $\mathsf{Reduce}(S, p)$ in combination with a full coefficient reduction (see also \textsection\ref{section_algorithmic_opt}) is the same as computing $\prem(p, S_T^=)$.
It is therefore possible to replace $\mathsf{Reduce}(S, p)$ with $\prem(p, S_T^=)$ in the following algorithms.
Our approach adds some flexibility, as we can choose to omit a full reduction in an implementation.
In particular, if a polynomial does not reduce to zero, we can determine that without performing a full $\prem$ reduction.
We apply this multiple times in our implementation, most prominently in Algorithm (\ref{algo_ressplit}).
However, if a polynomial reduces to zero, $\mathsf{Reduce}$ has no advantage over $\prem$.

Later, we will use the following facts about the \textsf{Reduce} algorithm.

\begin{myremark}\label{rem_reduceprop} Let $q=\mathsf{Reduce}(S,p)\neq0$.
 \begin{enumerate}
  \item If $S_{\ld(q)}$ is an equation, then $\mdeg(q)<\mdeg(S_{\ld(q)})$.\label{rem_reduceprop1}
  \item $\mathsf{Reduce}(S, \ini(\mathsf{Reduce}(S,p)))\neq0$.\label{rem_reduceprop2}
  \item $\ld(q)\leq\ld(p)$ and if $\ld(q)=\ld(p)$, then $\mdeg(q)\leq\mdeg(p)$.\label{rem_reduceprop3}
 \end{enumerate}
\end{myremark}

The result of the {\sf Reduce} algorithm does not need to be a canonical normal form, however, the algorithm recognizes polynomials that vanish on all solutions:

\begin{mycorollary}\label{reduce0}
 Let $p \in R$ with $\ld(p)=x$. $\mathsf{Reduce}(S, p)=0$ implies $\pma(p)=0\ \forall\ \ma\in\sol(S_{\leq x})$.
\end{mycorollary}
 \begin{proof}
  For all $\ma \in \sol(S_{\le x})$, it holds that $\pma(p) = 0$ if and only if $\pma(\mathsf{Reduce}(S, p)) = 0$. The statement follows from $\pma(\mathsf{Reduce}(S, p)) = \pma(0) = 0$.\qed
 \end{proof}

The converse of this corollary doesn't hold in general.
Thus, we provide two weaker statements in the following remark.

\begin{myremark}
 Let $p$ and $x$ as in Corollary (\ref{reduce0}).
 \begin{enumerate}
  \item If $(S_Q)_{\le x}=\emptyset$, i.e., $S_{\le x} = (S_T)_{\le x}$ is simple, then $\mathsf{Reduce}(S, p) \neq 0$ implies $\exists\ \ma \in \sol(S_{\le x})$ such that $\pma(p)\neq0$.
  \item If $(S_Q)^=_{<x}=\emptyset$ and $\textsf{Reduce}(S,p)\neq0$ hold, then either $\sol(S_{<x})=\emptyset$ or $\exists\ \ma\in\sol(S_{<x} \cup \{ (S_T)_x \})$ such that $\pma(p)\neq0$.
 \end{enumerate}
\end{myremark}

\proof We only prove the second part, as the first part easily follows.

Let $(S_Q)^=_{<x}=\emptyset$, $\textsf{Reduce}(S,p)\neq0$ and $|\sol(S_{<x})|>0$.
First, as $\ld((S_T)_x)=x$ and $\mdeg((S_T)_x)>0$, for each $\ma \in \sol(S_{<x})$, the univariate polynomial $\pma[<x]((S_T)_x) \in \overline{F}[x]$ has positive degree.
Thus $|\sol(S_{<x} \cup \{ (S_T)_x \})|>0$.

Let $\pma(p)=0\ \forall\ \ma \in \sol(S_{<x} \cup \{ (S_T)_x \})$ (*).
Then $(S_T)_x$ is an equation and $\deg_x(p)\ge\deg_x((S_T)_x)$ and therefore $p \neq \mathsf{Reduce}(S, p)$.
In fact, (*) further implies $\ld(\mathsf{Reduce}(S,p)) < x$, as otherwise $\deg_x(\mathsf{Reduce}(S,p))\ge\deg_x((S_T)_x)$ would hold.
By repeating the previous arguments, we can inductively conclude $\ld(\mathsf{Reduce}(S,p))=1$.
As $\pma(p)=0$, we conclude $\mathsf{Reduce}(S,p)=0$, a contradiction. \qed

The first part of this remark in conjunction with Corollary (\ref{reduce0}) implies \cite[Thm.~4]{wang_simple}.

\begin{myexample}
  Reduce $q_1:=x^2+y^2x+x+y$ modulo the simple system on the left.

  
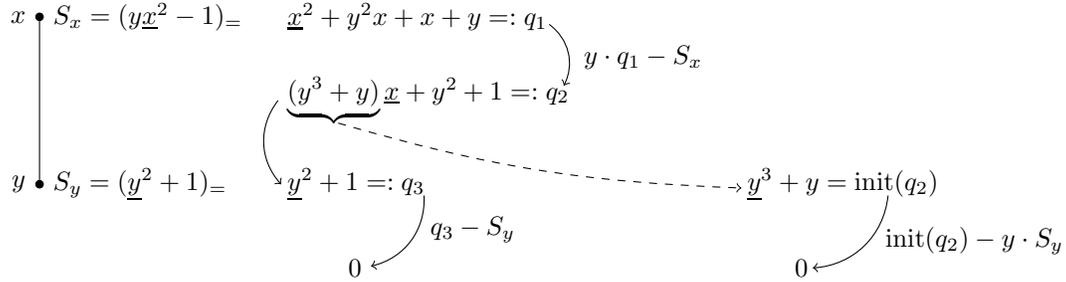
\begin{tikzpicture}[dot/.style={circle,draw=black,fill=black,inner sep=1pt}]
    \node[dot,label=180:$x$] (x) at (0,0) {};
    \phantom{\node[dot,below=of x] (no) {};}
    \node[dot,label=180:$y$,below=of no] (y) {};
    \phantom{\node[dot,right=3cm of x] (x2) {};
    \node[dot,below=of x2] (no2) {};
    \node[dot,below=of no2] (y2) {};
    \node[dot,below=of y2] (O2) {};
    \node[dot,right=6cm of y2] (y3) {};
    \node[dot,right=6cm of O2] (O3) {};}
    \draw (x) -- (y);
    \node[right=0pt of x] (S1x) { $S_x = (y\underline{x}^2-1)_=$ };
    \node[right=0pt of y] (S1y) { $S_y = (\underline{y}^2+1)_=$ };
    \node[right=0pt of x2] (r3) { $\underline{x}^2+y^2x+x+y =: q_1$ };
    \node[right=0pt of no2] (r4) { $\underbrace{(y^3+y)}\underline{x} + y^2 + 1 =: q_2$ };
    \draw[->] ($(r3)+(1.75cm,-0.1cm)$) to [bend left=40] ($(r4)+(1.8cm,0.2cm)$) {}; \node[] at ($(r3)+(3cm,-.54cm)$) { $y\cdot q_1-S_x$ };
    \node[right=0pt of y2] (r5) { $\underline{y}^2+1 =: q_3$ };
    \draw[->] ($(r4)+(-2cm,0)$) to [bend right=40] ($(r5)+(-1cm,0)$) {};
    \node[right=23pt of O2] (r6) { $0$ };
    \draw[->] ($(r5)+(0.9cm,-0.15cm)$) to [bend left=40] (r6) {}; \node[] at ($(r5)+(1.55cm,-0.6cm)$) { $q_3-S_y$ };
    \node[right=0pt of y3] (r7) { $\underline{y}^3+y=\ini(q_2)$ };
    \draw[->] ($(r4)+(-1.2cm,-.3cm)$) to [bend right=7,dashed] ($(r7)+(-1.35cm,-.05cm)$) {};
    \node[right=18pt of O3] (r8) { $0$ };
    \draw[->] ($(r7)+(0.6cm,-0.15cm)$) to [bend left=40] ($(r8)+(.15cm,0)$) {}; \node[] at ($(r7)+(1.75cm,-0.75cm)$) { $\ini(q_2)-y\cdot S_y$ };
  \end{tikzpicture}

  In the first reduction step, $q_1$ is pseudo-reduced modulo $S_x$.
  The result $q_2$ still has leader $x$, but a main degree smaller than $S_x$.
  We determine that the initial of $q_2$ reduces to $0$ and remove the highest power of $x$ from $q_2$.
  The resulting polynomial $q_3$ now pseudo-reduces to $0$ modulo $S_y$, i.e. $\mathsf{Reduce}(\{S_x,S_y\}, q_1)=0$.
\end{myexample}

Now, we examine all splitting methods needed during the algorithm.
We will use the following one-liner as subalgorithm for the splitting subalgorithms.

\begin{myalgorithm}[{\sf Split}]
 \textit{Input:} A system $S$, a polynomial $p \in R$ \\
 \textit{Output:} The disjoint decomposition $\left( S \cup \left\{ p_{\neq} \right\}, S \cup \left\{ p_{=} \right\} \right)$ of $S$.\\
 \textit{Algorithm:} \begin{algorithmic}[1]
             \RETURN $\left( \left(S_T, S_Q \cup \{ p_{\neq} \}\right), \left(S_T, S_Q \cup \{ p_{=} \}\right) \right)$
            \end{algorithmic}
\end{myalgorithm}

For a better understanding of the following splitting subalgorithms we first need to explain how they are applied in the main algorithm.
Each step of the algorithm treats a system $S$ as follows.
An equation or inequation $q$ is chosen and removed from the queue $S_Q$.
Then we reduce $q$ modulo $S_T$.
For the simplicity properties to hold w.r.t. $q$ it is necessary to add inequations to $S$.
To accomplish this, we pass $S$ together with $q$ to the splitting subalgorithms.
Each such subalgorithm returns two systems.
The first system $S_1$ contains an additional inequation.
The second system $S_2$ contains a complementary equation, $q$ is added back into the queue of $S_2$, and $S_2$ is put aside for later treatment.
In each case $(S_1 \cup \{q\}, S_2)$ is a disjoint decomposition of the original system $S \cup \{q\}$.
Then $S_1$ and $q$ may be subjected to further splitting algorithms and eventually $q$ is added into the candidate simple system.

The first splitting algorithm we consider is \textsf{InitSplit}, which is concerned with property (\ref{simple_system})(\ref{simple_nonzeroinitials}).

\begin{myalgorithm}[\sf InitSplit]\label{algo_initsplit}
 \textit{Input:} A system $S$, an equation or inequation $q$ with $\ld(q)=x$. \\
 \textit{Output:} Two systems $S_1$ and $S_2$, where $\left(S_1 \cup \{q\}, S_2\right)$ is a disjoint decomposition of $S \cup \{q\}$. Moreover, $\pma(\ini(q))\neq0$ holds for all $\ma \in \sol(S_1)$ and $\pma(\ini(q))=0$ for all $\ma \in \sol(S_2)$.
 \textit{Algorithm:} \begin{algorithmic}[1]
    \STATE $(S_1, S_2) \gets \textrm{\sf Split}(S, \ini(q))$
    \STATE $(S_2)_Q \gets (S_2)_Q \cup \left\{ q \right\}$
    \RETURN $(S_1, S_2)$
   \end{algorithmic}
\end{myalgorithm}

For the further splitting algorithms, we need some preparation.
In Definition (\ref{simple_system}) we consider a multivariate polynomial $p$ as the family of univariate polynomials $\pma[<\ld(p)](p)$.
For ensuring triangularity and square-freeness, we have to compute the gcd (greatest common divisor) of two polynomials, which in general depends on $\ma$.
Subresultants provide a generalization of the \textsc{Euclid}ean algorithm and enable us to take the tuple $\ma$ into account.

\begin{mydefinition}\label{not_prs}
 Let $p, q \in R$ with $\ld(p)=\ld(q)=x$, $\deg_x(p)=d_p > \deg_x(q)=d_q$.
 We denote by $\prs(p,q,x)$ the \textbf{subresultant polynomial remainder sequence} (see \cite{habicht}, \cite[Chap.~7]{mishra}, \cite[Chap.~3]{Yap}) of $p$ and $q$ w.r.t. $x$, and  by $\prs_i(p,q,x)$, $i<d_q$ the \emph{regular} polynomial of degree $i$ in $\prs(p,q,x)$ if it exists, or $0$ otherwise.
 Furthermore, $\prs_{d_p}(p,q,x):=p$, $\prs_{d_q}(p,q,x):=q$ and $\prs_i(p,q,x):=0$, $d_q<i<d_p$.

 Define $\res_i(p,q,x) := \ini\left(\prs_i\left(p,q,x\right)\right)$ for $0<i<d_p$, $\res_{d_p}(p,q,x):=1$ and $\res_0(p,q,x):=\prs_0\left(p,q,x\right)$.
 Note that $\res_0(p,q,x)$ is the usual resultant.
 \footnote{These definitions are slightly different from the ones cited in the literature (\cite[Chap.~7]{mishra}, \cite[Chap.~3]{Yap}), since we only use the regular subresultants.
 However, it is easy to see that all theorems from \cite[Chap.~7]{mishra} we refer to still hold for $i<d_q$.}
\end{mydefinition}

The initials of the subresultants provide conditions to determine the degrees of all possible gcds.
Using these conditions, we describe the splittings necessary to determine degrees of polynomials within one system.

\begin{mydefinition}\label{def_fibrcard}
 Let $S$ be a system and $p_1, p_2 \in R$ with $\ld(p_1)=\ld(p_2)=x$. If $|\sol(S_{<x})|>0$, we call
 \[
   i := \min\left\{ i \in \mathbb{Z}_{\ge0} \mid \exists\ \ma \in \sol(S_{<x}) \mbox{ such that } \deg_x(\gcd(\pma[<x](p_1),\pma[<x](p_2))) = i \right\}
 \]
 the \textbf{fiber cardinality} of $p_1$ and $p_2$ w.r.t. $S$. Moreover, if $(S_Q)_{<x}^==\emptyset$, then
 \[
   i^\prime := \min \{ i \in \mathbb{Z}_{\ge0} \mid \textrm{\sf Reduce}(\res_j(p_1,p_2,x), S_T)=0\ \forall\ j<i \textrm{ and } \textrm{\sf Reduce}(\res_i(p_1,p_2,x), S_T)\neq0\}
 \]
 is the \textbf{quasi fiber cardinality} of $p_1$ and $p_2$ w.r.t. $S$. A disjoint decomposition $(S_1, S_2)$ of $S$ such that \begin{enumerate}
   \item $\deg_x(\gcd(\pma[<x](p_1),\pma[<x](p_2))) = i\ \forall\ \ma \in \sol\left((S_1)_{<x}\right)$
   \item $\deg_x(\gcd(\pma[<x](p_1),\pma[<x](p_2))) > i\ \forall\ \ma \in \sol\left((S_2)_{<x}\right)$
  \end{enumerate}
 is called $i$-th \textbf{fibration split} of $p_1$ and $p_2$ w.r.t. $S$. A polynomial $r \in R$ with $\ld(r)=x$ such that $\deg_x(r)=i$ and
 \[\pma[<x](r) \sim \gcd(\pma[<x](p_1),\pma[<x](p_2))\ \forall\ \ma \in \sol\left((S_1)_{<x}\right)\] is called $i$-th \textbf{conditional greatest common divisor} of $p_1$ and $p_2$ w.r.t. $S$, where $p \sim q$ if and only if $p \in (\overline{F} \setminus \{0\})q$.
 Furthermore, $q \in R$ with $\ld(q)=x$ and $\deg_x(q)=\deg_x(p_1)-i$ such that \[\pma[<x](q) \sim \frac{\pma[<x](p_1)}{\gcd(\pma[<x](p_1),\pma[<x](p_2))}\ \forall\ \ma \in \sol\left((S_1)_{<x}\right)\] is called the $i$-th \textbf{conditional quotient} of $p_1$ by $p_2$ w.r.t. $S$.
 By replacing $\pma[<x](p_2)$ in the above definition with $\frac{\partial}{\partial x}(\pma[<x](p_1))$, we get an $i$-th \textbf{square-free split} and $i$-th \textbf{conditional square-free part} of $p_1$ w.r.t. $S$.
\end{mydefinition}

\begin{myexample}
 Consider the system $S := \{ (\underline{x}^3+y)_= \}$ and the polynomial $q := \underline{x}^2+x+y+1$ with $y<x$.
 Compute $res_0(S_x, q)=\underline{y}^3+7y^2+5y+1$, $res_1(S_x, q)=-\underline{y}$ and $res_2(S_x,q)=1$.
 The fiber cardinality is of $S_x$ and $q$ w.r.t. $S$ is $0$.
 A zeroth fibration split is given by $S_1 := S \cup \{ (res_0(S_x, q))_{\neq} \}$ and $S_2 := S \cup \{ (res_0(S_x, q))_= \}$.
 The fiber cardinality w.r.t. $S_2$ is $1$.
 A first fibration split is given by $S_{2,1} := S_2 \cup \{ (-\underline{y})_{\neq} \}$ and $S_{2,2} := S \cup \{ (-\underline{y})_= \}$.
 Note in this case that $\sol(S_{2,1})=\sol(S_2)$ and $\sol(S_{2,2})=\emptyset$.
 A zeroth conditional quotient of $S_x$ and $q$ is $S_x$.
 A first conditional gcd and first conditional quotient are $-y\underline{x} + 2y + 1$ and $y^2\underline{x}^2+(2y^2+y)x+4y^2+4y+1$, respectively.
\end{myexample}

It is in general hardly possible to compute the fiber cardinality directly.
However, in the case where the quasi fiber cardinality is strictly smaller than the fiber cardinality, the corresponding fibration split will lead to one inconsistent system, and one where the quasi fiber cardinality is increased.

\begin{mylemma}\label{lm_ressplit}
 Let $|\sol(S_{<x})|>0$ and $(S_Q)_{<x}^==\emptyset$. For $p_1$, $p_2$ as in Definition (\ref{def_fibrcard}) with $\pma(\ini(p_1))\neq0\ \forall\ \ma\in \sol(S_{<x})$ and $\mdeg(p_1)>\mdeg(p_2)$, let $i$ be the fiber cardinality of $p_1$ and $p_2$ w.r.t. $S$ and $i^\prime$ the corresponding quasi fiber cardinality.
 Then \[i^\prime \leq i\] where the equality holds if and only if $\left|\sol\left(S_{<x} \cup \{\res_{i^\prime}(p_1,p_2,x)_{\neq}\}\right)\right|>0$.
\end{mylemma}
 \begin{proof}
 Let $\ma \in \sol(S_{<x})$, $\mdeg(p_1)>\mdeg(p_2)$, $d_{p_1} := \deg_x(p_1) = \deg_x(\pma[<x](p_1))$, $d_{p_2} := \deg_x(p_2)$ and $d_{p_2,\ma} := \deg_x(\pma[<x](p_2))$. If $i < \max(d_{p_1},d_{p_2,\ma})-1 = d_{p_1}-1$, then \cite[Thm.~7.8.1]{mishra} implies
 \begin{equation}\label{lm_ressplit_proof1} \pma[<x](\prs_i(p_1,p_2,x))\sim\prs_i(\pma[<x](p_1),\pma[<x](p_2),x)\end{equation}
 and
 \begin{equation}\label{lm_ressplit_proof2}\pma(\res_i(p_1,p_2,x))=0 \Longleftrightarrow \res_i(\pma[<x](p_1),\pma[<x](p_2),x)=0\enspace.\end{equation}
 Conditions (\ref{lm_ressplit_proof1}) and (\ref{lm_ressplit_proof2}) by definition also hold for the trivial cases $d_{p_2} \le i \le d_{p_1}$.

 For all indices $j<i^\prime$, Corollary (\ref{reduce0}) and the fact $\mathsf{Reduce}(\res_j(p_1,p_2,x), S_T)=0$ imply $\pma(\res_j(p_1,p_2,x))=0$.
 By (\ref{lm_ressplit_proof1}) and (\ref{lm_ressplit_proof2}), $\res_j(\pma[<x](p_1),\pma[<x](p_2),x)=0$ follows.
 We apply \cite[Thm.~7.10.5]{mishra} successively and get $\prs_j(\pma[<x](p_1),\pma[<x](p_2),x)=0$.
 Thus, \begin{equation} \label{lm_ressplit_proof3} \deg_x(\gcd(\pma[<x](p_1),\pma[<x](p_2)))\ge i^\prime\end{equation} holds.
 This implies $i^\prime \le i$.

 Equality in (\ref{lm_ressplit_proof3}) holds if and only if there exists $\ma \in \sol(S_{<x})$ such that $\phi_\ma(\res_{i^\prime}(p_1,p_2,x))\neq0$.
 Therefore, $i=i^\prime$ if and only if $\left|\sol\left(S_{<x} \cup \{\res_{i^\prime}(p_1,p_2,x)_{\neq}\}\right)\right|>0$.\qed
 \end{proof}

The above lemma doesn't apply if both polynomials have the same degree.
In this case, both polynomials must have non-vanishing initials, as shown in the following corollary.

\begin{mycorollary}\label{cor_ressplit}
 Let $|\sol(S_{<x})|>0$ and $(S_Q)_{<x}^==\emptyset$. For polynomials $p_1$, $p_2$ as in Definition (\ref{def_fibrcard}) with $\pma(\ini(p_1))\neq0$ and $\pma(\ini(p_2))\neq0\ \forall\ \ma\in \sol(S_{<x})$, let $i$ be the fiber cardinality of $p_1$ and $p_2$ w.r.t. $S$ and $i^\prime$ the quasi fiber cardinality of $p_1$ and $\prem(p_2,p_1,x)$ w.r.t. $S$. Then \[i^\prime \leq i\] with equality if and only if $\left|\sol\left(S_{<x} \cup \{\res_{i^\prime}(p_1,\prem(p_2,p_1,x),x)_{\neq}\}\right)\right|>0$.
\end{mycorollary}
 \begin{proof}
  Let $\ma\in \sol(S_{<x})$.
  By the assumption on the initials, \cite[Corr.~7.5.6]{mishra} implies $\pma[<x](\prem(p_2,p_1,x))=\prem(\pma[<x](p_2),\pma[<x](p_1),x)$.
  The univariate polynomials $\pma[<x](p_1)$ and $\pma[<x](p_2)$ have the same gcd as $\pma[<x](p_1)$ and $\prem(\pma[<x](p_2),\pma[<x](p_1),x)$.
  We can therefore replace $p_2$ with $\prem(p_2,p_1,x)$ in Lemma (\ref{lm_ressplit}).\qed
 \end{proof}

The following algorithm computes the quasi fiber cardinality of two polynomials.

\begin{myalgorithm}[\sf ResSplit]\label{algo_ressplit}
 \textit{Input:} A system $S$ with $(S_Q)_{<x}^==\emptyset$, two polynomials $p, q \in R$ with $\ld(p)=\ld(q)=x$, $\mdeg(p)>\mdeg(q)$ and $\pma(\ini(p))\neq0$ for all $\ma \in \sol(S_{<x})$.\\
 \textit{Output:} The quasi fiber cardinality $i$ of $p$ and $q$ w.r.t. $S$ and an $i$-th fibration split $(S_1,S_2)$ of $p$ and $q$ w.r.t. $S$.\\
 \textit{Algorithm:} \begin{algorithmic}[1]
    \STATE \label{ressplit_cond_nontrivial_split} $i \gets \min \{ i \in \mathbb{Z}_{\ge0} \mid \textrm{\sf Reduce}(\res_j(p,q,x), S_T)=0\ \forall\ j<i \textrm{ and } \textrm{\sf Reduce}(\res_i(p,q,x), S_T)\neq0$\}
    \RETURN $(i, S_1, S_2) := \left(i, \textrm{\sf Split}(S, \res_i(p,q,x)) \right)$
   \end{algorithmic}
\end{myalgorithm}
\begin{proof}[Correctness]
 Assume $|\sol((S_l)_{<x})|>0$, $l=1,2$, as the statement is trivial otherwise.

 Let $\ma \in \sol((S_1)_{<x})$.
 The polynomial $g:=\prs_i(\pma[<x](p),\pma[<x](q),x)$ is not identically zero, due to $(\ini(g))_{\not=}=(\res_i(p,q,x))_{\neq} \in (S_1)_Q$.
 The degree of $g$ is $i$ and $g\sim\gcd(\pma[<x](p),\pma[<x](q))$, as discussed in the proof of Lemma (\ref{lm_ressplit}).

 Let $\ma \in\sol((S_2)_{<x})$.
 \cite[Thm.~7.10.5]{mishra} and $(\ini(g))_==(\res_i(p,q,x))_= \in (S_2)_Q$ imply $g\equiv0$.
 Therefore, $\deg_x(\gcd(\pma[<x](p),\pma[<x](q))) > i$.\qed
\end{proof}

We apply the fiber cardinality and fibration split to compute a greatest common divisor of an existing polynomial in $S_T$ and another polynomial.

\begin{myalgorithm}[\sf ResSplitGCD]\label{algo_ressplitgcd}
 \textit{Input:} A system $S$ with $(S_Q)_{<x}^==\emptyset$, where $(S_T)_x$ is an equation, and an equation $q_=$ with $\ld(q)=x$. Furthermore, $\mdeg(q)<\mdeg((S_T)_x)$.\\
 \textit{Output:} Two systems $S_1$ and $S_2$ and an equation $\widetilde{q}_=$ such that:
   \renewcommand{\theenumi}{\alph{enumi}}
   \renewcommand{\labelenumi}{\theenumi)}
   \begin{enumerate}
    \item\label{algo_ressplitgcd_cond1} $S_2 = \widetilde{S_2} \cup \{ q \}$ where $\left(S_1, \widetilde{S_2}\right)$ is an $i$-th fibration split of $(S_T)_x$ and $q$ w.r.t. $S$,
    \item\label{algo_ressplitgcd_cond2} $\widetilde{q}$ is an $i$-th conditional gcd of $(S_T)_x$ and $q$ w.r.t. $S$,
   \end{enumerate}
   \renewcommand{\theenumi}{\arabic{enumi}}
   \renewcommand{\labelenumi}{\theenumi)}
   where $i$ is the quasi fiber cardinality of $p$ and $q$ w.r.t. $S$.\\
 \textit{Algorithm:} \begin{algorithmic}[1]
      \STATE $(i, S_1, S_2) \gets \textrm{\sf ResSplit}\left(S, (S_T)_x, q\right)$
      \STATE \label{algo_ressplitgcd_line2}$(S_2)_Q \gets (S_2)_Q \cup \{ q \}$
      \RETURN $S_1, S_2, \prs_i((S_T)_x, q, x)_{=}$
   \end{algorithmic}
\end{myalgorithm}
\begin{proof}[Correctness]
 Property \ref{algo_ressplitgcd_cond1}) follows from Algorithm (\ref{algo_ressplit}) and line \ref{algo_ressplitgcd_line2}.
 Property \ref{algo_ressplitgcd_cond2}) was already shown in the correctness proof of Algorithm (\ref{algo_ressplit}).\qed
\end{proof}

Note that $i>0$ is required in this case, as $i=0$ would yield an inconsistency.
Therefore, before calling \textsf{ResSplitGCD}, we will always ensure this condition in the main algorithm by incorporating the resultant of two equations into the system.

The following algorithm is similar.
But instead of the gcd, it returns the first input polynomial divided by the gcd.

\begin{myalgorithm}[\sf ResSplitDivide]\label{algo_ressplitdivide}
 \textit{Input:} A system $S$ with $(S_Q)_{<x}^==\emptyset$ and two polynomials $p$, $q$ with $\ld(p)=\ld(q)=x$ and $\pma(\ini(p))\neq0$ for all $\ma \in \sol(S_{<x})$.
                 Furthermore, if $\mdeg(p)\le\mdeg(q)$ then $\pma(\ini(q))\neq0$.\\ 
 \textit{Output:} Two systems $S_1$ and $S_2$ and a polynomial $\widetilde{p}$ such that:
   \renewcommand{\theenumi}{\alph{enumi}}
   \renewcommand{\labelenumi}{\theenumi)}
  \begin{enumerate}
   \item\label{algo_ressplitdivide_cond1} $S_2 = \widetilde{S_2} \cup \{ q \}$ where $\left(S_1,\widetilde{S_2}\right)$ is an $i$-th fibration split of $p$ and $q^\prime$ w.r.t. $S$,
   \item\label{algo_ressplitdivide_cond2} $\widetilde{p}$ is an $i$-th conditional quotient of $p$ by $q^\prime$ w.r.t. $S$,
  \end{enumerate} where $i$ is the quasi fiber cardinality of $p$ and $q^\prime$ w.r.t. $S$, with $q^\prime=q$ for $\mdeg(p)>\mdeg(q)$ and $q^\prime=\prem(q,p,x)$ otherwise.\\
   \renewcommand{\theenumi}{\arabic{enumi}}
   \renewcommand{\labelenumi}{\theenumi)}
 \textit{Algorithm:} \begin{algorithmic}[1]
  \IF{$\mdeg(p)\le\mdeg(q)$}
   \RETURN $\mathsf{ResSplitDivide}(S, p, \prem(q,p,x))$
  \ELSE
   \STATE $(i, S_1, S_2) \gets \textrm{\sf ResSplit}\left(S, p, q\right)$
   \IF{$i>0$}
     \STATE $\widetilde{p} \gets \pquo(p, \prs_i(p, \prem(q,p,x), x),x)$
   \ELSE
     \STATE $\widetilde{p} \gets p$
   \ENDIF
   \STATE $(S_2)_Q \gets (S_2)_Q \cup \{ q \}$ \label{algo_ressplitdivide_line7}
   \RETURN $S_1, S_2, \widetilde{p}$
  \ENDIF
 \end{algorithmic}
\end{myalgorithm}
 \begin{proof}[Correctness]
  According to Corollary (\ref{cor_ressplit}), we can without loss of generality assume $\mdeg(p)>\mdeg(q)$.

  Property \ref{algo_ressplitgcd_cond1}) follows from Algorithm (\ref{algo_ressplit}) and line \ref{algo_ressplitdivide_line7}.
  For all $\ma \in \sol(S_1)$, the following holds:
  If $i=0$, then $\deg_x(\gcd(\pma[<x](p), \pma[<x](q^\prime)))=0$ and thus $\pma[<x](p)$ shares no roots with $\pma[<x](q^\prime)$.
  Now let $i>0$.
  Formula (\ref{prempquo}) implies
    \[m\cdot p = \widetilde{p} \cdot \prs_i\left(p, q^\prime, x\right)+\prem\left(p,\prs_i\left(p, q^\prime, x\right),x\right)\mbox{ .}\]
  Due to \cite[Cor.~7.5.6]{mishra} and (\ref{lm_ressplit_proof1}), (\ref{lm_ressplit_proof2}) there exist $k_1, k_2 \in F \setminus \{0\}$ such that
   \[\begin{array}{rl} & \underbrace{\pma(m)}_{\neq 0}\cdot \pma[<x](p) \\
   = & \pma[<x](\widetilde{p}) \cdot \pma[<x]\left(\prs_i\left(p, q, x\right)\right)+\pma[<x]\left(\prem(p,\prs_i\left(p, q, x\right),x)\right)\\
   = & \pma[<x](\widetilde{p}) \cdot k_1 \prs_i\left(\pma[<x](p), \pma[<x](q), x\right)+k_2\prem(\pma[<x](p),\underbrace{\prs_i\left(\pma[<x](p), \pma[<x](q), x\right)}_{\mbox{\scriptsize divides } \pma[<x](p)},x)\\
   =& \pma[<x](\widetilde{p})\cdot k_1\gcd(\pma[<x](p), \pma[<x](q)) + 0 \mbox{ .}\end{array}\]
  Thus, we obtain property \ref{algo_ressplitdivide_cond2}) from
   \[\pma[<x](\widetilde{p}) \sim \frac{\pma[<x](p)}{\gcd(\pma[<x](p), \pma[<x](q))}\]
  and $\deg_x(\pma[<x](\widetilde{p})) = \deg_x(\pma[<x](p))-\deg_x(\gcd(\pma[<x](p), \pma[<x](q)))=\deg_x(p)-i$.\qed
 \end{proof}

Applying the last algorithm to a polynomial $p$ and $\frac{\partial}{\partial\ld(p)}p$ yields an algorithm to make $p$ square-free.
We present it separately for better readability of the main algorithm.

\begin{myalgorithm}[\sf ResSplitSquareFree]\label{algo_ressplitsquarefree}
 \textit{Input:} A system $S$ with $(S_Q)_{<x}^==\emptyset$ and a polynomial $p$ with $\ld(p)=x$ and $\pma(\ini(p))\neq0$ for all $\ma \in \sol(S_{<x})$.\\
 \textit{Output:} Two systems $S_1$ and $S_2$ and a polynomial $r$ such that:
  \renewcommand{\theenumi}{\alph{enumi}}
  \renewcommand{\labelenumi}{\theenumi)}
  \begin{enumerate} \label{algo_ressplitsquarefree_cond1}
   \item $S_2 = \widetilde{S_2} \cup \{ p \}$ where $\left(S_1,\widetilde{S_2}\right)$ is an $i$-th square-free split of $p$ w.r.t. $S$,
   \item \label{algo_ressplitsquarefree_cond2}$r$ is an $i$-th conditional square-free part of $p$ w.r.t. $S$,
  \end{enumerate}
  \renewcommand{\theenumi}{\arabic{enumi}}
  \renewcommand{\labelenumi}{\theenumi)}
  where $i$ is the quasi fiber cardinality of $p$ and $\frac{\partial}{\partial x}p$ w.r.t. $S$.\\
 \textit{Algorithm:} \begin{algorithmic}[1]
  \STATE $(i, S_1, S_2) \gets \textrm{\sf ResSplit}\left(S, p, \frac{\partial}{\partial x}p\right)$
  \IF{$i>0$}
    \STATE $r \gets \pquo\left(p, \prs_i\left(p, \frac{\partial}{\partial x}p, x\right),x\right)$
  \ELSE
    \STATE $r \gets p$
  \ENDIF
  \STATE $(S_2)_Q \gets (S_2)_Q \cup \{ p \}$
  \RETURN $S_1, S_2, r$
 \end{algorithmic}
\end{myalgorithm}
 \begin{proof}[Correctness]
  Since $\pma[<x](\frac{\partial}{\partial x}p)=\frac{\partial}{\partial x}\pma[<x](p)$, an $i$-th square-free split of $p$ is an $i$-th fibration split of $p$ and $\frac{\partial}{\partial x}p$.
  The rest follows from the proof of Algorithm (\ref{algo_ressplitdivide}).\qed
 \end{proof}

In all \textsf{ResSplit}-based algorithms, $(S_Q)_{<x}^==\emptyset$ is required.
This ensures that all equations of a smaller leader than $x$ will be respected by reduction modulo $S_T$.
The order in which polynomials are treated by the main algorithm must therefore be restricted.

\begin{mydefinition}[\sf Select]\label{Select}
 Let $\mathbb{P}_{\textrm{finite}}(M)$ be the set of all finite subsets of a set $M$. A \textbf{selection strategy} is a map \[\begin{array}{rcl} \textrm{\sf Select}: \mathbb{P}_{\textrm{finite}}(\{ p_=, p_{\neq} \mid p \in R \}) & \longrightarrow & \{ p_=, p_{\neq} \mid p \in R \}:\\ Q & \longmapsto & q \in Q\end{array}\] with the following properties: \begin{enumerate}
  \item\label{Select_Axiom1} If $\textrm{\sf Select}(Q) = q_=$ is an equation, then $Q_{<\ld(q)}^= = \emptyset$.
  \item\label{Select_Axiom2} If $\textrm{\sf Select}(Q) = q_{\neq}$ is an inequation, then $Q_{\leq\ld(q)}^= = \emptyset$.
 \end{enumerate}
\end{mydefinition}

We demonstrate that these conditions are necessary for termination of our approach, by giving an example where we violate them.

\begin{myexample}
 Consider $R := F[a,x]$ with $a<x$ and the system $S$ with $S_T := \emptyset$ and $S_Q := \left\{ (x^2-a)_= \right\}$.
 To insert $(x^2-a)_=$ into $S_T$, we need to apply the \textsf{ResSplitSquareFree} algorithm:
 We calculate $\res_0(x^2-a, 2x, x) = -4a$, $\res_1(x^2-a, 2x, x)=2$ and $\res_2(x^2-a, 2x, x)=1$ according to Definition (\ref{not_prs}).
 The quasi fiber cardinality is $0$ and we get the two new systems $S_1$, $S_2$ with \[(S_1)_T = \{ (x^2-a)_= \}, (S_1)_Q = \{ (-4a)_{\neq} \} \mbox{\ \ and\ \ } (S_2)_T = \emptyset, (S_2)_Q = \{ (x^2-a)_=, (-4a)_= \}\enspace.\]

 We now consider what happens with $S_2$: If we select $(x^2-a)_=$ as the next equation to be treated, in violation of the properties in Definition (\ref{Select}), \textsf{ResSplitSquareFree} will split up $S_2$ into $S_{2,1}$, $S_{2,2}$ with \[(S_{2,1})_T = \{ (x^2-a)_= \}, (S_{2,1})_Q = \{ (-4a)_{\neq}, (-4a)_= \}\] and \[(S_{2,2})_T = \emptyset, (S_{2,2})_Q = \{ (x^2-a)_=, (-4a)_=, (-4a)_= \}\enspace.\]
 As $S_2 = S_{2,2}$, this will lead to an endless loop.
\end{myexample}

The following trivial algorithm inserts a new equation into $S_T$.
It will be replaced with a different algorithm in \textsection\ref{section_differential} when the differential \textsc{Thomas} decomposition is considered.

\begin{myalgorithm}[\sf InsertEquation]
 \textit{Input:} A system $S$ and an equation $r_=$ with $\ld(r)=x$ satisfying $\pma(\ini(r))\neq0$ and $\pma[<x](r)$ square-free for all $\ma \in \sol(S_{<x})$. \\
 \textit{Output:} A system $S$ where $r_=$ is inserted into $S_T$. \\
 \textit{Algorithm:} \begin{algorithmic}[1]
  \IF{$(S_T)_x$ is not empty}
   \STATE $S_T \gets (S_T \setminus \{ (S_T)_x \})$
  \ENDIF
  \STATE $S_T \gets S_T \cup \{ r_= \}$
  \RETURN $S$
 \end{algorithmic}
\end{myalgorithm}

Now we present the main algorithm.
The general structure is as follows:
In each iteration, a system $S$ is selected from a list $P$ of unfinished systems.
An equation or inequation $q$ is chosen from the queue $S_Q$ according to the selection strategy.
Then $q$ is reduced modulo $S_T$ and incorporated into the candidate simple system $S_T$ with the splitting algorithms as described above.
In doing so, the algorithm may add new systems $S_i$ to $P$.
As soon as the algorithm produces a system containing an equation $c_=$ for $c\in F\setminus\{0\}$ or the inequation $0_{\neq}$, this system is discarded.

\begin{myalgorithm}[\sf Decompose]\label{algo_decompose} The algorithm is printed on page \pageref{algo_decompose_float}.
 \begin{Decompose}
 \label{algo_decompose_float}
 \begin{minipage}[top]{\textwidth}
 {\textbf{Algorithm \ref{algo_decompose} (\textsf{Decompose})}}\vspace{-5pt}\\
 \noindent\begin{tabular*}{\textwidth}{c}\hline\end{tabular*} 
 \textit{Input:} A system $S^\prime$ with $({S^\prime})_T=\emptyset$. \\
 \textit{Output:} A {\sc Thomas} decomposition of $S^\prime$. \\
 \textit{Algorithm:}
 \begin{algorithmic}[1]
  \STATE $P \gets \{ S^\prime \}$; $\mathit{Result} \gets \emptyset$
  \WHILE{$|P|>0$}\label{main_loop}
   \STATE \label{choose_S} Choose $S \in P$; $P \gets P \setminus \{ S \}$
   \IF{$|S_Q|=0$}\label{SQempty}
    \STATE $\mathit{Result} \gets \mathit{Result} \cup \{ S \}$
   \ELSE
    \STATE \label{remove_from_Q} $q \gets \textrm{\sf Select}(S_Q)$; $S_Q \gets S_Q \setminus \{ q \}$
    \STATE \label{reduce_q} $q \gets \textrm{\sf Reduce}(q, S_T)$; $x \gets \ld(q)$
    \IF{$q \notin \{ 0_{\neq}, c_{=} \ |\ c \in F \setminus \{0\} \}$}\label{omit_sys}
    \IF{$x \neq 1$}
    \IF{$q$ is an equation}
     \IF{$(S_T)_x$ is an equation}
      \IF{$\textrm{\sf Reduce}(\res_0((S_T)_x, q, x), S_T)=0$}\label{reduceres}
       \STATE \label{resplitgcd} $(S, S_1, p) \gets \textrm{\sf ResSplitGCD}(S, q, x)$; $P \gets P \cup \{ S_1 \}$
       \STATE \label{q_1}$S \gets \textrm{\sf InsertEquation}(S, p_=)$
      \ELSE
       \STATE \label{insertres}$S_Q \gets S_Q \cup \{ q_=, \res_0((S_T)_x, q, x)_= \}$
      \ENDIF
     \ELSE
      \IF{$(S_T)_x$ is an inequation\footnote{Remember that $(S_T)_x$ might be empty, and thus neither an equation nor an inequation.}}\label{begin_eqineq}
       \STATE \label{removeineq} $S_Q \gets S_Q \cup \{ (S_T)_x \}$; $S_T \gets S_T \setminus \{ (S_T)_x \}$
      \ENDIF
      \STATE $(S, S_2) \gets \textrm{\sf InitSplit}(S, q)$; $P \gets P \cup \left\{ S_2 \right\}$
      \STATE $(S,S_3,p) \gets \textrm{\sf ResSplitSquareFree}\left(S, q\right)$; $P \gets P \cup \{ S_3 \}$
      \STATE \label{q_2}$S \gets \textrm{\sf InsertEquation}(S, p_=)$
     \ENDIF
    \ELSIF{$q$ is an inequation}
     \IF{$(S_T)_x$ is an equation}
      \STATE \label{dividebyinequation} $(S, S_4, p) \gets \textrm{\sf ResSplitDivide}\left(S, (S_T)_x, q\right)$; $P \gets P \cup \{ S_4 \}$
      \STATE \label{q_3}$S \gets \textrm{\sf InsertEquation}(S, p_=)$
     \ELSE\label{begin_ineqineq}
      \STATE $(S, S_5) \gets \textrm{\sf InitSplit}(S, q)$; $P \gets P \cup \{ S_5 \}$
      \STATE $(S, S_6, p) \gets \textrm{\sf ResSplitSquareFree}\left(S, q\right)$; $P \gets P \cup \{ S_6 \}$
      \IF{$(S_T)_x$ is an inequation}
       \STATE $(S, S_7, r) \gets \textrm{\sf ResSplitDivide}\left(S, (S_T)_x, p \right)$; $P \gets P \cup \{ S_7 \}$
       \STATE \label{q_4}$(S_T)_x \gets (r\cdot p)_{\neq}$
      \ELSIF{$(S_T)_x$ is empty}
       \STATE \label{q_5}$(S_T)_x \gets p_{\neq}$
      \ENDIF\label{end_ineqineq}
     \ENDIF
    \ENDIF
    \ENDIF
    \STATE \label{add_S}$P \gets P \cup \{ S \}$
    \ENDIF
   \ENDIF
  \ENDWHILE
  \RETURN $\mathit{Result}$
 \end{algorithmic}
 \end{minipage}
 \end{Decompose}
\end{myalgorithm}

We demonstrate the algorithm with a simple example.
Note, that we will omit systems which are obviously inconsistent.

\begin{myexample}
 Let $S = (S_T, S_Q) := (\emptyset, \{ (x^2+x+1)_=, (x+a)_{\neq} \})$ with $a<x$.
 According to \textsf{Select}, $q := (x^2+x+1)_=$ is chosen.
 As $\ini(q)=1$ and $\res_0(q, \frac{\partial}{\partial x}q, x)=1$, the original system $S$ is replaced by $\left(\{(x^2+x+1)_=\}, \{(x+a)_{\neq}\}\right)$.

 Now, $q := (x+a)_{\neq}$ is selected and $\mathsf{ResSplitDivide}(S, (S_T)_x, q)$ computes $\res_0((S_T)_x, q, x) = \prem((S_T)_x, q, x) = a^2-a+1$, $\res_1((S_T)_x, q, x)=\ini(q)=1$, and $\res_2((S_T)_x, q, x)=1$.
 As $S_T$ contains no equation of leader $a$, none of these polynomials can be reduced.
 Then, we decompose $S$ into
  \[S := (\underbrace{\{ (x^2+x+1)_=, (a^2-a+1)_{\neq} \}}_{= S_T}, \underbrace{\{\}}_{=S_Q})\mbox{ ,}\]
 which is already simple, and
  \[S_1 := ( \underbrace{\{ (x^2+x+1)_= \}}_{=(S_1)_T}, \underbrace{\{ (x+a)_{\neq}, (a^2-a+1)_= \}}_{=(S_1)_Q} )\enspace.\]
 We replace $S_1$ by \[S_1 := \left( \{ (x^2+x+1)_=, (a^2-a+1)_= \}, \{ (x+a)_{\neq} \} \right)\]
 and apply $\mathsf{ResSplitDivide}(S_1, ((S_1)_T)_x, q)$ to $S_1$ again.
 This time, $\mathsf{Reduce}(a^2-a+1, (S_1)_T)=0$ holds and $S_1$ is replaced with \[S_1 := ( \{ \underbrace{(x-a+1)_=}_{\pquo(x^2+x+1, x+a, x)}, (a^2-a+1)_= \}, \{ 1_{\neq} \} )\enspace.\]
 Finally, a {\sc Thomas} decomposition of $S$ is:
 \[\left(\{ (x^2+x+1)_=, (a^2-a+1)_{\neq} \}, \{ (x-a+1)_=, (a^2-a+1)_= \} \right)\enspace.\]
\end{myexample}

 \begin{proof}[Correctness]
  First, note that it is easily verified that the input specifications of all subalgorithms are fulfilled (in particular, for lines \ref{resplitgcd} and \ref{dividebyinequation}, cf. Remark (\ref{rem_reduceprop})(\ref{rem_reduceprop1})).
  
  The correctness of the \textsf{Decompose} algorithm is proved by verifying two loop invariants:
  \begin{enumerate}
   \item $P \cup \mathit{Result}$ is a disjoint decomposition of the input $S^\prime$.
   \item For all systems $S \in P \cup \mathit{Result}$, $S_T$ is triangular and
   \begin{enumerate}
     \item $\pma[<x](p)$ is square-free and 
     \item $\pma(\ini(p))\neq0$ 
   \end{enumerate}
   for all $p \in S_T$ with $\ld(p)=x$ and all $\ma \in \sol((S_T)_{<x} \cup (S_Q)_{<x})$.
  \end{enumerate}

  We begin with proving the first loop invariant.
  Assume that $P \cup \mathit{Result}$ is a disjoint decomposition of $S^\prime$ at the beginning of the main loop.
  It suffices to show that all systems we add to $P$ or $\mathit{Result}$ add up to a disjoint decomposition of the system $S$, that is chosen in line \ref{choose_S}.
  If $S_Q=\emptyset$ holds in line \ref{SQempty}, the algorithm just moves $S$ from $P$ to $\mathit{Result}$.

  In line \ref{insertres}, adding $\res_0((S_T)_x,q,x)_=$ to $S$ does not change the solutions of $S$, as $\pma[<x]((S_T)_x)=0$ and $\pma[<x](q)=0$ for each $\ma \in F[\ y\mid y<x\ ]$ implies $\pma(\res_0(p,q,x))=0$ (cf.\  \cite[Lemma 7.2.3]{mishra}).

  Note now that if $(S, S_i)$ is the output of any of the {\sf ResSplitGcd}, {\sf InitSplit}, {\sf ResSplitSquareFree} and {\sf ResSplitDivide} algorithms, then $(S \cup \{q\}, S_i)$ is a disjoint decomposition of $S_0 \cup \{q\}$, where $S_0$ is the input of the respective algorithm.
  It remains to be shown that the actions performed in lines \ref{q_1}, \ref{q_2}, \ref{q_3}, \ref{q_4} and \ref{q_5} are equivalent to putting $q$ back into the system $S$.

  Let $\ma \in \sol(S_{<x})$. In the context of line \ref{q_1}, Algorithm (\ref{algo_ressplitgcd}) guarantees
   \[\pma[<x](p)=0 \Longleftrightarrow \pma[<x]((S_T)_x)=0 \textrm{ and } \pma[<x](q)=0 \mbox{ .}\]
  In the context of line \ref{q_3}, Algorithm (\ref{algo_ressplitdivide}) ensures that
   \[\pma[<x](p)=0 \Longleftrightarrow \pma[<x]((S_T)_x)=0 \textrm{ and } \pma[<x](q)\neq0\enspace.\]
  In lines \ref{q_2}, \ref{q_4} and \ref{q_5}, $p$ has the same solutions as $q$, due to Algorithm (\ref{algo_ressplitsquarefree}) and
   \[\pma[<x](p) \sim \frac{\pma[<x](q)}{\gcd(\pma[<x](q), \pma[<x](\frac{\partial}{\partial x}q))} = \frac{\pma[<x](q)}{\gcd(\pma[<x](q), \frac{\partial}{\partial x}\pma[<x](q))}\enspace.\]
  In addition, in line \ref{q_4},
   \[\pma[<x](r) \sim \frac{\pma[<x]((S_T)_x)}{\gcd(\pma[<x]((S_T)_x), \pma[<x](p))}
  \Longrightarrow
   \pma[<x](r\cdot p) \sim \lcm(\pma[<x]((S_T)_x), \pma[<x](p))\enspace.\]
  This concludes the proof of the first loop invariant.

  Now, we prove the second loop invariant.
  At the beginning, the loop invariant holds because $S^\prime_T=\emptyset$ holds for the input system $S^\prime$.
  Assume that the second loop invariant holds at the beginning of the main loop.
  
  One easily checks that all steps in the algorithm allow only one polynomial $(S_T)_x$ in $S_T$ for each leader $x$, thus triangularity obviously holds.

  We show that all polynomials added to $S_T$ have non-zero initial and are square-free.
  For $\sol(S_{<x})=\emptyset$, the statement is trivially true.
  So, let $\ma \in \sol(S_{<x})$.

  For the equation $p_=$ added as conditional gcd of $(S_T)_x$ and $q$ in line \ref{q_1}, it holds that $\pma[<x](p)$ is a divisor of $\pma[<x]((S_T)_x)$. As $\pma[<x]((S_T)_x)$ is square-free by assumption, so is $\pma[<x](p)$.
  The inequation added to $S$ in {\sf ResSplitGCD} is by Definition (\ref{not_prs}) the initial of $p_=$.

  The equation $p_=$ inserted into $S_T$ in line \ref{q_2} and the inequation $p_{\neq}$ inserted in line \ref{q_5} are square-free due to Algorithm (\ref{algo_ressplitsquarefree}) and their initials are non-zero as $p$ is either identical to $q$, or it is a pseudo quotient of $q$ by $\prs_i\left(q,\frac{\partial}{\partial x}q, x\right)$ for some $i>0$.
  On the one hand, if $p$ equals $q$, the call of \textsf{InitSplit} for $q$ ensures a non-zero initial for $p$.
  On the other hand, the polynomial $\prs_i\left(q,\frac{\partial}{\partial x}q, x\right)$ has initial $\res_i\left(q,\frac{\partial}{\partial x}q, x\right)$, which is added as an inequation by {\sf ResSplitSquareFree}.
  This implies that the initial of the pseudo-quotient is also non-zero.

  The equation $p_=$ that replaces the old equation $(S_T)_x$ in line \ref{q_3} is the quotient of $(S_T)_x$ by an inequation.
  It is square-free, because $\pma[<x](p)$ is a divisor of $\pma[<x]((S_T)_x)$, which is square-free by assumption.
  Again, $p$ is either identical to $(S_T)_x$ or a pseudo quotient of $(S_T)_x$ by $\prs_i\left((S_T)_x, q, x\right)$ for some $i>0$ and, using the same arguments as in the last paragraph, the initial of $p$ does not vanish.

  Finally, consider the inequation $(r\cdot p)_{\neq}$ added in line \ref{q_4} as a least common multiple of $\left((S_T)_x\right)_{\not=}$ and $p_{\not=}$.
  The inequation $\pma[<x](p)$ is square-free and has non-vanishing initial for the same reasons as before.
  Due to $\pma[<x](r) \sim \frac{\pma[<x]((S_T)_x)}{\gcd(\pma[<x]((S_T)_x), \pma[<x](p))}$, the polynomials $\pma[<x](r)$ and $\pma[<y](p)$ have no common divisors.
  As $\pma[<x](r)$ divides $\pma[<x]((S_T)_x)$, using the same arguments as before, $\pma[<x](r)$ is square-free and has a non-vanishing initial.
  This completes the proof of the second loop invariant.

  It is obvious that a system $S$ with $S_Q=\emptyset$ for which these loop invariants hold is simple. Thus the algorithm returns the correct result if it terminates.\qed
 \end{proof}

We now start showing termination.
The system $S$ chosen from $P$ is treated in one of three ways:
It is either discarded, added to $\mathit{Result}$, or replaced in $P$ by at least one new system.
To show that $P$ is empty after finitely many iterations, we define an order on the systems and show that it is well-founded.
Afterwards we prove termination by detailing that the algorithm produces descending chains of systems.

\begin{mydefinition}\label{def_comp_order}
For transitive and asymmetric\footnote{A relation $\prec$ is asymmetric, if $S\prec S'$ implies $S'\not\prec S$ for all $S,S'$. Asymmetry implies irreflexivity.} partial orders $<_i$ for $i=1,\ldots,m$, we define the \textbf{composite order} $\mbox{``}<\mbox{''}:=[<_1,\ldots,<_m]$ as follows:
$a<b$ if and only if there exists $i \in \{ 1, \ldots, m \}$ such that $a<_ib$ and neither $a<_jb$ nor $b<_ja$ for $j<i$.
The composite order is clearly transitive and asymmetric.
An order $<$ is called \textbf{well-founded}, if each $<$-descending chain becomes stationary.
\end{mydefinition}

The following trivial statement will be used repeatedly:
  
\begin{myremark}\label{composite_well_founded}
If each $<_i$ is well-founded, then so is the composite ordering $<$, using the notation from Definition (\ref{def_comp_order}).
\end{myremark}

Now we define the orders and show their well-foundedness:
  
\begin{defrem}\label{def_orders}
   
  Define $\prec$ as the composite order $[\prec_1,\prec_2,\prec_3,\prec_4]$ of the four orders defined below. It is well-founded since the $\prec_i$ are.

   \begin{enumerate}

   \item \label{order_1}For $i=1,\ldots,n$ define $\prec_{1,x_i}$ by $S\prec_{1,x_i} S'$ if and only if $\mdeg\left((S_T)^=_{x_i}\right) < \mdeg\left((S^\prime_T)^=_{x_i}\right)$, with $\mdeg\left((S_T)^=_{x_i}\right):=\infty$ if $(S_T)^=_{x_i}$ is empty.
   Define the composite order $\prec_1$ as $[\prec_{1,x_1},\ldots,\prec_{1,x_n}]$.
   Since degrees can only decrease finitely many times, the orders $\prec_{1,x_i}$ are clearly well-founded and, thus, $\prec_1$ is.

   \item \label{order_2}Define the map $\mu$ from the set of all systems over $R$ to $\{1,x_1,\ldots,x_n,x_\infty\}$, where $\mu(S)$ is minimal such that there exists an equation $p \in (S_Q)^=_{\mu(S)}$ with $\mathsf{Reduce}(S_T, p)\neq0$, or $\mu(S)=x_\infty$ if no such equation exists. Then, $S \prec_2 S^\prime$ if and only if $\mu(S) < \mu(S^\prime)$ with $1<x_i$ and $x_i<x_\infty$ for $i\in\{1,\ldots,n\}$.
   The ordering $\prec_2$ is well-founded since $<$ is well-founded on the finite set $\{1,x_1, \ldots, x_n,x_\infty\}$.

  \item \label{order_3}$S \prec_3 S'$ if and only if there is $p_{\not=}\in R^{\not=}$ and a finite (possibly empty) set $L\subset R^{\not=}$ with $\ld(q)<\ld(p)\ \forall\ q\in L$ such that $S_Q \uplus \{p_{\not=}\} = S'_Q \uplus L$ holds.
  We show well-foundedness by induction on the highest appearing leader $x$ in $(S_Q)^{\not=}$:
  For $x=1$ we can only make a system $S$ $\prec_3$-smaller by removing one of the finitely many inequations in $(S_Q)^{\not=}$.
  Now assume that the statement is true for all indeterminates $y<x$.
  By the induction hypothesis we can only $\prec_3$-decrease $S$ finitely many times without changing $(S_Q)_x^{\not=}$.
  To further $\prec_3$-decrease $S$, we have to remove an inequation from $(S_Q)^{\not=}_x$.
  As $(S_Q)^{\neq}_x$ is finite, this process can only be repeated finitely many times until $(S_Q)^{\neq}_x=\emptyset$. Now, the highest appearing leader in $(S_Q)^{\neq}$ is smaller than $x$ and by the induction hypothesis, the statement is proved.

   \item \label{order_4}$S \prec_4 S^\prime$ if and only if $|S_Q| < |S^\prime_Q|$.
  \end{enumerate}
  
\end{defrem}

 \begin{proof}[Termination]

  We will tacitly use the fact that reduction never makes polynomials bigger in the sense of Remark (\ref{rem_reduceprop})(\ref{rem_reduceprop3}).
  
  We denote the system chosen from $P$ in line \ref{choose_S} by $\widehat{S}$ and the system added to $P$ in line \ref{add_S} by $S$.
  We prove that the systems $S, S_1, \ldots, S_7$ generated from $\widehat{S}$ are $\prec$-smaller than $\widehat{S}$.
  For $i=1,\ldots,4$ we will use the notation $S \osim_i S^\prime$ if neither $S \prec_i S^\prime$ nor $S^\prime \prec_i S$ holds.

  For $j=1,\ldots,7$, $((S_j)_T)^= = (\widehat{S}_T)^=$ and thus $S_j \osim_1 \widehat{S}$.
  The properties of {\sf Select} in Definition (\ref{Select}) directly require, that there is no equation in $(\widehat{S}_Q)^=$ with a leader smaller than $x$.
  However, the equation added to the system $S_j$ returned from {\sf InitSplit} (\ref{algo_initsplit}) is the initial of $q$, which has a leader smaller than $x$ and does not reduce to $0$ (cf.~Remark (\ref{rem_reduceprop})(\ref{rem_reduceprop2})).
  Furthermore, the equations added in one of the subalgorithms based on {\sf ResSplit} (\ref{algo_ressplit}) have a leader smaller than $x$ and do not reduce to $0$.
  In each case $S_j \prec_2 \widehat{S}$ is proved.

  It remains to show $S \prec \widehat{S}$.
  If $q$ is reduced to $0_=$, then it is omitted from $S_Q$ and so $S\prec_4\widehat{S}$.
  As the system is otherwise unchanged, $S\osim_i\widehat{S}, i=1,2,3$ and therefore $S\prec\widehat{S}$ holds.
  If $q$ is reduced to $c_{\neq}$ for some $c \in F \setminus \{0\}$, then $S \prec_3 \widehat{S}$ and $S \osim_i \widehat{S}, i=1,2$, since the only change was the removal of an inequation from $S_Q$.
  Otherwise, exactly one of the following cases will occur:

  Lines \ref{resplitgcd}-\ref{q_1} set $(S_T)_x$ to $p_=$ of smaller degree than $(\widehat{S}_T)_x$ and \ref{begin_eqineq}-\ref{q_2} add $(S_T)_x$ as a new equation.
  In both cases we get $S \prec_1 \widehat{S}$.
  
  In line \ref{insertres}, $S_T=\widehat{S}_T$ implies $S \osim_1 \widehat{S}$.
  The polynomial $q$ is chosen according to {\sf Select} (cf.\ (\ref{Select})(\ref{Select_Axiom1})), which implies $(\widehat{S}_Q)^=_{<x}=\emptyset$ and $(S_Q)^=_{<x}=\{\res_0((S_T)_x, q, x)_=\}$.
  Line \ref{reduceres} ensures $\mbox{\sf Reduce}(\res_0((S_T)_x, q, x), S) \neq 0$ and, thus, $S \prec_2 \widehat{S}$ follows.
  
  Consider lines \ref{dividebyinequation}-\ref{q_3}.
  If the degree of $(S_T)_x$ is smaller than the degree of $(\widehat{S}_T)_x$, then $S \prec_1 \widehat{S}$.
  In case the degree doesn't change, we have $S \osim_1 \widehat{S}$ and $(S_Q)^==(\widehat S_Q)^=$ guarantees $S \osim_2 \widehat{S}$.
  However, $q$ is removed from $S_Q$ and replaced by an inequation of smaller leader, which implies $S \prec_3 \widehat{S}$.
  
  In \ref{begin_ineqineq}-\ref{end_ineqineq}, obviously $S \osim_i \widehat{S}, i=1,2$.
  As before, $q$ is removed from $S_Q$ and replaced by an inequation of smaller leader, which once more implies $S \prec_3 \widehat{S}$.\qed
 \end{proof}

\subsection{Notes to Applications of Simple Systems}\label{sect_compare_simple_and_regular}

In this subsection, we shortly present some examples where simple systems are necessary and any weaker decomposition into triangular systems is not sufficient.

The properties of simple systems (cf.\ Definition (\ref{simple_system})) correspond exactly to the following fibration structure on the solution sets (cf.~\cite{PleskenCounting}).
Let $S$ be a simple system and $\Pi_i: \overline{F}^i \to \overline{F}^{i-1}: (a_1,\ldots,a_i) \mapsto (a_1,\ldots,a_{i-1})$.
Furthermore, for any solution $\mathbf{a} \in \sol(S_{\leq x_i})$, let $s_{i,\mathbf{a}} = {\Pi_i}^{-1}(\{ \Pi_i(\mathbf{a}) \})$.
Then, if $S_{x_i}$ is an equation, $|s_{i,\mathbf{a}}|=\mdeg(S_{x_i})$ holds.
If $S_{x_i}$ is an inequation, then $s_{i,\mathbf{a}}=\overline{F} \setminus \tilde{s}_{i,\mathbf{a}}$ with $|\tilde{s}_{i,\mathbf{a}}|=\mdeg(S_{x_i})$.
If $S_{x_i}$ is empty, then $s_{i,\mathbf{a}}=\overline{F}$.
The cardinalities of $s_{i,\mathbf{a}}$ or $\tilde{s}_{i,\mathbf{a}}$ are constant for each $i$, i.e. independent of the choice of the solution $\mathbf{a} \in \sol(S_{\leq x_i})$ (cf.\ Remark (\ref{exist_sol})).
We can examine solution sets of arbitrary systems by decomposing them disjointly into simple systems.
Further analysis of this fibration structure, especially in the context of algebraic varieties, is a topic of future research.

We already saw such a fibration structure in Example (\ref{ex_elliptic_curve}).
In this case, other triangular decompositons like a decomposition into regular chains would have only resulted in a single system consisting of the polynomial $p$ from the input.

A special case occurs when all polynomials in the input and output can be factored into linear polynomials.
If we compute the counting polynomial as introduced by \cite{PleskenCounting} (which requires the disjointness of the decomposition and the fibration structure), we can substitute the cardinality of a finite field $F$ (of sufficiently large characteristic) into the counting polynomial of a \textsc{Thomas} decomposition computed over $\mathbb{Q}$.
This yields the exact number of distinct solutions over $F$.
For example, the counting polynomial of a \textsc{Thomas} decomposition of $\{\det(M)_{\neq}\}$ for a generic $n\times n$ matrix $M=(x_{ij})_{1\le i,j\le n}$ yields the well-known formula for the cardinality of $\operatorname{GL}_n(F)$ for any finite field $F$.
Furthermore, we can automatically reproduce the results in \cite[Ex.~V.4]{PleskenCountingGroupsAndRings}, where pairs of matrices $(A,B)$ with given ranks of $A$, $B$, and $A+B$ are counted.

\cite{PleskenBruhat} gave another example concerning the \textsc{Gauss}-\textsc{Bruhat}-decomposition and the LU-decomposition.
The cells of these decompositions of $M$ as above can be identified with certain simple systems in the \textsc{Thomas} decomposition of $\{\det(M)_{\neq}\}$ for suitable rankings on the $x_{ij}$.

We clearly see that simple systems are necessary for these applications to expose the aforementioned fibration structure and count solutions. 
A disjoint decomposition into triangular systems with weaker properties does not suffice.

\section{Differential Thomas Decomposition}\label{section_differential}

The differential \textsc{Thomas} decomposition is concerned with manipulations of polynomial differential equations and inequations.
The basic idea for our construction of this decomposition is twofold.
On the one hand, a combinatorial calculus developed by {\sc Janet} finds unique reductors and all integrability conditions by completing systems to involution.
On the other hand, the algebraic {\sc Thomas} decomposition makes the necessary splits for regularity of initials and ensures disjointness of the solution sets.

Initially, we recall some basic definitions from differential algebra.
Then, we summarize the {\sc Janet} division and its relevance.
Its combinatorics lead us to substitute the algebraic algorithm {\sf InsertEquation} by its differential analog.
Afterwards, we review a differential generalization of the algebraic reduction algorithm and present the algorithm {\sf Reduce} utilized for differential reduction.
Replacing the insertion and reduction from the previous section with these differential counterparts yields the differential \textsc{Thomas} decomposition algorithm.

\subsection{Preliminaries from Differential Algebra}\label{differential_preliminaries}

Let $\Delta=\{\partial_1,\ldots,\partial_n\}$ be a non-empty set of derivations and $F$ be a $\Delta$\textbf{-ring}.
This means any $\partial_j\in\Delta$ is a linear operator $\partial_j:F\to F$ which satisfies the \textsc{Leibniz} rule.
Given a \textbf{differential indeterminate} $u$, the \textbf{polynomial $\Delta$-ring} $F\{u\}:=F\left\lbrack\ u_{\mathbf{i}}\mid\mathbf{i}\in\Z_{\ge0}^n\ \right\rbrack$ is defined as the polynomial ring infinitely generated by the algebraically independent set $\langle u\rangle_\Delta:=\{u_{\mathbf{i}}\mid\mathbf{i}\in\Z_{\ge0}^n\}$.
The operation of $\partial_j\in\Delta$ on $\langle u\rangle_\Delta$ is defined by $\partial_j u_{\mathbf{i}}=u_{\mathbf{i}+e_j}$ and this operation extends linearly and via the \textsc{Leibniz} rule to $F\{u\}$.
Let $U=\{u^{(1)},\ldots,u^{(m)}\}$ be a set of differential indeterminates.
The multivariate polynomial $\Delta$-ring is given by $F\{U\}:=F\{u^{(1)}\}\ldots\{u^{(m)}\}$.
Its generators, the elements of $\langle U\rangle_\Delta:=\left\{u^{(j)}_{\mathbf{i}}\mid\mathbf{i}\in\Z_{\ge0}^n, j\in\{1,\ldots,m\}\right\}$, are called \textbf{differential variables}.
From now on let $F$ be a computable $\Delta$\textbf{-field} of characteristic zero.

The differential structure of $F$ uniquely extends to the differential structure of its algebraic closure $\overline{F}$ \cite[\textsection II.2, Lemma 1]{Kol}.
Let $E:=\bigoplus_{j=1}^m\overline{F}[[z_1,\ldots,z_n]]$ where $\overline{F}[[z_1,\ldots,z_n]]$ denotes the ring of formal power series in $z_1,\ldots,z_n$.
Then $E$ is isomorphic to $\overline{F}^{\langle U\rangle_\Delta}$ via
\[
  \alpha:
  \bigoplus_{j=1}^m\overline{F}[[z_1,\ldots,z_n]]
  \to\overline{F}^{\langle U\rangle_\Delta}:
  \left(
    \sum_{\mathbf{i}\in\Z_{\ge0}^n}a_{\mathbf{i}}^{(1)}\frac{z^{\mathbf{i}}}{\mathbf{i}!},
    \ldots,
    \sum_{\mathbf{i}\in\Z_{\ge0}^n}a_{\mathbf{i}}^{(m)}\frac{z^{\mathbf{i}}}{\mathbf{i}!}
  \right)
  \mapsto\left(u^{(j)}_{\mathbf{i}}\mapsto a_{\mathbf{i}}^{(j)}\right)
\]
where $z^\mathbf{i}:=z_1^{i_1}\cdot\ldots\cdot z_n^{i_n}$ and $\mathbf{i}!:=i_1!\cdot\ldots\cdot i_n!$.

We define solutions in $E$, consistent with the algebraic case:
For $e\in E$, let
\[
  \phi_{e}: F\{U\}\to\overline{F}: u^{(j)}_{\mathbf{i}}\mapsto\alpha(e)(u^{(j)}_{\mathbf{i}})
\]
be the $F$-algebra homomorphism evaluating the differential variables at $e$.
A \textbf{differential equation} or \textbf{inequation} for $m$ functions $U=\{u^{(1)},\ldots,u^{(m)}\}$ in $n$ indeterminates is an element $p\in F\{U\}$, written $p_=$ or $p_{\not=}$, respectively.
A \textbf{solution} of $p_=$ or $p_{\not=}$ is an $e\in E$ with $\phi_e(\langle p\rangle_\Delta)=\{0\}$ or $\phi_e(\langle p\rangle_\Delta)\not=\{0\}$, respectively.
Here, $\langle p\rangle_\Delta$ denotes the differential ideal in $F\{U\}$ generated by $p$.
Furthermore, $e\in E$ is called a solution of a set $P$ of equations and inequations, if it is a solution of each element in $P$.
The set of solutions of $P$ is denoted by $\sol(P):=\sol_E(P)\subseteq E$.

In differential algebra one usually considers solutions in a universal $\Delta$-field, while we consider power series solutions.
As the universal differential field we can take the universal closure $\widehat{F}$ of $\overline{F}$.
There is a strong link between these two concepts.
On the one hand, \cite{Seidenberg58,Seidenberg69} has shown that every finitely differentially generated differential field is differentially isomorphic to a differential field of meromorphic functions in $n$ variables.
On the other hand, $\overline{F}[[z_1,\ldots,z_n]] \hookrightarrow \overline{F}((z_1,\ldots,z_n)) \hookrightarrow \widehat{F}$.
Here, the first map is the natural embedding into the quotient field and the second is an embedding given by the definition of the universal $\Delta$-field \cite[\textsection II.2 and \textsection III.7]{Kol}, as $\overline{F}((z_1,\ldots,z_n))$ is a finitely generated $\Delta$-field extension of $\overline{F}$.
Thus, any power series solution can be considered as a solution in the universal differential field.

A finite set of equations and inequations is called a \textbf{(differential) system} over $F\{U\}$.
We will be using the same notation for systems as in the algebraic {\sc Thomas} decomposition introduced in \textsection\ref{algebraic_definition_notation} and \textsection\ref{algebraic_algorithms}, in particular a system $S$ is represented by a pair $(S_T,S_Q)$.
However, the candidate simple system $S_T$ will also reflect a differential structure based on the combinatorics from the following section.

\subsection{Janet Division}\label{Janet}

In this subsection we will focus on a combinatorial approach called \textsc{Janet} division (cf.\ \cite{GB1}).
It manages the infinite set of differential variables and guarantees inclusion of all integrability conditions in a differential system.
For this purpose, it partitions the set of differential variables into ``free'' variables and finitely many ``cones'' of dependent variables.
We present an algorithm for inserting new equations into an existing set of equations and adjusting this cone decomposition accordingly.
An overview of modern development on \textsc{Janet} division can be found in \cite{GerI} and \cite{Seiler} and the original ideas were formulated by \cite{Janet}.

A (differential) \textbf{ranking} $<$ is defined as a total order on the differential variables and $1$ with $1<u\ \forall\ u\in U$, such that
\begin{enumerate}
  \item $u<\partial_j u$ and
  \item $u<v$ implies $\partial_j u<\partial_j v$\label{rankingproperty2}
\end{enumerate}
for all $u,v\in \langle U\rangle_\Delta\mbox{, }\partial_j\in\Delta$.
From now on let $<$ be an arbitrary and fixed differential ranking.
For any finite set of differential variables, a differential ranking induces a ranking as defined for the algebraic case in \textsection\ref{algebraic_definition_notation}.
Thereby, in accordance to the algebraic part, define the largest differential variable $\ld(p)$ appearing in a differential polynomial $p\in F\{U\}$ as \textbf{leader}, which is set to $1$ for $p\in F$.
Furthermore, define $\mdeg(p)$ and $\ini(p)$ as the degree in the leader and the coefficient of $\ld(p)^{\mdeg(p)}$, respectively.

\begin{myexample}\label{example_janet_combinatoric}

  Consider two derivations $\Delta=\{\partial_x, \partial_t\}$ and one differential indeterminate $u$.

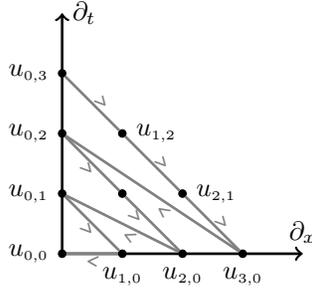
  \parpic[l]{
    \begin{tikzpicture}[domain=0:4, smooth, x=\picturesize, y=\picturesize, place/.style={circle,draw=black,fill=black,inner sep=1pt}]
      \draw[color=black, line width=1pt, ->] (0.0,0.0) -- (0.0,4.0) node[right]{$\partial_{t}$};
      \draw[color=black, line width=1pt, ->] (0.0,0.0) -- (4.0,0.0) node[above]{$\partial_{x}$};
      \node[place,label=180:$u\der{0,0}$] (u00) at (0.0,0.0) {};
      \node[place,label=-90:$u\der{1,0}$] (u10) at (1.0,0.0) {};
      \node[place,label=180:$u\der{0,1}$] (u01) at (0.0,1.0) {};
      \node[place,label=-90:$u\der{2,0}$] (u20) at (2.0,0.0) {};
      \node[place,label=180:{
      }] (u11) at (1.0,1.0) {};
      \node[place,label=180:$u\der{0,2}$] (u02) at (0.0,2.0) {};
      \node[place,label=-90:$u\der{3,0}$] (u30) at (3.0,0.0) {};
      \node[place,label=0:$u\der{2,1}$] (u21) at (2.0,1.0) {};
      \node[place,label=0:$u\der{1,2}$] (u12) at (1.0,2.0) {};
      \node[place,label=180:$u\der{0,3}$] (u03) at (0.0,3.0) {};
      \draw[color=gray, line width=1pt] (u00)--(u10);\node (00mid10) at ($ (u00)!.5!(u10) $) {};\node[below=-0.27 of 00mid10] (00smaller10) {\tiny{\color{gray}\textsf{<}}};
      \draw[color=gray, line width=1pt] (u10)--(u01);\node (10mid01) at ($ (u10)!.4!(u01) $) {};\node[above=0.005 of 10mid01,rotate=135] (10smaller01) {\tiny{\color{gray}\textsf{<}}};
      \draw[color=gray, line width=1pt] (u01)--(u20);\node (01mid20) at ($ (u01)!.6!(u20) $) {};\node[below=-0.27 of 01mid20,rotate=-22.5] (01smaller20) {\tiny{\color{gray}\textsf{<}}};
      \draw[color=gray, line width=1pt] (u20)--(u11);\node (20mid11) at ($ (u20)!.6!(u11) $) {};\node[above=+0.02 of 20mid11,rotate=135] (20smaller11) {\tiny{\color{gray}\textsf{<}}};
      \draw[color=gray, line width=1pt] (u11)--(u02);\node (11mid02) at ($ (u11)!.4!(u02) $) {};\node[above=+0.02 of 11mid02,rotate=135] (11smaller02) {\tiny{\color{gray}\textsf{<}}};
      \draw[color=gray, line width=1pt] (u02)--(u30);\node (02mid30) at ($ (u02)!.6!(u30) $) {};\node[below=-0.27 of 02mid30,rotate=-33.69] (02smaller30) {\tiny{\color{gray}\textsf{<}}};
      \draw[color=gray, line width=1pt] (u30)--(u21);\node (30mid21) at ($ (u30)!.5!(u21) $) {};\node[below=-0.00 of 30mid21,rotate=135] (30smaller21) {\tiny{\color{gray}\textsf{<}}};
      \draw[color=gray, line width=1pt] (u21)--(u12);\node (21mid12) at ($ (u21)!.5!(u12) $) {};\node[below=-0.00 of 21mid12,rotate=135] (21smaller12) {\tiny{\color{gray}\textsf{<}}};
      \draw[color=gray, line width=1pt] (u12)--(u03);\node (12mid03) at ($ (u12)!.5!(u03) $) {};\node[below=-0.00 of 12mid03,rotate=135] (12smaller03) {\tiny{\color{gray}\textsf{<}}};
    \end{tikzpicture}
  }

  In this setting, any partial differential equation with constant coefficients in one dependent variable and two independent variables can be represented as a differential polynomial in $\C\{u\}$.

  The ranking $<$ is defined by $u_{i_1,i_2}<u_{j_1,j_2}$ if and only if either $i_1+i_2<j_1+j_2$ or $i_1+i_2=j_1+j_2$ and $i_2<j_2$ holds.
  Thus, the smallest differential variables are:
  $u\der{0,0}<u\der{1,0}<u\der{0,1}<u\der{2,0}<u\der{1,1}<u\der{0,2}<u\der{3,0}$.
  When considering the set of differential variables as a grid in the first quadrant of a plane, the picture on the left illustrates this ranking.
  
  Consider $(\underline{u\der{0,1}}+u\der{0,0}u\der{1,0})_=$ representing the inviscid \textsc{Burger}'s equation $\frac{\partial u}{\partial t} + u \frac{\partial u}{\partial x} = 0$.
  As in the algebraic part, we indicate an equation in the picture by attaching it to its leader.
  However, contrary to the algebraic part, a differential equation does not only affect its leader, but also the derivatives of its leader.
  This is because property \ref{rankingproperty2} of a differential ranking implies $\partial\ld(p)=\ld(\partial p)\ \forall\partial\in\Delta, p\in F\{U\}$.
  For example $\partial_t (\underline{u\der{0,1}}+u\der{0,0}u\der{1,0}) =$ $\underline{u\der{0,2}}+u\der{0,1}u\der{1,0}+u\der{0,0}u\der{1,1}$.
  In the diagram we illustrate this by drawing a cone with apex $u\der{0,1}$.

  \parpic[r]{\textnormal{
    \begin{tikzpicture}[domain=0:4, smooth, x=\picturesize, y=\picturesize, place/.style={circle,draw=black,fill=black,inner sep=1pt}]
      \draw[color=black, line width=1pt, ->] (0.0,0.0) -- (0.0,4.0) node[right]{$\partial_{t}$};
      \draw[color=black, line width=1pt, ->] (0.0,0.0) -- (4.0,0.0) node[above]{$\partial_{x}$};
      \node[place,label=180:$u\der{0,0}$] (u00) at (0.0,0.0) {};
      \node[place,label=-90:$u\der{1,0}$] (u10) at (1.0,0.0) {};
      \node[place,label=180:$\underline{u\der{0,1}}$] (u01) at (0.0,1.0) {};
      \node[place,label=-90:$\underline{u\der{2,0}}$] (u20) at (2.0,0.0) {};
      \node[place,label=90:] (u11) at (1.0,1.0) {};\node[above=-0.15 of u11] {$\underline{u\der{1,1}}$};
      \node[place,label=180:$\underline{u\der{0,2}}$] (u02) at (0.0,2.0) {};
      \node[place,label=-90:$\underline{u\der{3,0}}$] (u30) at (3.0,0.0) {};
      \node[place,label=90:] (u21) at (2.0,1.0) {};\node[above=-0.15 of u21] {$\underline{u\der{2,1}}$};
      \node[place,label=90:] (u12) at (1.0,2.0) {};\node[above=-0.15 of u12] {$\underline{u\der{1,2}}$};
      \node[place,label=180:$\underline{u\der{0,3}}$] (u03) at (0.0,3.0) {};
      \node (eq1) at (-3.0,4.0) {$(\underline{u\der{0,1}}+u\der{0,0}u\der{1,0})_=$};
      \draw [thick,->] (eq1) to [bend left=-15] (u01);
      \node (eq2) at (-3.0,0.0) {$(u\der{2,0})_=$};
      \draw [thick,->] (eq2) to [bend left=15] (u20);
      \draw[color=gray, line width=1pt] (u01)+(-0.1,-0.1) rectangle +(0.1,0.1);
      \draw[color=gray, line width=1pt,->] (u01)+(-0.1,-0.1) -- +(3.7,-0.1);
      \draw[color=gray, line width=1pt,->] (u01)+(-0.1,-0.1) -- +(-0.1,2.7);
      \draw[color=gray, line width=1pt] (u20)+(-0.1,-0.1) rectangle +(0.1,0.1);
      \draw[color=gray, line width=1pt,->] (u20)+(-0.1,-0.1) -- +(1.7,-0.1);
      \draw[color=gray, line width=1pt,->] (u20)+(-0.1,0.1) -- +(1.7,0.1);
    \end{tikzpicture}}
  }
  
  Assume that we are only interested in solutions of the inviscid \textsc{Burger}'s equation which are linear in $x$.
  So, we add the second equation $(\underline{u\der{2,0}})_=$ to our system.
  This second equation also affects the derivatives of its leader.
  In particular, $(\underline{u\der{0,1}}+u\der{0,0}u\der{1,0})_=$ and $(\underline{u\der{2,0}})_=$ both affect the differential variable $u\der{2,1}$ and its derivatives.
  This contradicts the triangularity of the system.
  According to the involutive approach as suggested by \textsc{Janet}, we don't allow certain equations to be derived by certain partial derivations.
  In this example, we allow $(\underline{u\der{2,0}})_=$ to be derived only by $\partial_x$.
  In the diagram we illustrate this by drawing a (degenerate) cone with apex $u\der{2,0}$ in direction of $\partial_x$.
  Thus, the differential consequence $(\partial_t\underline{u\der{2,0}})_=$ is not yet considered and, so, we have to add it as a separate equation for further treatment.

\end{myexample}

A set $W$ of differential variables is \textbf{closed} under the action of $\Delta'\subseteq\Delta$ if $\partial_i w\in W$ for all $\partial_i\in\Delta'$ and $w\in W$.
The smallest set containing a differential variable $w$, which is closed under $\Delta'$, is called a \textbf{cone} and denoted by $\langle w\rangle_{\Delta'}$.
In this case, we call the elements of $\Delta'$ \textbf{reductive derivations}\footnote{ In~\cite{Ger3} and \cite[Chap.~7]{Seiler} the reductive derivations are called multiplicative variables and in \cite{thomasalg_casc} they are called admissible derivations.}.
The $\Delta^{\prime}$-closed set generated by a set $W$ of differential variables is defined as
\[
  \langle W\rangle_{\Delta^{\prime}}:=\bigcap_{\stackrel{W_i\supseteq W}{W_i\, \Delta^{\prime}\mbox{\scriptsize -closed}}} W_i\ \ \subseteq \ \ \langle U\rangle_{\Delta}\enspace.
\]

For a finite set $W=\{w_1,\ldots,w_r\}$, the {\sc Janet} \textbf{division} algorithmically assigns reductive derivations to the elements of $W$ such that the cones generated by the $w\in W$ are disjoint (cf.\ \cite{GYB} for a fast algorithm).
We call these derivations {\sc Janet}-reductive.
The derivation $\partial_l\in\Delta$ is assigned to the cone generated by $w=u^{(j)}_{\mathbf{i}}\in W$ as reductive derivation, if and only if
\[
  \mathbf{i}_l=\max\left\{\mathbf{i}'_l\mid u^{(j)}_{\mathbf{i}'}\in W,\mathbf{i}'_k=\mathbf{i}_k\mbox{ for all }1\le k<l\right\}
\]
holds \cite[Ex.~3.1]{GerI}.
We remark that $j$ is fixed in this definition, i.e., when constructing cones we only take into account other differential variables belonging to the same differential indeterminate.
Furthermore, the assignment of reductive derivations to $w\in W$ in general depends on the whole set $W$.
The reductive derivations assigned to $w$ are denoted by $\Delta_W(w)\subseteq\Delta$ and we call the cone $\langle w\rangle_{\Delta_W(w)}$ the {\sc Janet} \textbf{cone} of $w$ with respect to $W$.
This construction ensures disjointness of cones but not necessarily that the union of cones equals $\langle W\rangle_\Delta$.
The problem is circumvented by enriching $W$ to its {\sc Janet} \textbf{completion} $\widetilde W\supseteq W$.
This completion $\widetilde W$ is successively created by adding any 
\[
  \tilde w=\partial_i w_j\not\in\biguplus_{w\in\widetilde W}\langle w\rangle_{\Delta_{\widetilde W}(w)}  
\]
to $\widetilde W$, where $w_j\in\widetilde W$ and $\partial_i\in\Delta\setminus\Delta_{\widetilde W}(w_j)$. This leads to the disjoint {\sc Janet} \textbf{decomposition}
\[
  \langle W\rangle_\Delta=\biguplus_{w\in\widetilde W}\langle w\rangle_{\Delta_{\widetilde W}(w)}
\]
that algorithmically separates a $\Delta$-closed set $\langle W\rangle_\Delta$ into finitely many cones $\langle w\rangle_{\Delta_{\widetilde W}(w)}$.
For details see \cite[Def.\ 3.4]{GerI} and \cite[Cor.\ 4.11]{GB1}.

We extend the {\sc Janet} decomposition from differential variables to differential polynomials according to their leaders.
To be precise, $\Delta_T(q):=\Delta_{\ld(T)}(\ld(q))$ for finite $T\subset F\{U\}$ and $q\in T$.
We call a derivative of an equation by a finite (possibly empty) sequence of derivations a \textbf{prolongation}.
If all these derivations are reductive, the derivative is called \textbf{reductive prolongation} of $q$ with respect to $T$.
Otherwise it is called \textbf{non-reductive prolongation}.
  
A differential polynomial $p\in F\{U\}$ is called \textbf{reducible} modulo $q\in F\{U\}$, if there exists $\mathbf{i}\in\Z_{\ge 0}^n$ such that $\partial_1^{\mathbf{i}_1}\cdot\ldots\cdot\partial_n^{\mathbf{i}_n}\ld(q)=\ld(\partial_1^{\mathbf{i}_1}\cdot\ldots\cdot\partial_n^{\mathbf{i}_n}q)=\ld(p)$ and $\mdeg(\partial_1^{\mathbf{i}_1}\cdot\ldots\cdot\partial_n^{\mathbf{i}_n}q)\le\mdeg(p)$.
For $\mathbf{i}\neq(0,\ldots,0)$ the condition on the main degree always holds.
We now restrict ourselves to reductive prolongations:
For a finite set $T\subset F\{U\}$, we call a differential polynomial $p\in F\{U\}$ {\sc Janet}-\textbf{reducible} modulo $q\in T$ w.r.t. $T$, if $p$ is reducible modulo $q$ and $\partial_1^{\mathbf{i}_1}\cdot\ldots\cdot\partial_n^{\mathbf{i}_n}q$ is a reductive prolongation of $q$ w.r.t. $T$, with $\mathbf{i} \in \mathbb{Z}_{\ge 0}^n$ from the reducibility conditions.
We also say that $p$ is {\sc Janet}-\textbf{reducible} modulo $T$ if there is a $q\in T$ such that $p$ is {\sc Janet}-reducible modulo $q$ w.r.t. $T$.

A set of differential variables $T\subset\langle U\rangle_{\Delta}$ is called \textbf{minimal}, if for any set $S\subset\langle U\rangle_{\Delta}$ with \[\biguplus_{t\in T}\langle t\rangle_{\Delta_T(t)}=\biguplus_{s\in S}\langle s\rangle_{\Delta_S(s)}\] the condition $T\subseteq S$ holds \cite[Def.~4.2]{GB2}.
We call a set of differential polynomials minimal, if the corresponding set of leaders is minimal.

At each step of the algorithm we assign reductive derivations to the equations in $(S_T)^=$.
When an equation $p$ is not reducible modulo $(S_T)^=$, it is added to $(S_T)^=$.
Then, we remove all polynomials from $S_T$ that have a leader which is derivative of $\ld(p)$.
This will later ensure minimality.
In addition, when adding a new equation to $(S_T)^=$, all non-reductive prolongations are put into the queue.
This is formalized in the following algorithm.

\begin{myalgorithm}[{\sf InsertEquation}]\ \\
\label{algo_differential_insertequation}
 \textit{Input:} A system $S'$ and a polynomial $p_=\in F\{U\}$ not reducible modulo $(S'_T)^=$.\\
 \textit{Output:} A system $S$, where $(S_T)^=\subseteq (S_T')^=\cup\{p_=\}$ is maximal satisfying
\[
\left(\ld(S_T)\setminus\{\ld(p)\}\right)\cap\langle\ld(p)\rangle_{\Delta}=\emptyset,
\]
\[
S_Q=S_Q'\cup(S'_T\setminus S_T)\cup\{(\partial_i q)_=\mid q\in (S_T)^=,\partial_i\not\in\Delta_{((S_T)^=)}(q)\}\enspace.
\]
 \textit{Algorithm:}
  \begin{algorithmic}[1]
    \STATE $S\gets S'$
    \STATE $S_T \gets S_T\cup\{p_=\}$
    \FOR{$q\in S_T\setminus\{p\}$}
      \IF{$\ld(q)\in\langle\ld(p)\rangle_{\Delta}$}
        \STATE $S_Q\gets S_Q\cup\{q\}$
        \STATE $S_T\gets S_T\setminus\{q\}$
     \ENDIF 
    \ENDFOR
    \STATE Reassign reductive derivations to $(S_T)^=$
    \STATE $S_Q\gets S_Q\cup \{(\partial_i q)_=\mid q\in (S_T)^=,\partial_i\notin\Delta_{((S_T)^=)}(q)\}$
    \RETURN $S$
  \end{algorithmic}
\end{myalgorithm}

Correctness and termination are obvious.
We remark that a non-reductive prolongation might be added to $S_Q$ several times.
An implementation should remember which prolongations have been added before to avoid redundant computations.

\subsection{Differential Simple Systems}\label{differential reduction}

In this subsection, we extend the algebraic reduction algorithm to its differential counterpart.
Finally, we can define differential simple systems at the end of this subsection.

The {\sc Janet} partition of the dependent differential variables into cones provides a mechanism to find the unique reductor for the differential reduction in a fast way (cf.\ \cite{GYB}).
We prolong this reductor and afterwards apply a pseudo reduction algorithm.

For a valid pseudo-reduction, we need to ensure that initials (and initials of the prolongations) of equations are non-zero.
Let $r\in F\{U\}$ with $x=\ld(r)$ and define the \textbf{separant} $\sepa(r):=\frac{\partial r}{\partial x}$.
One easily checks that the initial of any non-trivial prolongation of $r$ is $\sepa(r)$ and the separant of any square-free equation $r$ is non-zero (cf.\ \cite[\textsection I.8, Lemma 5]{Kol} or \cite[\textsection3.1]{Hubert2}).
So, by making sure that the equations have non-vanishing initials and are square-free, as in the algebraic case, we ensure that we can reduce modulo all prolongations of $r$.
This provides the correctness of the following reduction algorithm. \footnote{
In differential algebra, one usually distinguishes a (full) differential reduction as used here and a partial (differential) reduction.
Partial reduction only employs \emph{proper} derivations of equations for reduction (cf.\ \cite[\textsection I.9]{Kol} or \cite[\textsection 3.2]{Hubert2}).
This is useful for separation of differential and algebraic parts of the algorithm and for the use of {\sc Rosenfeld}'s Lemma (cf.\ \cite{Rosenfeld}), which is the theoretical basis for the \textsc{Rosenfeld}-\textsc{Gr\"obner} algorithm (cf.\ \cite{BLOP,BLOP2,Hubert2}.)
}

\begin{myalgorithm}[{\sf Reduce}]\label{algo_differential_reduce}\ \\
 \textit{Input:} A differential system $S$ and a polynomial $p \in F\{U\}$.\\
 \textit{Output:} 
 A polynomial $q$ that is not \textsc{Janet}-reducible modulo $S_T$
 with $\phi_e(p)=0$ if and only if $\phi_e(q)=0$ for each $e\in\sol(S)$.\\
 \textit{Algorithm:} \begin{algorithmic}[1]
                      \STATE $x \gets \ld(p)$
                       \WHILE{exists $q_=\in(S_T)^=$ and $\mathbf{i}\in\Z_{\ge0}^n$ with $\mathbf{i}_j=0$ for $\partial_j\not\in\Delta_{(S_T)^=}(q)$ such that $\partial_1^{\mathbf{i}_1}\cdot\ldots\cdot\partial_n^{\mathbf{i}_n}\ld(q)=\ld(p)$ and  $\mdeg(\partial_1^{\mathbf{i}_1}\cdot\ldots\cdot\partial_n^{\mathbf{i}_n}p) \geq \mdeg(q)$ hold}
                       \STATE $p \gets \prem(p, \partial_1^{\mathbf{i}_1}\cdot\ldots\cdot\partial_n^{\mathbf{i}_n}q, x)$
                       \STATE $x \gets \ld(p)$
                      \ENDWHILE
                      \IF{$\textrm{\sf Reduce}(S, \ini(p)) = 0$}
                       \RETURN $\textrm{\sf Reduce}(S, p-\ini(p)x^{\mdeg(p)})$
                      \ELSE
                       \RETURN $p$
                      \ENDIF
                     \end{algorithmic}
\end{myalgorithm}

A polynomial $p\in F\{U\}$ \textbf{reduces to $q$ modulo $S_T$} if $\mathsf{Reduce}(S, p)=q$.
A polynomial $p\in F\{U\}$ is called \textbf{reduced\footnote{
There is a fine difference between not being reducible and being reduced.
In the case of not being reducible the initial of a polynomial can still reduce to zero and iteratively the entire polynomial.
}
modulo $S_T$} if it reduces to itself.
The properties of the algebraic reduction algorithm from Remark (\ref{rem_reduceprop}) also apply for this reduction algorithm.

Termination of the reduction algorithm is provided by \textsc{Dickson}'s Lemma (cf.\ \cite[Chap.~2, Thm.~5]{CLO} or \cite[\textsection 0.17, Lemma 15]{Kol}), which states that the ranking $<$ is well-founded on the set of leaders, i.e., a strictly $<$-descending chain of leaders is finite.
  
\begin{myexample}\label{example_janet_reduce}
  We continue Example (\ref{example_janet_combinatoric}) and take care of the differential consequence $(\underline{u\der{2,1}})_=$.
  \parpic[r]{\textnormal{
    \begin{tikzpicture}[domain=0:4, smooth, x=1.0cm, y=1.0cm, place/.style={circle,fill=black,inner sep=1pt}]
      \draw[color=black, line width=1pt, ->] (0.0,0.0) -- (0.0,2.4) node[right]{\tiny$\partial_{t}$};
      \draw[color=black, line width=1pt, ->] (0.0,0.0) -- (4.0,0.0) node[right]{\tiny$\partial_{x}$};
      \node[place] (u00) at (0.0,0.0) {};
      \node[] (u00p) at (0.0,-0.3) {\tiny{$0$}};
      \node[place] (u10) at (1.0,0.0) {};
      \node[place] (u01) at (0.0,1.0) {};
      \node[place] (u20) at (2.0,0.0) {};
      \node[place] (u11) at (1.0,1.0) {};
      \node[place] (u02) at (0.0,2.0) {};
      \draw[color=gray, line width=1pt] (u01)+(-0.1,-0.1) rectangle +(0.1,0.1);
      \draw[color=gray, line width=1pt,->] (u01)+(-0.1,-0.1) -- +(3.7,-0.1);
      \draw[color=gray, line width=1pt,->] (u01)+(-0.1,-0.1) -- +(-0.1,1.2);
      \draw[color=gray, line width=1pt] (u20)+(-0.1,-0.1) rectangle +(0.1,0.1);
      \draw[color=gray, line width=1pt,->] (u20)+(-0.1,-0.1) -- +(1.7,-0.1);
      \draw[color=gray, line width=1pt,->] (u20)+(-0.1,0.1) -- +(1.7,0.1);
      \node[] (u20p) at (1.9,-0.4) {\tiny{$3u\der{1,0}\underline{u\der{2,0}}$}}
        edge[->,bend left=50,color=black, line width=1pt] node[below,color=black]{\tiny{\begin{tabular}{c}reduce\\modulo $p_2$\end{tabular}}} (u00p);
      \node[place] (u30) at (3.0,0.0) {};
      \node[] (u30p) at (3.3,-0.6) {\tiny{\begin{tabular}{c}
                                    $u\der{0,0}\underline{u\der{3,0}}$\\
                                    $+3u\der{1,0}u\der{2,0}$
                                    \end{tabular}}
                                  }
        edge[->,bend left=40,color=black, line width=1pt] node[below,color=black]{\tiny{\begin{tabular}{c}reduce\\modulo $\partial_x p_2$\end{tabular}}} (u20p);
      \node[place] (u21) at (2.0,1.0) {};
      \node[] (u21p) at (2.4,1.2) {\tiny{$\underline{u\der{2,1}}$}}
        edge[->,bend left=30,color=black, line width=1pt] node[right,color=black]{\tiny{\begin{tabular}{c}reduce\\modulo $\partial_x^2 p_1$\end{tabular}}} (u30p);
      \node[place] (u12) at (1.0,2.0) {};
      \node[] (p1) at (0.3,1.3) {\small{$p_1$}};
      \node[] (p2) at (2.0,0.3) {\small{$p_2$}};
      \fill[color=black] (u21) circle (2pt);
      \fill[color=black] (u30) circle (2pt);
      \fill[color=black] (u20) circle (2pt);
      \fill[color=black] (u00) circle (2pt);
    \end{tikzpicture}}
  } 
  We reduce $(\underline{u\der{2,1}})_=$ modulo the system $S$ with
  \[S_T:=\left\{p_1:=(\underline{u\der{0,1}}+u\der{0,0}u\der{1,0})_=, p_2:=(\underline{u\der{2,0}})_=\right\}\mbox{ .}\]
  First, we observe, that $\ld(\underline{u\der{2,1}})=u\der{2,1}$ is in the cone generated by $\ld(p_1)$ and $\ld(\partial_x^2 p_1)=\ld(\underline{u\der{2,1}})$.
  Thus, we reduce $(\underline{u\der{2,1}})$ modulo $\partial_x^2 p_1$ and the pseudo reduction yields $u\der{0,0}\underline{u\der{3,0}}+3u\der{1,0}u\der{2,0}$.
  Second, we reduce $u\der{0,0}\underline{u\der{3,0}}+3u\der{1,0}u\der{2,0}$ modulo $\partial_xp_2$, because $u\der{3,0}$ lies in the cone generated by $(\underline{u\der{2,0}})_=$.
  This results in $3u\der{1,0}\underline{u\der{2,0}}$ and a third reduction step modulo $p_2$ produces zero.
  As a result, the only differential consequence is already implied by the system.
  In this desirable situation, there are no further integrability conditions, which motivates the definition of involutivity below.
\end{myexample}

Now, we define differential simple systems.
We demand algebraic simplicity, involutivity of differential equations as seen in the previous Example (\ref{example_janet_reduce}), and some minimality conditions.

\begin{mydefinition}[Differential Simple Systems]
A differential system $S$ is (\textsc{Janet}-) \textbf{involutive}, if
all non-reductive prolongations of $(S_T)^=$ reduce to zero modulo $(S_T)^=$.\\
A system $S$ is called \textbf{differentially simple} or \textbf{simple}, if
\begin{enumerate}
  \item $S$ is algebraically simple (in the finitely many differential variables that appear in it),
  \item $S$ is involutive,
  \item $S^=$ is minimal,
  \item no inequation in $S^{\neq}$ is reducible modulo $S^=$.
\end{enumerate}
A disjoint decomposition of a system into differentially simple subsystems is called \textbf{(differential) }{\sc Thomas}\textbf{ decomposition}.
\end{mydefinition}

As in the algebraic case, every simple system has a solution in $E$.

\subsection{Differential Decomposition Algorithm}\label{differential algorithm}

The differential \textsc{Thomas} decomposition algorithm is a modification of the algebraic \textsc{Thomas} decomposition algorithm.
We have already introduced the new algorithms \textsf{InsertEquation} (\ref{algo_differential_insertequation}) for adding new equations into the systems and \textsf{Reduce} (\ref{algo_differential_reduce}) for reduction, that can replace their counterparts in the algebraic algorithm.
This subsection provides the necessary correctness and termination proofs for this modified algorithm.
It then demonstrates this algorithm with examples.

\begin{myalgorithm}[{\sf DifferentialDecompose}]\label{algo_differential_decompose}\ \\
 \textit{Input:} A differential system $S^\prime$ with $({S^\prime})_T=\emptyset$. \\
 \textit{Output:} A differential {\sc Thomas} decomposition of $S^\prime$. \\
 \textit{Algorithm:} The algorithm is obtained by replacing the two subalgorithms {\sf\nobreak InsertEquation} and {\sf\nobreak Reduce} in Algorithm (\ref{algo_decompose}) with their differential counterparts (\ref{algo_differential_insertequation}) and (\ref{algo_differential_reduce}), respectively.
\end{myalgorithm}

\begin{proof}[Correctness]

The correctness proof of the algebraic decomposition Algorithm (\ref{algo_decompose}) also holds verbatim for the differential case.
Therefore, we do not need to show that the output is algebraically simple.
We will prove three loop invariants for any system $S \in P \cup Result$:
  
\begin{enumerate}
  \item $(S_T)^=$ is minimal.
  \item No inequation in $(S_T)^{\not=}$ is \textsc{Janet}-reducible modulo $S_T$.
  \item Let $r$ be any non-reductive prolongation of $(S_T)^=$.
        Then $r$ reduces to zero by using both conventional differential reductions\footnote{i.e. modulo any prolongation} of $(S_Q)^=$ \emph{and} reductions modulo reductive prolongations of $(S_T)^=$.
\end{enumerate}

The first loop invariant is a purely combinatorial matter, which is proved by \cite{Ger2} for an algorithm using exactly the same combinatorial approach.

Proving the second loop invariant is equally simple.
On the one hand, a newly added inequation $q$ in $S_T$ is not \textsc{Janet}-reducible modulo $(S_T)^=$, since algorithm \textsf{Reduce} (\ref{algo_differential_reduce}) is applied to it before insertion.
On the other hand, algorithm {\sf InsertEquation} (\ref{algo_differential_insertequation}) removes all inequations from $S_T$ which are divisible by a newly added equation and places them into $S_Q$.

The third loop invariant clearly holds at the beginning of the algorithm, because $S_T$ is empty.

We claim that reduction of an equation $q_=\in S_Q$ by $(S_T)^=$ in line \ref{reduce_q} of Algorithm (\ref{algo_decompose}) does not affect the loop invariant, i.e. any non-reductive prolongation $r$ reducing to zero beforehand reduces to zero afterwards.
We prove this claim by performing a single reduction step on $q$, which generalizes by an easy induction.
Let $q':=\prem(q,p,x)=m\cdot q-\pquo(q,p,x)\cdot p$ be a pseudo remainder identity
(see (\ref{prempquo}) on page \pageref{prempquo})
reducing $q$ to $q'$ modulo $p$.
Then a pseudo remainder identity $\prem(r,q,x)=m'\cdot r-\pquo(r,q,x)\cdot q$ describing a reduction of $r$ modulo $q$ might simply be rewritten as the iterated identity
\[
  \underbrace{m\cdot \prem(r,q,x)}_{\prem(\prem(r,p,x),q',x)}=m\cdot 
  m'\cdot r-\pquo(r,q,x)\cdot q'-\pquo(r,q,x)\cdot \pquo(q,p,x)\cdot p.
\]
Using the \textsc{Leibniz} rule the same holds for reduction modulo partial derivatives of $q$.
This holds especially for an equation $q_=\in S_Q$ reducing to $0$ modulo $(S_T)^=$ in line \ref{reduce_q}, which can be removed from $S_Q$ without violating the loop invariant.

Now, we consider line \ref{q_2}, where \textsf{InsertEquation} inserts the square-free part $p_=$ of $q_=$ into $S_T$ and show that this does not violate the third loop invariant. 
First, the non-reductive prolongations in $\{(\partial_i r)_=\mid r\in (S_T)^= \partial_i\notin\Delta_{((S_T)^=)}(r)\}$ are added to $S_Q$ as equations.
Thus, any of these reduce to $0$ modulo $(S_Q)^=$.
Second, moving equations from $S_T$ back into $S_Q$ in \textsf{InsertEquation} does not change the loop invariant either, because their reductive prolongations can still be used for reduction afterwards.
Third, every non-reductive prolongation that reduced to zero using $q_=\in (S_Q)^=$ still reduces to zero after \textsf{InsertEquation}.
This holds for two reasons.
On the one hand, everything that reduces to zero modulo $q_=$, also reduces to zero modulo $p_=$. 
Write $m\cdot q=p\cdot q_1$ with $\ld(m)<x$ and $\pma(m)\not=0\ \forall\mathbf{a}\in\sol(S_{<\ld(q)})$.
Then $p$ algebraically pseudo-reduces $q$ to zero.
Any derivative $\partial q$ of $q$ is reduced to zero modulo $p_=$ and $(\partial p)_=$, since $\partial(m\cdot q)=(\partial p)\cdot q_1+p\cdot(\partial q_1)$ for any $\partial\in\Delta$.
Inductively, the same holds for repeated derivatives of $q_=$.
Therefore, $p_=$ implies all constraints given by $q_=$.
On the other hand, all reduction steps modulo $p_=$ are either \textsf{Janet}-reductions modulo $p_=$ w.r.t. $S_T$ or differential reductions modulo non-reductive prolongations of $p_=$.
The latter equations have been added to $S_Q$.

When computing the gcd of two equations in line \ref{resplitgcd}, the gcd of $q$ and $(S_T)_x$ will be inserted into $S_T$ and reduces everything to zero that both $q$ and $(S_T)_x$ did.
As above, the non-reductive prolongations are covered by inserting them into $S_Q$ and the reductive prolongations are implied.

Dividing an equation $(S_T)_x$ by an inequation $q_{\not=}$ in lines \ref{dividebyinequation} and \ref{q_3} also influences $(S_T)^=$.
The new equation $p_=$, being a divisor of $(S_T)_x$, reduces everything to zero that $(S_T)_x$ and its non-reductive prolongations did by the same arguments as before.

This proves the third loop invariant.
When the algorithm terminates, $S_Q$ is empty and thus all non-reductive prolongations from $(S_T)^=$ \textsc{Janet}-reduce to zero modulo $(S_T)^=$.
The system is therefore involutive.

Furthermore, the first loop invariant implies minimality and the second loop invariant implies that no inequation is reducible by an equation, since for an involutive set reducibility is equivalent to \textsc{Janet}-reducibility. \qed
\end{proof}

Our main tool for proving the termination of the algorithm is using six orders on differential systems.
These are similar to the four orders used to show the termination of the algebraic decomposition algorithm.
We use {\sc Dickson}'s lemma as main tool to show the well-foundedness of these orders.

\begin{defrem}
Define the orders $\prec_{1a}$, $\prec_{1b}$, $\prec_{1c}$, $\prec_{2}$, $\prec_{3}$, and $\prec_{4}$ as follows.
\begin{itemize}
  \item[$\prec_{1a}$:] For $V\subseteq\langle U\rangle_{\Delta}$ there is a unique minimal set $\nu(V)\subseteq V$ with $V\subseteq \langle \nu(V)\rangle_{\Delta}$ \cite[Chap.~2, \textsection4, exercise 7 and 8]{CLO}, called \textbf{canonical differential generators} of $V$.
  For a system $S$, define $\nu(S)$ as $\nu(\ld((S_T)^=))$. 
  For systems $S, S'$ we define $S\prec_{1a}S'$ if and only if 
 $\min_<(\nu(S)\setminus\nu(S'))<\min_<(\nu(S')\setminus\nu(S))$. An empty set is assumed to have $x_\infty$ as minimum, which is $<$-larger than all differential variables.
  By {\sc Dickson}'s lemma, $\prec_{1a}$ is well-founded.

 \item[$\prec_{1b}$:] 
 For systems $S, S'$ define $S\prec_{1b}S'$ if and only if $S\osim_{1a} S'$ and $\min_<\left(\ld((S_T)^=)\setminus\ld(({S'}_T)^=)\right)<$\\$\min_<\left(\ld(({S'}_T)^=)\setminus\ld((S_T)^=)\right)$.
 Minimality of $(S_T)^=$ at each step of the algorithm and the constructivity property of the {\sc Janet} division \cite[Prop.~4.13]{GB1} imply well-foundedness of $\prec_{1b}$ \cite[Thm.~4.14]{GB1}.

 \item[$\prec_{1c}$:] 
 For systems $S$ and $S'$ with $S\osim_{1a}S'$ and $S\osim_{1b}S'$, both $(S_T)^=$ and $({S'}_T)^=$ have the same leaders $x_1,\ldots,x_l$.
 Define $S\prec_{1c,x_k} S'$ if and only if $\mdeg((S_T)^=_{x_i})<\mdeg((S'_T)^=_{x_i})$.
 This order is clearly well-founded.
 For these systems define $S\prec_{1c}S'$ as $[\prec_{1c,x_1},\ldots,\prec_{1c,x_l}]$, which is again well-founded as a composite order.

 \item[$\prec_2$:] This is defined identical to the algebraic $\prec_2$.
 We remark, that in this case the set of possible leaders is $\{1\}\cup\langle U\rangle_{\Delta}$.
 To show well-foundedness of the differential ordering $\prec_2$ we use that $<$ is well-founded on the set of leaders as implied by {\sc Dickson}'s lemma.
 This way, $<$ is extended to a well-founded ordering on $\{1,x_{\infty}\}\cup\langle U\rangle_{\Delta}$ with $1<y$ and $y<x_{\infty}$ for all $y\in\langle U\rangle_{\Delta}$.

 \item[$\prec_3$:] This is verbatim the same condition and proof of well-foundedness as in the algebraic case.
 However, in the latter proof, we do a {\sc Noether}ian induction \cite[III.6.5, Prop.~7]{Bou68} instead of an ordinary induction.

 \item[$\prec_4$:] This is identical to the algebraic case.
\end{itemize}
Remark (\ref{composite_well_founded}) provides the well-foundedness of the composite order $\prec:=[{\prec_{1a},} {\prec_{1b},} {\prec_{1c},} {\prec_{2},} {\prec_{3},} {\prec_{4}}]$.
\end{defrem}

\begin{proof}[Termination]

We prove termination the same way as in the algebraic case.
All arguments where systems get $\prec_2$, $\prec_3$, or $\prec_4$ smaller apply verbatim here.

In the algebraic case a system $\prec_1$-decreases if and only if either an equation is added to $S_T$ or the degree of an existing equation in $S_T$ is decreased.
We adapt this argument to the differential case:
On the one hand, inserting a new equation with a leader that is not yet present in $\ld((S_T)^=)$ decreases either $\prec_{1a}$ or $\prec_{1b}$.
On the other hand, if an existing equation in $(S_T)^=$ is replaced by one with the same leader and lower degree, the system $\prec_{1c}$-decreases.

Thus, like in the algebraic termination proof, we have a strictly decreasing chain of systems and, thus, termination is proved. \qed
\end{proof}

In the following examples, we use jet notation for differential polynomials, e.g., $u_{x,x,y}:=u_{2,1}$ in the case $\Delta=\{\partial_x,\partial_y\}$ and $U=\{u\}$.

We give an example taken from \cite[pp.\ 597-600]{BC}:

\begin{myexample}[Cole-Hopf Transformation]
For $F:=\R(x,t)$, $\Delta=\{\frac{\partial}{\partial x},\frac{\partial}{\partial t}\}$, and $U=\{\eta,\zeta\}$ consider the heat equation $h=(\eta_t+\eta_{xx})_=$ and {\sc Burger}'s equation $b=(\zeta_t+\zeta_{xx}+2\zeta_x\cdot\zeta)_=$.

First, we claim that any power series solution for the heat equation with a non-zero constant term can be transformed to a solution of \textsc{Burger}'s equation using the {\sc Cole-Hopf} transformation $\lambda: \eta\mapsto\frac{\eta_x}{\eta}$.
A differential {\sc Thomas} decomposition for an orderly ranking with $\zeta_x>\eta_t$ of
\[
  \{h_=,\underbrace{(\eta\cdot\zeta-\eta_x)_=}_{\Leftrightarrow \zeta=\lambda(\eta)},\eta_{\not=}\}
\]
consists of the single system
\[
  S=\{(\underline{\eta_x}-\eta\cdot\zeta)_=,(\eta\cdot\underline{\zeta_x}+\eta_t+\eta\cdot\zeta^2)_=,\underline{\eta}_{\not=}\}
\]
and one checks that $\mathsf{Reduce}(S,b)=0$ holds.
This implies that $\lambda$ maps any non-zero solution of the heat equation to a solution of {\sc Burger}'s equation.

In addition we claim that $\lambda$ is surjective.
For the proof we choose an elimination ranking (cf. \cite[\textsection 8.1]{Hubert2} or \cite{BoulierElimination}) with $\eta\gg\zeta$, i.e., $\eta_\mathbf{i}>\zeta_\mathbf{j}$ for all $\mathbf{i},\mathbf{j}\in \mathbb{Z}_{\ge 0}$.
We compute a differential {\sc Thomas} decomposition of $\{h_=,b_=,(\eta\cdot\zeta-\eta_x)_=,\eta_{\not=}\}$.
It consists of the single system
\[
  S=\{(\underline{\eta_x}-\eta\cdot\zeta)_=,(\eta\cdot\zeta_x+\underline{\eta_t}+\eta\cdot\zeta^2)_=,b_=,\underline{\zeta}_{\not=}\}\enspace.
\]
The elimination ordering guarantees that the only constraint for $\zeta$ is {\sc Burger}'s equation $b_=$.
As $S$ is simple, for any solution $f\in\sol(b_=)$ there exists a solution $(g,f)\in\sol(S)$ (cf.\ (\ref{exist_sol})), implying that $\lambda$ is surjective.
\end{myexample}

Elements of the $\Delta$-field $F$ are not subjected to splittings and assumed to be non-zero.
However, we are able to model the elements of $F$ as differential indeterminates.
For example for $F=\C(x)$ with $\Delta=\{\frac{\partial}{\partial x}\}$, we can study a differential polynomial ring over $\C\{X\}$ instead and replace $x$ by $X$ in all equations and inequations.
We subject $X$ to the relation $\frac{\partial}{\partial x}X=1$ for $X$ being ``generic'' or  $(\frac{\partial}{\partial x}X-1)\cdot\frac{\partial}{\partial x}X=0$, if we allow specialization of $X$.
Both these cases are considered in examples (\ref{ex_model_field_elements1}) and (\ref{ex_model_field_elements2}), respectively, and will be subject of further study.

\begin{myexample}\label{ex_model_field_elements1}
For $F:=\C(x)$, $\Delta=\{\frac{\partial}{\partial x},\frac{\partial}{\partial t}\}$ and $U=\{u\}$ consider the special case
\begin{equation}\label{autonomousequation}
\left(u_t-u_{xx}-x\cdot u_x-u\right)_=
\end{equation}
of the \textsc{Fokker-Planck} equation.
We add an auxiliary differential indeterminate $X$ to $U$ and instead examine the equation
\begin{equation}\label{autonomousequation2}
\left(u_t-u_{xx}-X\cdot u_x-u\right)_=, \left(X_x-1\right)_=, \left(X_t\right)_=
\end{equation}
in the $\Delta$-ring $\C\{X,u\}$. An elimination ranking $X\gg u$ splits the system (\ref{autonomousequation2}) into two simple systems:
\begin{enumerate}[(i)]
\item\label{systemone} $\displaystyle
\left(u_x\cdot(-\underline{u_{xxx}}+u_{xt}-2u_x)-u_{xx}\cdot(u_t-u_{xx}-u)\right)_=,\\
\left(u_x\cdot(-\underline{u_{xxt}}+u_{tt}-u_t)-u_{xt}\cdot(u_t-u_{xx}-u)\right)_=,\\
\left(u_t-u_{xx}-\underline{X}\cdot u_x-u\right)_=,
\left(\underline{u_x}\right)_{\not=}$
\item $\displaystyle \left(\underline{u_x}\right)_=, \left(\underline{u_t}-u\right)_= , \left(\underline{X_x}-1\right)_=, \left(\underline{X_t}\right)_=$
\end{enumerate}
Due to the ranking, the first two equations in (\ref{systemone}) generate $\left(F\{u\}[\Delta]\cdot(u_t-u_{xx}-x\cdot u_x-u)\right)\cap\C\{u\}$, i.e., they have constant coefficients.
These two equations are the derivatives of $\frac{u_t-u_{xx}-u}{u_x}-x$, which is clearly equivalent to (\ref{autonomousequation}) in the case $u_x\neq0$.
\end{myexample}

The next example sketches an approach to treat equations with variable coefficients and find submanifolds where solutions behave differently.

\begin{myexample}\label{ex_model_field_elements2}
For $F:=\C(x,y)$, $\Delta=\{\frac{\partial}{\partial x},\frac{\partial}{\partial y}\}$ and $U=\{u\}$ consider
\[
(xy-1)\cdot u(x,y)=0
\]
and determine solutions on $\C^2$ and its submanifolds.
A differential \textsc{Thomas} decomposition over $F\{u\}$ simply reproduces this equation, because $(xy-1)\in F\setminus\{0\}$.
However, we can model a search for solutions on submanifolds by adding two differential indeterminates $X$ and $Y$ to $U$ and consider the equations.
In order to allow splitting the manifold $\C^2$, we add two differential indeterminates $X$ and $Y$ to $U$ which model the $\Delta$-field elements $x$ and $y$.
Thus, we have to consider the additional equations
$(X_x\cdot(X_x-1))_=, (X_y)_=, (Y_y\cdot(Y_y-1))_=, (Y_x)_=$
together with the modified equation
$((XY-1)\cdot u)_=$.
A differential \textsc{Thomas} decomposition with $X,Y \ll u$ yields three systems:
\begin{enumerate}[(i)]
\item
$\displaystyle (XY-1)_=,\enskip (X_x)_=,\enskip (X_y)_=,\enskip (X)_{\neq}$\label{ex_subman_1}
\item
$\displaystyle (u)_=,\enskip (Y_x)_=,\enskip (Y_y\cdot(Y_y-1))_=,\enskip (X_x\cdot(X_x-1))_=,\enskip (X_y)_=,\enskip (X)_{\neq},\enskip (XY-1)_{\neq}$\label{ex_subman_2}
\item
$\displaystyle (u)_=,\enskip (Y_x)_=,\enskip (Y_y\cdot(Y_y-1))_=,\enskip (X)_=$\label{ex_subman_3}
\end{enumerate}
System (\ref{ex_subman_1}) allows an arbitrary function $u$ on the submanifold $M\subset\C^2$ defined by $xy-1=0$ as a solution.
The other systems (\ref{ex_subman_2}) and (\ref{ex_subman_3}) determine $u\equiv0$ as the only solution on $\C^2\setminus M$.
\end{myexample}

\section{Implementation}\label{implementation}

In this section, we describe our implementation of the decomposition algorithm.
First, we list some other implementations of triangular decomposition algorithms.
Second, we give some typical optimizations to make the computations feasible.
Third, we describe our implementation in {\sc Maple}.
Fourth, we give benchmarks to get a more detailed and practical comparison between different decomposition algorithms.

\subsection{Implementations of Similar Decomposition Algorithms}\label{other_implementation}

The \textsf{RegularChains} package by \cite{regularchains} is shipped with recent versions of \textsc{Maple}.
It contains the \textsf{Triangularize} command, which implements a decomposition of an algebraic variety given by a set of equations by means of regular chains.
If the input also contains inequations, the resulting decomposition is represented by regular systems instead.
It is possible to make these decompositions disjoint using the \textsf{MakePairwiseDisjoint} command.

The $\epsilon$\textsf{psilon} package by \cite{epsilon} implements different kinds of triangular decompositions in \textsc{Maple}.
It is the only software package besides our own that implements the algebraic \textsc{Thomas} decomposition.
It closely resembles the approach that \cite{Tho1,Tho2} suggested, i.e., polynomials of higher leader are considered first.
All polynomials of the same leader are combined into one common consequence, resulting in new conditions of lower leader.
These are not taken into account right away and will be treated in later steps.
Contrary to our approach, one cannot reduce modulo an \emph{unfinished} system.
Therefore, one needs extra inconsistency checks to avoid spending too much time on computations with inconsistent systems.
$\epsilon$\textsf{psilon} implements such checks in order to achieve good performance.

The \textsc{Maple} packages \textsf{diffalg} by \cite{diffalg} and \textsf{DifferentialAlgebra} by Boulier and Cheb-Terrab deal with ordinary and partial differential equations as described by \cite{BLOP}.
They compute a radical decomposition of a differential ideal, i.e., a description of the vanishing ideal of the \textsc{Kolchin} closure \cite[\textsection IV.1]{Kol} of the set of solutions.
Computation of integrability conditions is driven by reduction of $\Delta$-polynomials \cite[Sect. 2]{Rosenfeld}, which are the analogon of $s$-polynomials in differential algebra.
Just like in \textsf{RegularChains}, this approach usually does not give disjoint solution sets, although in principle disjointness might be achieved.
The \textsf{diffalg} package has been superseded by \textsf{DifferentialAlgebra} in \textsc{Maple} 14.
\textsf{DifferentialAlgebra} is based on the \textsf{BLAD}-libraries by \cite{blad} which have been designed as a set of stand-alone \textsf{C}-libraries with an emphasis on usability for non-mathematicians and extensive documentation.

\subsection{Algorithmic Optimizations}\label{section_algorithmic_opt}

In this subsection, we describe algorithmic optimizations helpful for a reasonably fast implementation of the \textsf{Decompose} algorithm.

In our algorithm, pseudo remainder sequences for the same pairs of polynomials are usually needed several times in different branches.
As these calculations are expensive in general, our implementation always keeps the results in memory and reuses them when the same pseudo remainder sequence is requested again to \emph{avoid repeated computations}.

Coefficient growth is a common problem in elimination.
Polynomials should be represented as compact as possible.
Once we know that the initial of a polynomial is non-zero, the \emph{content} of a polynomial (in the univariate sense) is non-zero, too.
Thus, every time an initial is added to the system as an inequation, we can divide the polynomial by its content.
Additionally, the multivariate content, which is an element of $F$, can be removed.

The reduction algorithms (\ref{algo_reduce}) and (\ref{algo_differential_reduce}) do not recognize that non-leading coefficients are zero.
However, we can reduce the coefficients modulo the polynomials of lower leader, in addition to reduction of the polynomial itself.
Thereby, in some cases the sizes of coefficients decrease, in other cases they increase.
The latter is partly due to multiplying the whole polynomials with the initials of the reductors.
Finding a good heuristic for this \emph{coefficient reduction} is crucial for efficiency.

\emph{Factorization} of a polynomial improves computation time in many cases.
More precisely, the system $S \uplus \{ (p\cdot q)_= \}$ decomposes disjointly into $(S \cup \{ p_= \}, S \cup \{ p_{\neq}, q_= \})$ and the system $S \cup \{ (p\cdot q)_{\neq} \}$ is equivalent to $S \cup \{ p_{\neq}, q_{\neq} \}$.
In most cases, the computation of two smaller problems resulting from a factorization is cheaper than the computation of the big, original problem.
This idea extends to factorizations over an extension of the base field:
Let $Y_i:=\left\{x_j\mid x_j<x_i,(S_T)^=_{x_j}\not=\emptyset\right\}$ and $Z_i:=\left\{x_j\mid x_j<x_i, (S_T)^=_{x_j}=\emptyset\right\}$.
Assume that $(S_T)^=_{x_i}$ is irreducible over the field $F_i:=F(Z_i)[Y_i]/\langle (S_T)^=_{<x_i}\rangle$ for all $i\in\{1,\ldots,n\}$, where $\langle (S_T)^=_{<x_i}\rangle$ is the ideal generated by $(S_T)^=_{<x_i}$ in the polynomial ring $F(Z_i)[Y_i]$.
Factorization over $F_n$ instead of $F$ may split the polynomial into more factors, but it is not clear whether this improves runtime.
Preliminary tests show that factorization over $F$ should be preferred for $F=\mathbb{Q}$.

In the algebraic algorithm, polynomials need not be square-free when they are inserted into the candidate simple system.
Efficiency can sometimes be improved by postponing the computation of the square-free split as long as possible.
However, this is not possible for the differential case.
Differential polynomials need to be made square-free to ensure that their separant is non-zero, i.e. non-trivial prolongations have a non-zero initial.

In the differential case, application of \emph{criteria} can decrease computation time by avoiding useless reductions of non-reductive prolongations.
\textsc{Janet}'s combinatorial approach already avoids many reductions of $\Delta$-polynomials, as used in other approaches (see \cite{GY}).
In addition, we use the involutive criteria 2-4 (cf. \cite{GB1,GerI,AH}), which together are equivalent to the chain criterion.
Applicability of this criterion in the non-linear differential case was shown in \cite[\textsection 4, Prop.\ 5]{BLOP}.

The axioms of a selection strategy (see Definition (\ref{Select})) already strongly limit the choice for the polynomial considered in the current step.
However, the remaining freedom is another important aspect for the speed of an actual implementation.
We will describe different selection strategies in \textsection\ref{section_selection_strat} and compare them in the benchmarks.

As described up to now, the algorithm often keeps on computing with inconsistent systems.
We want to optimize the algorithm to detect the inconsistencies as early as possible.
This allows the algorithm to discard inconsistent systems as early as possible.
One of the problems are selection strategies that postpone the costly treatment of inequations.
A test to detect whether inequations in $S_Q$ reduce to zero is comparably cheap.

Another possible improvement is parallelization, since the main loop in \textsf{Decompose} (\ref{algo_decompose}) can naturally be used in parallel for different systems.

\subsection{Selection strategies}\label{section_selection_strat}

We consider our two main approaches to selection strategies (see Definition (\ref{Select})).
\begin{enumerate}
 \item The ``equations first'' strategies:
  \textsf{Select} only chooses an inequation if $Q$ does not contain any equations.
  Among the equations or inequations, it prefers the ones with smallest leader.
 \item The ``leader first'' strategies:
  \textsf{Select} always chooses an equation or inequation with the smallest leader occurring in $Q$.
  If there are both equations and inequations with that leader, it chooses an equation.
\end{enumerate}
In both approaches, if the above criteria do not yield a unique choice, we compare the leader of the initial and choose the smaller one.
We apply the last test recursively to the initial of the initial and so on.
At this point, it is still possible that we fail to make a unique selection.
However, these cases are rare and there does not seem to be a considerable performance advantage for any choice.
Therefore, it suffices to make an arbitrary (but preferably unique) choice.

In our experimental observation ``leader first'' strategies usually produce decompositions with less systems, while ``equations first'' strategies are more efficient (cf.\ \textsection\ref{sect_benchmarks}).

\subsection{Implementation in {\sc Maple}}

Both the algebraic and the differential case of the \textsc{Thomas} decomposition algorithm have been implemented in the \textsc{Maple} computer algebra system.
Packages can be downloaded from our web page (\cite{hp_thomasdecomp}), documentation and example worksheets are available there.

The main reason for choosing \textsc{Maple} for the implementation is the collection of solvers for polynomial equations, ODEs, and PDEs already present.
Furthermore, fast algorithms exist for polynomial factorization over finitely generated field extensions of $\mathbb{Q}$ and for gcd computation.

The \textsf{AlgebraicThomas} package includes procedures to compute a \textsc{Thomas} decomposition, reduce polynomials modulo simple systems and compute counting polynomials (cf.~\cite{PleskenCounting}).
Furthermore, it can represent the complement and intersection of solution sets as decompositions into simple systems.
Finally, a comprehensive \textsc{Thomas} decomposition can be computed, this topic will be discussed in a later publication.

\begin{myexample}
\begin{maplegroup}
We demonstrate how to use the \textsf{AlgebraicThomas} package by
computing a decomposition of the system in example
(\ref{ex_mitternacht}).

\end{maplegroup}
\begin{maplegroup}
\begin{mapleinput}
\mapleinline{active}{1d}{with(AlgebraicThomas):}{%
}
\end{mapleinput}

\end{maplegroup}
\begin{maplegroup}
\begin{mapleinput}
\mapleinline{active}{1d}{p := a*x^2 + b*x + c;}{%
}
\end{mapleinput}

\mapleresult
\begin{maplelatex}
\mapleinline{inert}{2d}{p := x^2*a+x*b+c;}{%
\[
p := x^{2}\,a + x\,b + c
\]
}
\end{maplelatex}

\end{maplegroup}
\begin{maplegroup}
\begin{mapleinput}
\mapleinline{active}{1d}{S := AlgebraicThomasDecomposition([p], [x,c,b,a]);}{%
}
\end{mapleinput}

\mapleresult
\begin{maplelatex}
\mapleinline{inert}{2d}{S := [[x^2*a+x*b+c = 0, 4*c*a-b^2 <> 0, a <> 0], [2*x*a+b = 0,
4*c*a-b^2 = 0, a <> 0], [x*b+c = 0, b <> 0, a = 0], [c = 0, b = 0, a =
0]];}{%
\maplemultiline{
S := [[x^{2}\,a + x\,b + c=0, \,4\,c\,a - b^{2}\neq 0, \,a\neq 0]
, \,[2\,x\,a + b=0, \,4\,c\,a - b^{2}=0, \,a\neq 0],  \\
[x\,b + c=0, \,b\neq 0, \,a=0], \,[c=0, \,b=0, \,a=0]] }
}
\end{maplelatex}

\end{maplegroup}
\begin{maplegroup}
Information about leader and main degree can optionally be included in the
output.

\end{maplegroup}
\begin{maplegroup}
\begin{mapleinput}
\mapleinline{active}{1d}{map(printSystem, S, ["PT", "LR"]);}{%
}
\end{mapleinput}

\mapleresult
\begin{maplelatex}
\mapleinline{inert}{2d}{[[[x^2*a+x*b+c = 0, x^2], [4*c*a-b^2 <> 0, c], [a <> 0, a]],
[[2*x*a+b = 0, x], [4*c*a-b^2 = 0, c], [a <> 0, a]], [[x*b+c = 0, x],
[b <> 0, b], [a = 0, a]], [[c = 0, c], [b = 0, b], [a = 0, a]]];}{%
\maplemultiline{
[[[x^{2}\,a + x\,b + c=0, \,x^{2}], \,[4\,c\,a - b^{2}\neq 0, \,c
], \,[a\neq 0, \,a]],  \\
[[2\,x\,a + b=0, \,x], \,[4\,c\,a - b^{2}=0, \,c], \,[a\neq 0, \,
a]],  \\
[[x\,b + c=0, \,x], \,[b\neq 0, \,b], \,[a=0, \,a]], \,[[c=0, \,c
], \,[b=0, \,b], \,[a=0, \,a]]] }
}
\end{maplelatex}

\end{maplegroup}
\begin{maplegroup}
It is possible to include inequations in the input to exclude some
degenerate cases:

\end{maplegroup}
\begin{maplegroup}
\begin{mapleinput}
\mapleinline{active}{1d}{q := a<>0;}{%
}
\end{mapleinput}

\mapleresult
\begin{maplelatex}
\mapleinline{inert}{2d}{q := a <> 0;}{%
\[
q := a\neq 0
\]
}
\end{maplelatex}

\end{maplegroup}
\begin{maplegroup}
\begin{mapleinput}
\mapleinline{active}{1d}{T := AlgebraicThomasDecomposition([p, q], [x,c,b,a]);}{%
}
\end{mapleinput}

\mapleresult
\begin{maplelatex}
\mapleinline{inert}{2d}{T := [[x^2*a+x*b+c = 0, 4*c*a-b^2 <> 0, a <> 0], [2*x*a+b = 0,
4*c*a-b^2 = 0, a <> 0]];}{%
\[
T := [[x^{2}\,a + x\,b + c=0, \,4\,c\,a - b^{2}\neq 0, \,a\neq 0]
, \,[2\,x\,a + b=0, \,4\,c\,a - b^{2}=0, \,a\neq 0]]
\]
}
\end{maplelatex}

\end{maplegroup}
\begin{maplegroup}
\end{maplegroup}
\end{myexample}

Features for the differential package \textsf{DifferentialThomas} include arbitrary differential rankings, using special functions implemented in \textsc{Maple} as differential field elements, computation of power series solutions, and a direct connection to the solvers of \textsc{Maple} for differential equations.
\begin{myexample}
\begin{maplegroup}
We treat the following control theoretic example taken from \cite{diop}.

\end{maplegroup}
\begin{maplegroup}
\emptyline
\end{maplegroup}
\begin{maplegroup}
\begin{mapleinput}
\mapleinline{active}{1d}{with(DifferentialThomas):
ComputeRanking([t],[x2,x1,y,u],"EliminateFunction");}{%
}
\end{mapleinput}

\end{maplegroup}
\begin{maplegroup}
\emptyline
\end{maplegroup}
\begin{maplegroup}
This creates the differential polynomial ring
$\mathbb{Q}\{x^{(2)},x^{(1)},y,u\}$ for
$\Delta=\{\frac{\partial}{\partial t}\}$. Here $u$ indicates the
input, $x^{(1)}$ and $x^{(2)}$ the state, and $y$ the output of the
system. The chosen ranking ``$<$'' is the elimination ranking with
$x^{(2)}\gg x^{(1)}\gg y\gg u$, i.e.,
$x^{(2)}_\mathbf{i}>x^{(1)}_\mathbf{j}>y_\mathbf{k}>u_\mathbf{l}$ for
all $\mathbf{i},\mathbf{j},\mathbf{k},\mathbf{l}\in \mathbb{Z}_{\ge
0}$.

\end{maplegroup}
\begin{maplegroup}
\emptyline
\end{maplegroup}
\begin{maplegroup}
\begin{mapleinput}
\mapleinline{active}{1d}{L:=[x1[1]-u[0]*x2[0],x2[1]-x1[0]-u[0]*x2[0],y[0]-x1[0]]:}{%
}
\end{mapleinput}

\end{maplegroup}
\begin{maplegroup}
\emptyline
\end{maplegroup}
\begin{maplegroup}
We follow \cite[Ex.~1]{diop} and compute the external
trajectories of a differential ideal generated by $L$, i.e. intersect
this differential ideal with $\mathbb{Q}\{y,u\}$. 

\end{maplegroup}
\begin{maplegroup}
\emptyline
\end{maplegroup}
\begin{maplegroup}
\begin{mapleinput}
\mapleinline{active}{1d}{res:=DifferentialThomasDecomposition(L,[]);}{%
}
\end{mapleinput}

\mapleresult
\begin{maplelatex}
\mapleinline{inert}{2d}{res := [DifferentialSystem, DifferentialSystem];}{%
\[
\mathit{res} := [\mathit{DifferentialSystem}, \,\mathit{
DifferentialSystem}]
\]
}
\end{maplelatex}

\end{maplegroup}
\begin{maplegroup}
\emptyline
\end{maplegroup}
\begin{maplegroup}
We show the equations and inequations of the differential systems not
involving $x^{(1)}$ or $x^{(2)}$.
The chosen ranking guarantees that the systems shown determine the
external trajectories of the system:

\end{maplegroup}
\begin{maplegroup}
\emptyline
\end{maplegroup}
\begin{maplegroup}
\begin{mapleinput}
\mapleinline{active}{1d}{PrettyPrintDifferentialSystem(res[1]):remove(a->has(a,[x1,x2]),\%);}{
}
\end{mapleinput}

\mapleresult
\begin{maplelatex}
\mapleinline{inert}{2d}{[-u(t)*diff(y(t),`$`(t,2))+diff(y(t),t)*u(t)^2+diff(y(t),t)*diff(u(t)
,t)+y(t)*u(t)^2 = 0, u(t) <> 0];}{%
\[
[ - \mathrm{u}(t)\,({\frac {d^{2}}{dt^{2}}}\,\mathrm{y}(t)) + (
{\frac {d}{dt}}\,\mathrm{y}(t))\,\mathrm{u}(t)^{2} + ({\frac {d}{
dt}}\,\mathrm{y}(t))\,({\frac {d}{dt}}\,\mathrm{u}(t)) + \mathrm{
y}(t)\,\mathrm{u}(t)^{2}=0, \,\mathrm{u}(t)\neq 0]
\]
}
\end{maplelatex}

\end{maplegroup}
\begin{maplegroup}
\begin{mapleinput}
\mapleinline{active}{1d}{PrettyPrintDifferentialSystem(res[2]):remove(a->has(a,[x1,x2]),\%);}{
}
\end{mapleinput}

\mapleresult
\begin{maplelatex}
\mapleinline{inert}{2d}{[diff(y(t),t) = 0, u(t) = 0];}{%
\[
[{\frac {d}{dt}}\,\mathrm{y}(t)=0, \,\mathrm{u}(t)=0]
\]
}
\end{maplelatex}

\end{maplegroup}
\begin{maplegroup}
\emptyline
\end{maplegroup}
\begin{maplegroup}
These systems, having disjoint solution sets, are identical to the
ones found in \cite{diop}.

\end{maplegroup}

\end{myexample}

\subsection{Benchmarks}\label{sect_benchmarks}

In this subsection, we compare our two \textsc{Maple} packages to the other implementations of triangular decompositions mentioned in \textsection\ref{other_implementation} using benchmarks.
Not all of the implementations compute equivalent results.
This should be considered when comparing the timings.
We omitted examples where all tested systems took less than one second to complete the computation or could not be computed by any software package.

All benchmarks have been performed with Linux x86-64 running on a third generation Opteron, $2.3$ GHz.
The time limit has been set to $3$ hours and available memory is limited to $4$ GB.
All times are given in seconds.
The polynomial multiplication in \textsc{Maple} 14 benefits from a new parallel implementation (cf.~\cite{mope}).
Nonetheless, we state the total CPU time in our benchmarks, as returned by \textsc{Maple}'s \texttt{time} command.

By default, both of our \textsc{Maple} packages behave as follows:
\begin{itemize}
 \item Polynomials are factorized over $\mathbb{Q}$.
 \item The content of polynomials is removed.
 \item The selection strategy is an ``equations first'' strategy, as decribed in \textsection\ref{section_selection_strat}.
 \item After reducing a polynomial, we always reduce its coefficients fully.
 \item Inequations in $S_Q$ are reduced for early inconsistency checks.
\end{itemize}
See \textsection\ref{section_algorithmic_opt} for details.

\subsubsection{Algebraic Systems}

For testing the \textsf{AlgebraicThomas} package, we used two sets of examples, namely, the test examples from the \texttt{polsys50} file in Wang's $\epsilon$psilon package (\cite{epsilon}) printed in table \ref{benchmark_t1} and the examples from \cite{ctd} as shown in table \ref{benchmark_t2}.

In contrast to Algorithm \ref{algo_decompose}, the implementation in the \textsf{AlgebraicThomas} package inserts equations or inequations into $S_T$ without making them square-free first.
It delays this computation as long as possible, sometimes until the end of the decomposition.
This avoids some expensive and unnecessary discriminant computations entirely.

\begin{table}
\caption{Comparison of algebraic decompositions 1: \texttt{polsys50} from \cite{epsilon}}
\label{benchmark_t1}
\begin{center}{\footnotesize
\begin{tabular}{|c||c|c|c|c||c|c|c|c|c|}
\hline
Name & RC1 & RC2 & RC3 & DW1 & DW2 & AT1 & AT2 & AT3 & AT4 \\                                                                                           
\hline
1 & $3.5$ & $3.7$ & $4.3$ & $0.4$ & $1.0$ & $3.0$ & $1.1$ & $1.8$ & $1.4$ \\
2 & $7.4$ & $6.7$ & $7.5$ & $7.6$ & $8.4$ & $7.1$ & $169.7$ & $95.8$ & $6.6$ \\
3 & $>3h$ & $>3h$ & $>4GB$ & $985.7$ & $1344.6$ & $7538.0$ & $>4GB$ & $>4GB$ & $194.6$ \\
4 & $>4GB$ & $>4GB$ & $>4GB$ & $>4GB$ & $>4GB$ & $0.2$ & $>4GB$ & $>4GB$ & $32.1$ \\
6 & $0.4$ & $0.4$ & $47.2$ & $0.1$ & $0.2$ & $0.2$ & $0.1$ & $0.2$ & $0.1$ \\
7 & $>3h$ & $>3h$ & $>3h$ & $7352.6$ & $>3h$ & $>4GB$ & $>4GB$ & $>4GB$ & $>4GB$ \\
12 & $0.5$ & $0.6$ & $0.5$ & $0.3$ & $0.4$ & $0.4$ & $0.6$ & $1.1$ & $0.4$ \\
14 & $0.5$ & $2.3$ & $0.6$ & $>3h$ & $>4GB$ & $1.5$ & $1.6$ & $1.4$ & $2.5$ \\
16 & $0.9$ & $0.9$ & $1.0$ & $1.4$ & $1.5$ & $1.8$ & $5.6$ & $>3h$ & $2.2$ \\
17 & $6.5$ & $6.4$ & $13.0$ & $4.7$ & $6.3$ & $75.5$ & $12076.5$ & $>3h$ & $12.6$ \\
18 & $0.3$ & $0.3$ & $3.7$ & $0.1$ & $0.1$ & $0.1$ & $0.1$ & $0.1$ & $0.1$ \\
19 & $419.9$ & $452.9$ & $>4GB$ & $0.4$ & $0.6$ & $0.4$ & $5842.5$ & $0.4$ & $0.3$ \\
21 & $1.6$ & $1.9$ & $2.1$ & $86.6$ & $>4GB$ & $4.5$ & $>3h$ & $4.4$ & $112.8$ \\
22 & $0.6$ & $0.6$ & $0.6$ & $1.2$ & $1.6$ & $1.5$ & $2.9$ & $32.4$ & $2.0$ \\
23 & $0.4$ & $0.7$ & $0.4$ & $0.1$ & $>4GB$ & $29.5$ & $>3h$ & $>4GB$ & $29.0$ \\
24 & $1.2$ & $1.1$ & $1.3$ & $1.3$ & $2.6$ & $1.0$ & $2.0$ & $4.5$ & $1.6$ \\
25 & $1.2$ & $8.5$ & $1.6$ & $>3h$ & $>4GB$ & $>4GB$ & $>3h$ & $>3h$ & $>3h$ \\
29 & $0.3$ & $0.5$ & $0.4$ & $0.3$ & $0.3$ & $0.3$ & $55.2$ & $0.3$ & $0.3$ \\
30 & $>4GB$ & $>4GB$ & $>4GB$ & $>4GB$ & $>3h$ & $45.3$ & $42.9$ & $40.8$ & $>4GB$ \\
31 & $>4GB$ & $>4GB$ & $>4GB$ & $>4GB$ & $>3h$ & $>3h$ & $>4GB$ & $>3h$ & $>3h$ \\
33 & $3.4$ & $3.6$ & $3.2$ & $1.3$ & $1.3$ & $3.5$ & $66.9$ & $15.2$ & $1.1$ \\
34 & $911.5$ & $916.9$ & $926.5$ & $>3h$ & $>4GB$ & $>4GB$ & $>4GB$ & $>3h$ & $>4GB$ \\
35 & $1.5$ & $1.5$ & $1.6$ & $1.2$ & $1.3$ & $1.7$ & $4.9$ & $7.3$ & $0.5$ \\
39 & $0.6$ & $0.7$ & $0.8$ & $1.2$ & $1.9$ & $0.6$ & $1.0$ & $8.0$ & $0.5$ \\
41 & $1.5$ & $1.5$ & $1.6$ & $1.5$ & $1.7$ & $7.0$ & $1.4$ & $0.6$ & $114.5$ \\
43 & $0.7$ & $0.7$ & $0.7$ & $3.1$ & $4.4$ & $0.2$ & $1.0$ & $0.2$ & $0.2$ \\
44 & $24.5$ & $17.2$ & $24.1$ & $3.4$ & $4.2$ & $1.2$ & $1.7$ & $0.8$ & $>4GB$ \\
47 & $1.3$ & $1.7$ & $1.4$ & $2.8$ & $6.6$ & $13.0$ & $>3h$ & $11.1$ & $92.4$ \\
49 & $0.3$ & $0.3$ & $0.3$ & $610.2$ & $32.1$ & $0.5$ & $>3h$ & $0.5$ & $0.5$ \\
\hline
\end{tabular}
}\end{center}
\end{table}

We compared \textsf{AlgebraicThomas} with the \textsf{RegularChains} package from \textsc{Maple} 14 and $\epsilon$psilon.
We also tested the \textsf{AlgebraicThomas} and \textsf{RegularChains} packages in different configurations.
The timings in \textsc{Maple} 14 of the following procedures are being compared:
\begin{itemize}
 \item (RC1) \texttt{RegularChains[Triangularize]}.
 \item (RC2) \texttt{RegularChains[Triangularize]} with the \texttt{'output'='lazard'} option set.
 \item (RC3) \texttt{RegularChains[Triangularize]} with the \texttt{'radical'='yes'} option set.
 \item (DW1) \texttt{epsilon[RegSer]}.
 \item (DW2) \texttt{sisys[simser]}.
 \item (AT1) \texttt{AlgebraicThomasDecomposition}.
 \item (AT2) \texttt{AlgebraicThomasDecomposition} with factorization disabled.
 \item (AT3) \texttt{AlgebraicThomasDecomposition} with a ``leader first'' selection strategy (cf.~\textsection\ref{section_selection_strat}).
 \item (AT4) \texttt{AlgebraicThomasDecomposition} with coefficient reduction disabled.
\end{itemize}

We compare table \ref{benchmark_t1} and \ref{benchmark_t2} within our own implementation.
We observe, that (AT2) is much slower than (AT1) and, thus, conclude that factorization is vital to make many computations feasible.
In a few examples, we see the relative advantage of the default selection strategy compared to the one used in (AT3).
Generally speaking, disabling coefficient reduction increases computation time for (AT4), but there are some strong counterexamples to this observation.
This indicates that different strategies for coefficient reduction, as seen in (AT1) and (AT4), should be investigated further.

The programs \texttt{sisys[simser]} (DW2) and \texttt{AlgebraicThomasDecomposition} (AT1-4) are the only ones that compute a \textsc{Thomas} decomposition.
All test examples that could be computed by (DW2) could also be computed by (AT1).
However, there are some examples that \texttt{RegularChains} (RC1) or \texttt{epsilon} (DW1) could treat, but we could not decompose into simple systems.
Moreover, the test examples indicate that (RC1) is in general faster than (AT1) in the positive-dimensional case.
Our evaluation suggests that this is due to the strict square-free property of simple systems.
In the zero-dimensional case, however, the situation is less clear, since there are examples where (RC1) is faster than (AT1) and vice versa.

\begin{table}
\caption{Comparison of algebraic decompositions 2: Test examples from \cite{ctd}}
\label{benchmark_t2}
\begin{center}{\footnotesize
\begin{tabular}{|c||c|c|c|c||c|c|c|c|c|}
\hline
Name & RC1 & RC2 & RC3 & DW1 & DW2 & AT1 & AT2 & AT3 & AT4 \\                                                                                           
\hline
AlkashiSinus & $0.6$ & $0.6$ & $0.8$ & $0.1$ & $7.1$ & $5.7$ & $2.6$ & $6.5$ & $3.6$ \\
Bronstein & $0.4$ & $0.5$ & $0.5$ & $0.2$ & $0.4$ & $0.3$ & $0.4$ & $1.1$ & $0.4$ \\
Cheaters-homotopy-easy & $0.7$ & $>3h$ & $532.5$ & $>4GB$ & $>4GB$ & $>3h$ & $>4GB$ & $>3h$ & $>3h$ \\
Cheaters-homotopy-hard & $0.7$ & $>3h$ & $559.8$ & $>4GB$ & $>4GB$ & $>3h$ & $>4GB$ & $>3h$ & $>3h$ \\
Gerdt & $1.4$ & $1.4$ & $1.4$ & $2.0$ & $2.2$ & $8.1$ & $3.2$ & $0.5$ & $1532.1$ \\
Hereman-2 & $0.8$ & $1.0$ & $0.8$ & $0.3$ & $0.4$ & $0.3$ & $1.2$ & $0.3$ & $0.5$ \\
Hereman-8-8 & $26.9$ & $31.6$ & $208.3$ & $>3h$ & $>3h$ & $>3h$ & $>3h$ & $>3h$ & $>4GB$ \\
KdV & $722.2$ & $707.1$ & $725.7$ & $>3h$ & $>3h$ & $>3h$ & $>3h$ & $>3h$ & $>3h$ \\
Lanconelli & $0.4$ & $0.6$ & $0.4$ & $0.2$ & $0.4$ & $0.4$ & $1.3$ & $0.3$ & $0.3$ \\
Lazard-ascm2001 & $1.2$ & $17.5$ & $1.4$ & $>3h$ & error & $>4GB$ & $>4GB$ & $>3h$ & $>3h$ \\
Leykin-1 & $5.6$ & $8.0$ & $5.8$ & $>3h$ & $>3h$ & $2.3$ & $>3h$ & $6.0$ & $1.4$ \\
Maclane & $2.4$ & $6.7$ & $2.6$ & $3576.5$ & $>4GB$ & $7.4$ & $17.4$ & $13.1$ & $7.0$ \\
MontesS10 & $0.5$ & $1.0$ & $0.6$ & $0.4$ & $17.7$ & $2.0$ & $2.2$ & $2.3$ & $1.5$ \\
MontesS11 & $0.2$ & $0.5$ & $0.2$ & $0.2$ & $>3h$ & $23.5$ & $>3h$ & $21.9$ & $12.4$ \\
MontesS12 & $0.5$ & $3.2$ & $0.5$ & $1.6$ & $>3h$ & $9.4$ & $31.0$ & $13.6$ & $115.4$ \\
MontesS13 & $0.3$ & $0.5$ & $0.3$ & $0.2$ & $0.5$ & $0.8$ & $1.8$ & $1.1$ & $0.9$ \\
MontesS14 & $0.7$ & $1.3$ & $0.8$ & $>3h$ & $>3h$ & $6.0$ & $>4GB$ & $14.5$ & $12.1$ \\
MontesS15 & $1.2$ & $1.8$ & $1.3$ & $0.6$ & $8.2$ & $4.8$ & $3.8$ & $6.7$ & $3.3$ \\
MontesS16 & $4.3$ & $3.4$ & $4.2$ & $1.5$ & $1.5$ & $2.5$ & $2.2$ & $3.2$ & $1.5$ \\
MontesS7 & $0.4$ & $0.5$ & $0.4$ & $0.2$ & $0.5$ & $0.7$ & $0.6$ & $2.5$ & $1.2$ \\
Neural & $0.5$ & $0.7$ & $0.6$ & $>3h$ & $>4GB$ & $1.4$ & $3050.9$ & $1.7$ & $1.2$ \\
Pavelle & $1.1$ & $15.9$ & $1.4$ & $>3h$ & $>4GB$ & $>3h$ & $>3h$ & $>3h$ & $>3h$ \\
Wang93 & $1.2$ & $1.3$ & $1.2$ & $1.6$ & $3.4$ & $4.9$ & $>3h$ & $3.4$ & $6.5$ \\
genLinSyst-3-2 & $0.3$ & $1.1$ & $0.3$ & $0.2$ & $0.2$ & $0.3$ & $0.2$ & $0.3$ & $0.2$ \\
genLinSyst-3-3 & $0.3$ & $4.5$ & $0.4$ & $1.2$ & $1.2$ & $6.0$ & $2.5$ & $4.8$ & $1.1$ \\
\hline
\end{tabular}
}\end{center}
\end{table}

\subsubsection{Differential Systems}

\begin{table}
\caption{Benchmarks for ODE systems}
\label{benchmark_ode}
\begin{center}{\footnotesize
\begin{tabular}{|c||c|c||c|c|c|c|}
\hline
Name & DA  & da  & DT1  & DT2  & DT3  & DT4 \\
\hline
 Diffalg4 & 2.9 & 2.9 & 852.5 & $>3h$ & 8932.4 & 36.0\\
 LLG3 & 0.5 & $>3h$ & 5.4 & 5.6 & 4.4 & 4.9\\
 LLG4 & 0.3 & 19.1 & 2.6 & 37.4 & 20.3 & 4.0\\
 ODE1 & 2.4 & 3.7 & 0.6 & 0.3 & 0.6 & 0.8\\
 ODE6 & 2.3 & 1.5 & 0.8 & 1.2 & 0.6 & 0.8\\
 ODE7 & $>3h$ & $>3h$ & 3.2 & 60.8 & 47.2 & 5.4\\
 kepler vs newton & 0.7 & 0.9 & 1.4 & 2.8 & 0.8 & 2.0\\
 keppler1 & 0.1 & 0.2 & 1.1 & 0.6 & 0.8 & 1.1\\
 keppler2 & 0.1 & 0.1 & 1.3 & 1.5 & 1.1 & 0.8\\
 keppler3 & 0.1 & 0.2 & 1.1 & 0.6 & 0.7 & 1.1\\
 murray1 & 0.1 & 0.4 & 1.9 & 2.4 & 1.7 & 2.2\\
 murray2 & 0.1 & 0.1 & 0.7 & 1.4 & 0.8 & 0.6\\
\hline
\end{tabular}
}\end{center}
\end{table}

We compared \textsf{DifferentialThomas} with the packages \textsf{diffalg} and \textsf{DifferentialAlgebra}.
Finding a suitable set of benchmark examples for the differential case was more difficult.
We are not aware of any sets of standard benchmarks.
Thus, we used a collection of examples, which we came across in our work.
These examples are published on our homepage (\cite{hp_thomasdecomp}).

The timings of the following procedures are being compared:
\begin{itemize}
 \item (DA) \texttt{DifferentialAlgebra[Rosenfeld\_Groebner]}.
 \item (da) \texttt{diffalg[Rosenfeld\_Groebner]\footnote{with \texttt{\_env\_diffalg\_uses\_DifferentialAlgebra:=false}}}.
 \item (DT1) \texttt{DifferentialThomasDecomposition}.
 \item (DT2) \texttt{DifferentialThomasDecomposition} with factorization disabled.
 \item (DT3) \texttt{DifferentialThomasDecomposition} with a ``leader first'' selection strategy (cf.~\textsection\ref{section_selection_strat})
 \item (DT4) \texttt{DifferentialThomasDecomposition} with coefficient reduction disabled.
\end{itemize}

\begin{table}
\caption{Benchmarks for PDE systems}
\label{benchmark_pde}
\begin{center}{\footnotesize
\begin{tabular}{|c||c|c||c|c|c|c|}
\hline
Name & DA  & da  & DT1  & DT2  & DT3  & DT4 \\
\hline
 Cyclic 5 variant1 & $>3h$ & 2.9 & 1.3 & 1.5 & 1.4 & 1.2\\
 Cyclic 5 variant2 & $>3h$ & $>3h$ & 2.3 & 2.6 & 0.7 & 2.3\\
 Diffalg2 & 0.5 & 0.3 & 1.4 & 1.5 & 41.5 & 1.2\\
 Diffalg3 & 0.2 & 0.5 & 0.9 & 0.6 & 1.6 & 0.7\\
 Ibragimov 2, 17.9c & $>3h$ & $>3h$ & 18.4 & $>3h$ & 40.8 & 12.3\\
 PDE6 & $>4GB$ & $>4GB$ & 11.1 & 23.0 & $>4GB$ & 16.5\\
 PDE7 & $>4GB$ & 116.7 & 91.3 & 83.3 & $>4GB$ & 41.4\\
 PDE8 & $>4GB$ & $>4GB$ & 6.8 & $>4GB$ & 14.2 & 7.5\\
 Riquier 1b & 0.1 & 0.1 & 1.9 & 1.9 & 1.9 & 2.0\\
 Riquier 3a & 0.3 & 0.2 & 0.8 & 1.1 & 0.5 & 0.6\\
 Riquier 3b & 0.6 & 0.6 & 1.4 & 1.7 & 1.4 & 1.2\\
 boulier & $>3h$ & 546.4 & 2.5 & 1690.9 & 2.6 & 1.5\\
 cyclic 6 & $>4GB$ & 571.9 & 160.7 & 159.8 & 349.6 & 154.1\\
 noon6 & $>4GB$ & 72.6 & 40.6 & 36.8 & 63.7 & 31.0\\
\hline
\end{tabular}
}\end{center}
\end{table}

We want to mention one further example not included in the benchmark table.
It is the test example $5$ of \texttt{diffalg}.
None of the packages in their default setting could compute this example.
Still, \texttt{diffalg} and \texttt{DifferentialAlgebra} were able to do so instantaneously by a change of ordering strategy.

The comparison between (DT1), (DT2) and (DT3) is similar to the algebraic case.
In particular, factorization should be enabled and the default selection strategy should be preferred.
In contrast to the algebraic implementation, the comparison of (DT1) and (DT4) is less conclusive.

All test examples which could be computed by \texttt{DifferentialAlgebra} or \texttt{diffalg} could also be computed by our default strategy (DT1).
For ODEs, the three packages show similar timings, but for PDEs, \texttt{DifferentialThomasDecomposition} appears to be faster.
This might be explained by the involutive approach, which we utilize to make the subsystems coherent.
A similar result can be found for the GINV-project (cf.\ \cite{GINV}, \cite{GINV_Bench}).

\section{Acknowledgments}\label{section_acknowledgements}

The contents of this paper profited very much from numerous useful comments and remarks by Wilhelm Plesken.
The authors thank him as well as Dongming Wang, Fran{\c{c}}ois Boulier and Fran{\c{c}}ois Lemaire for fruitful discussions.
The second author (V.P.G.) acknowledges the Deutsche Forschungsgemeinschaft for the financial support that made his stay in Aachen possible.
The presented results were obtained during his visits.
The contribution of the first and third author (T.B. and M.L.-H.) was partly supported by Schwerpunkt SPP 1489 of the Deutsche Forschungsgemeinschaft.

Finally, our gratitude goes to the anonymous referees for valuable comments, for pointing out informative references and for taking the time to verify the results of our benchmarks.

\bibliographystyle{elsarticle-harv}

\begin{thebibliography}{56}
\expandafter\ifx\csname natexlab\endcsname\relax\def\natexlab#1{#1}\fi
\expandafter\ifx\csname url\endcsname\relax
  \def\url#1{\texttt{#1}}\fi
\expandafter\ifx\csname urlprefix\endcsname\relax\def\urlprefix{URL }\fi

\bibitem[{Apel and Hemmecke(2005)}]{AH}
Apel, J., Hemmecke, R., 2005. Detecting unnecessary reductions in an involutive
  basis computation. J. Symbolic Comput. 40~(4-5), 1131--1149.
\newline\urlprefix\url{http://dx.doi.org/10.1016/j.jsc.2004.04.004}

\bibitem[{Aubry et~al.(1999)Aubry, Lazard, and Moreno~Maza}]{AubryTriTheo}
Aubry, P., Lazard, D., Moreno~Maza, M., 1999. On the theories of triangular
  sets. J. Symbolic Comput. 28~(1-2), 105--124, {Polynomial
  elimination---algorithms and applications}.
\newline\urlprefix\url{http://dx.doi.org/10.1006/jsco.1999.0269}

\bibitem[{B\"achler et~al.(2010)B\"achler, Gerdt, Lange-Hegermann, and
  Robertz}]{thomasalg_casc}
B\"achler, T., Gerdt, V., Lange-Hegermann, M., Robertz, D., 2010. Thomas
  decomposition of algebraic and differential systems. In: Computer Algebra in
  Scientific Computing. Tsakhkadzor, Armenia, pp. 31--54.

\bibitem[{B{\"a}chler and Lange-Hegermann(2008-2011)}]{hp_thomasdecomp}
B{\"a}chler, T., Lange-Hegermann, M., 2008-2011. \textsf{AlgebraicThomas} and
  \textsf{DifferentialThomas}: Thomas Decomposition for algebraic and
  differential systems.
  (\url{http://wwwb.math.rwth-aachen.de/thomasdecomposition/}).

\bibitem[{Blinkov et~al.(2003)Blinkov, Cid, Gerdt, Plesken, and
  Robertz}]{JanetPackage}
Blinkov, Y.~A., Cid, C.~F., Gerdt, V.~P., Plesken, W., Robertz, D., 2003. {The
  $\mathsf{MAPLE}$ Package {\sc Janet}: {\sc I}. Polynomial Systems. {\sc II}.
  Linear Partial Differential Equations}. In: Proc. 6th Int. Workshop on
  Computer Algebra in Scientific Computing, Passau, Germany. pp. 31--40 and
  41--54, (\url{http://wwwb.math.rwth-aachen.de/Janet}).

\bibitem[{Blinkov et~al.(2010{\natexlab{a}})Blinkov, Gerdt, and Robertz}]{GINV}
Blinkov, Y.~A., Gerdt, V.~P., Robertz, D., 2010{\natexlab{a}}.
  \textsf{GINV}-project. {Gröbner bases constructed by involutive algorithms},
  \href{http://invo.jinr.ru/ginv/index.html}{\tt
  http://invo.jinr.ru/ginv/index.html}.

\bibitem[{Blinkov et~al.(2010{\natexlab{b}})Blinkov, Gerdt, and
  Robertz}]{GINV_Bench}
Blinkov, Y.~A., Gerdt, V.~P., Robertz, D., 2010{\natexlab{b}}.
  \textsf{GINV}-project benchmarks. {Standard benchmarks between Gröbner bases
  systems}, \href{http://cag.jinr.ru/wiki/Benchmarks}{\tt
  http://cag.jinr.ru/wiki/Benchmarks}.

\bibitem[{Boulier(2004-2009)}]{blad}
Boulier, F., 2004-2009. \textsf{BLAD}: Bibliothèques Lilloises d'Algèbre
  Différentielle. (\url{http://www.lifl.fr/~boulier/BLAD/}).

\bibitem[{Boulier(2007)}]{BoulierElimination}
Boulier, F., 2007. Differential elimination and biological modelling. In:
  Gr\"obner bases in symbolic analysis. Vol.~2 of Radon Ser. Comput. Appl.
  Math. Walter de Gruyter, Berlin, pp. 109--137.

\bibitem[{Boulier and Hubert(1996-2004)}]{diffalg}
Boulier, F., Hubert, E., 1996-2004. \textsf{DIFFALG}: description, help pages
  and examples of use. Symbolic Computation Group, University of Waterloo,
  Ontario, Canada
  (\url{http://www-sop.inria.fr/members/Evelyne.Hubert/diffalg/}).

\bibitem[{Boulier et~al.(1995)Boulier, Lazard, Ollivier, and Petitot}]{BLOP2}
Boulier, F., Lazard, D., Ollivier, F., Petitot, M., 1995. {Representation for
  the radical of a finitely generated differential ideal}. In: {ISSAC'95:
  Proceedings of the 1995 International Symposium on Symbolic and Algebraic
  Computation}. ACM Press, New York, NY, USA, pp. 158--166,
  \url{http://hal.archives-ouvertes.fr/hal-00138020}.

\bibitem[{Boulier et~al.(2009)Boulier, Lazard, Ollivier, and Petitot}]{BLOP}
Boulier, F., Lazard, D., Ollivier, F., Petitot, M., 2009. Computing
  representations for radicals of finitely generated differential ideals. Appl.
  Algebra Engrg. Comm. Comput. 20~(1), 73--121.
\newline\urlprefix\url{http://dx.doi.org/10.1007/s00200-009-0091-7}

\bibitem[{Bourbaki(1968)}]{Bou68}
Bourbaki, N., 1968. Elements of mathematics. {T}heory of sets. Translated from
  the French. Hermann, Publishers in Arts and Science, Paris.

\bibitem[{Bouziane et~al.(2001)Bouziane, Kandri~Rody, and
  Ma{\^a}rouf}]{Bouziane}
Bouziane, D., Kandri~Rody, A., Ma{\^a}rouf, H., 2001. Unmixed-dimensional
  decomposition of a finitely generated perfect differential ideal. J. Symbolic
  Comput. 31~(6), 631--649.
\newline\urlprefix\url{http://dx.doi.org/10.1006/jsco.1999.1562}

\bibitem[{Buium and Cassidy(1999)}]{BC}
Buium, A., Cassidy, P.~J., 1999. Differential algebraic geometry and
  differential algebraic groups: from algebraic differential equations to
  diophantine geometry. \cite{Kol3}, 567--636.

\bibitem[{Chen et~al.(2007)Chen, Golubitsky, Lemaire, Moreno~Maza, and
  Pan}]{ctd}
Chen, C., Golubitsky, O., Lemaire, F., Moreno~Maza, M., Pan, W., 2007.
  Comprehensive triangular decomposition. In: Ganzha, V.~G., Mayr, E.~W.,
  Vorozhtsov, E.~V. (Eds.), CASC. Vol. 4770 of Lecture Notes in Computer
  Science. Springer, pp. 73--101.

\bibitem[{Cox et~al.(1992)Cox, Little, and O'Shea}]{CLO}
Cox, D., Little, J., O'Shea, D., 1992. Ideals, varieties, and algorithms.
  Undergraduate Texts in Mathematics. Springer-Verlag, New York, {An
  introduction to computational algebraic geometry and commutative algebra}.

\bibitem[{Dellière(2000)}]{delliere_wang_dyn}
Dellière, S., 2000. {D.M. Wang simple systems and dynamic constructible
  closure}. Rapport de Recherche No. 2000--16 de l'Université de Limoges.
\newline\urlprefix\url{http://www.unilim.fr/laco/rapports/2000/R2000_16.pdf}

\bibitem[{Diop(1992)}]{diop}
Diop, S., 1992. On universal observability. In: Proc. 31st Conference on
  Decision and Control (Tucson, Arizona).

\bibitem[{Gerdt(1999)}]{Ger3}
Gerdt, V.~P., 1999. Completion of linear differential systems to involution.
  In: Computer algebra in scientific computing---{CASC}'99 ({M}unich).
  Springer, Berlin, pp. 115--137.

\bibitem[{Gerdt(2002)}]{Ger2}
Gerdt, V.~P., 2002. On an algorithmic optimization in the computation of
  involutive bases. Programming and Computer Software 28~(2), 62--65.
\newline\urlprefix\url{http://dx.doi.org/10.1023/A:1014816631983}

\bibitem[{Gerdt(2005)}]{GerI}
Gerdt, V.~P., 2005. Involutive algorithms for computing {G}r\"obner bases. In:
  Computational commutative and non-commutative algebraic geometry. Vol. 196 of
  NATO Sci. Ser. III Comput. Syst. Sci. IOS, Amsterdam, pp. 199--225.

\bibitem[{Gerdt(2008)}]{GerdtSimple}
Gerdt, V.~P., 2008. On decomposition of algebraic {PDE} systems into simple
  subsystems. Acta Appl. Math. 101~(1-3), 39--51.
\newline\urlprefix\url{http://dx.doi.org/10.1007/s10440-008-9202-x}

\bibitem[{Gerdt and Blinkov(1998{\natexlab{a}})}]{GB1}
Gerdt, V.~P., Blinkov, Y.~A., 1998{\natexlab{a}}. Involutive bases of
  polynomial ideals. Math. Comput. Simulation 45~(5-6), 519--541,
  {Simplification of systems of algebraic and differential equations with
  applications}.

\bibitem[{Gerdt and Blinkov(1998{\natexlab{b}})}]{GB2}
Gerdt, V.~P., Blinkov, Y.~A., 1998{\natexlab{b}}. Minimal involutive bases.
  Math. Comput. Simulation 45~(5-6), 543--560, {Simplification of systems of
  algebraic and differential equations with applications}.

\bibitem[{Gerdt and Yanovich(2006)}]{GY}
Gerdt, V.~P., Yanovich, D.~A., 2006. Investigation of the effectiveness of
  involutive criteria for computing polynomial {J}anet bases. Programming and
  Computer Software 32~(3), 134--138.
\newline\urlprefix\url{http://dx.doi.org/10.1134/S0361768806030030}

\bibitem[{Gerdt et~al.(2001)Gerdt, Yanovich, and Blinkov}]{GYB}
Gerdt, V.~P., Yanovich, D.~A., Blinkov, Y.~A., 2001. Fast search for the
  {J}anet divisor. Programming and Computer Software 27~(1), 22--24.
\newline\urlprefix\url{http://dx.doi.org/10.1023/A:1007130618376}

\bibitem[{Gómez~Diaz(1994)}]{gomez-diaz_constructible_closure}
Gómez~Diaz, T., 1994. Quelques applications de l'évaluation dynamique. Ph.D.
  thesis, Univ. de Limoges, Limoges, France,
  \url{http://www-igm.univ-mlv.fr/~teresa/dynev/}.

\bibitem[{Habicht(1948)}]{habicht}
Habicht, W., 1948. Eine {V}erallgemeinerung des {S}turmschen
  {W}urzelz\"ahlverfahrens. Comment. Math. Helv. 21, 99--116.

\bibitem[{Hubert(2003{\natexlab{a}})}]{Hubert1}
Hubert, E., 2003{\natexlab{a}}. Notes on triangular sets and
  triangulation-decomposition algorithms. {I}. {P}olynomial systems. In:
  Symbolic and numerical scientific computation ({H}agenberg, 2001). Vol. 2630
  of Lecture Notes in Comput. Sci. Springer, Berlin, pp. 1--39.
\newline\urlprefix\url{http://dx.doi.org/10.1007/3-540-45084-X_1}

\bibitem[{Hubert(2003{\natexlab{b}})}]{Hubert2}
Hubert, E., 2003{\natexlab{b}}. Notes on triangular sets and
  triangulation-decomposition algorithms. {II}. {D}ifferential systems. In:
  Symbolic and numerical scientific computation ({H}agenberg, 2001). Vol. 2630
  of Lecture Notes in Comput. Sci. Springer, Berlin, pp. 40--87.
\newline\urlprefix\url{http://dx.doi.org/10.1007/3-540-45084-X_2}

\bibitem[{Janet(1929)}]{Janet}
Janet, M., 1929. {Le\c cons sur les syst\`emes d'\'equations aux d\'eriv\'ees
  partielles}. Cahiers Scientifiques IV. Gauthiers-Villars, Paris.

\bibitem[{Kolchin(1973)}]{Kol}
Kolchin, E.~R., 1973. Differential algebra and algebraic groups. Academic
  Press, New York, {Pure and Applied Mathematics, Vol. 54}.

\bibitem[{Kolchin(1999)}]{Kol3}
Kolchin, E.~R., 1999. Selected works of {E}llis {K}olchin with commentary.
  American Mathematical Society, Providence, RI, commentaries by Armand Borel,
  Michael F. Singer, Bruno Poizat, Alexandru Buium and Phyllis J. Cassidy,
  edited and with a preface by Hyman Bass, Buium and Cassidy.

\bibitem[{Lemaire et~al.(2005)Lemaire, Moreno~Maza, and Xie}]{regularchains}
Lemaire, F., Moreno~Maza, M., Xie, Y., 2005. {The \textsf{RegularChains}
  library in \textsc{Maple}}. SIGSAM Bull. 39~(3), 96--97.

\bibitem[{Li and Wang(1999)}]{WangLi}
Li, Z., Wang, D., 1999. Coherent, regular and simple systems in zero
  decompositions of partial differential systems. System Science and
  Mathematical Sciences 12, 43--60.

\bibitem[{Mishra(1993)}]{mishra}
Mishra, B., 1993. Algorithmic algebra. Texts and Monographs in Computer
  Science. Springer-Verlag, New York.

\bibitem[{Monagan and Pearce(2009)}]{mope}
Monagan, M., Pearce, R., 2009. Parallel sparse polynomial multiplication using
  heaps. In: {ISSAC '09: Proceedings of the 2009 International Symposium on
  Symbolic and Algebraic Computation}. ACM, New York, NY, USA, pp. 263--270.

\bibitem[{Moreno~Maza(1999)}]{MazaMEGA2000}
Moreno~Maza, M., 1999. On triangular decompositions of algebraic varieties.
  Tech. rep., presented at the MEGA-2000 Conference.

\bibitem[{Plesken(1982)}]{PleskenCountingGroupsAndRings}
Plesken, W., 1982. Counting with groups and rings. J. Reine Angew. Math. 334,
  40--68.
\newline\urlprefix\url{http://dx.doi.org/10.1515/crll.1982.334.40}

\bibitem[{Plesken(2009{\natexlab{a}})}]{PleskenCounting}
Plesken, W., 2009{\natexlab{a}}. Counting solutions of polynomial systems via
  iterated fibrations. Arch. Math. (Basel) 92~(1), 44--56.
\newline\urlprefix\url{http://dx.doi.org/10.1007/s00013-008-2785-7}

\bibitem[{Plesken(2009{\natexlab{b}})}]{PleskenBruhat}
Plesken, W., 2009{\natexlab{b}}. Gauss-{B}ruhat decomposition as an example of
  {T}homas decomposition. Arch. Math. (Basel) 92~(2), 111--118.
\newline\urlprefix\url{http://dx.doi.org/10.1007/s00013-008-2786-6}

\bibitem[{Riquier(1910)}]{Riquier}
Riquier, C., 1910. Les syst\`emes d'\'equations aux d\'eriv\'ees partielles.
  Gauthiers-Villars, Paris.

\bibitem[{Ritt(1950)}]{Ritt}
Ritt, J.~F., 1950. Differential {A}lgebra. American Mathematical Society
  Colloquium Publications, Vol. XXXIII. American Mathematical Society, New
  York, N. Y.

\bibitem[{Rosenfeld(1959)}]{Rosenfeld}
Rosenfeld, A., 1959. Specializations in differential algebra. Trans. Amer.
  Math. Soc. 90, 394--407.

\bibitem[{Seidenberg(1958)}]{Seidenberg58}
Seidenberg, A., 1958. {Abstract differential algebra and the analytic case}.
  Proc. Amer. Math. Soc. 9, 159--164.

\bibitem[{Seidenberg(1969)}]{Seidenberg69}
Seidenberg, A., 1969. {Abstract differential algebra and the analytic case II}.
  Proc. Amer. Math. Soc. 23, 689--691.

\bibitem[{Seiler(2010)}]{Seiler}
Seiler, W.~M., 2010. Involution. Vol.~24 of Algorithms and Computation in
  Mathematics. Springer-Verlag, Berlin, {The formal theory of differential
  equations and its applications in computer algebra}.
\newline\urlprefix\url{http://dx.doi.org/10.1007/978-3-642-01287-7}

\bibitem[{Thomas(1937)}]{Tho1}
Thomas, J.~M., 1937. Differential Systems. AMS Colloquium Publications vol XXI.

\bibitem[{Thomas(1962)}]{Tho2}
Thomas, J.~M., 1962. Systems and Roots. The William Byrd Press, INC, Richmond
  Virginia.

\bibitem[{Wang(1998)}]{wang_simple}
Wang, D., 1998. Decomposing polynomial systems into simple systems. J. Symbolic
  Comput. 25~(3), 295--314.
\newline\urlprefix\url{http://dx.doi.org/10.1006/jsco.1997.0177}

\bibitem[{Wang(2001)}]{WangMethods}
Wang, D., 2001. Elimination methods. Texts and Monographs in Symbolic
  Computation. Springer-Verlag, Vienna.

\bibitem[{Wang(2003)}]{epsilon}
Wang, D., 2003. $\epsilon$\textsf{psilon}: description, help pages and examples
  of use. (\url{http://www-spiral.lip6.fr/~wang/epsilon/}).

\bibitem[{Wang(2004)}]{WangPractice}
Wang, D., 2004. Elimination practice. Imperial College Press, London, software
  tools and applications, With 1 CD-ROM (UNIX/LINUX, Windows).

\bibitem[{Wu(2000)}]{Wu}
Wu, W.-T., 2000. Mathematics mechanization. Vol. 489 of Mathematics and its
  Applications. Kluwer Academic Publishers Group, Dordrecht, mechanical
  geometry theorem-proving, mechanical geometry problem-solving and polynomial
  equations-solving.

\bibitem[{Yap(2000)}]{Yap}
Yap, C.~K., 2000. Fundamental problems of algorithmic algebra. Oxford
  University Press, New York.

\end{thebibliography}

\def\cprime{$'$} \def\cprime{$'$}

\end{document}